\magnification=\magstep1

\let\iteM=\item
\input amstex.tex
\long\def\comment#1\endcomment{}
\def\lable#1{{\rom [#1]}}
\def\nolables{\def\lable##1{}}
\newcount\sectioncount
\newcount\commoncount\commoncount=1
\newcount\lemmacount\lemmacount=1
\newcount\defcount\defcount=1
\newcount\eqcount\eqcount=1
\newcount\theoremcount\theoremcount=1
\newcount\Ccount\Ccount=1

\def\firstsection#1{\sectioncount=#1 \advance\sectioncount by -1}
\firstsection1
\def\renewcount{\commoncount=1\eqcount=1
}
\renewcount

\def\newsection#1#2\par{\global\advance \sectioncount by 1%
\renewcount%
\specialhead\bigbf{\the\sectioncount}.\lable{#1}\bigbf\ #2
\expandafter\xdef\csname Section #1 \endcsname{%
\the\sectioncount}%
\endspecialhead
\relax}
\def\section#1{{\sl Section \csname Section #1 \endcsname}}

\def\newgrphno#1{
\the\sectioncount.\the\commoncount
\expandafter\xdef\csname Grph #1 \endcsname{%
\the\sectioncount.\the\commoncount}%
\global\advance\commoncount by1\relax}
\def\grph#1{{\csname Grph #1 \endcsname}}

\def\newtheorem#1{
\proclaim{\the\sectioncount.\the\commoncount. Theorem}\lable{#1}
\expandafter\xdef\csname Theorem #1 \endcsname{\the\sectioncount%
\the\commoncount}%
\global\advance\commoncount by1\it\relax}
\def\theorem#1{{\sl Theorem \csname Theorem #1 \endcsname}}

%

%

\def\newlemma#1{
\proclaim{\the\sectioncount.\the\commoncount.\lable{#1} Lemma}
\expandafter\xdef\csname Lemma #1 \endcsname{%
\the\sectioncount.\the\commoncount}%
\global\advance\commoncount by1\it\relax}
\def\lemma#1{{\sl Lemma \csname Lemma #1 \endcsname}}
\def\lemmano#1{\csname Lemma #1 \endcsname}

\def\newproposition#1{
\proclaim{\the\sectioncount.\the\commoncount.\lable{#1} Proposition}
\expandafter\xdef\csname Proposition #1 \endcsname{%
\the\sectioncount.\the\commoncount}%
\global\advance\commoncount by1\it\relax}
\def\proposition#1{{\sl Proposition \csname Proposition #1 \endcsname}}
\def\propositionno#1{\csname Proposition #1 \endcsname}

\def\newcorollary#1{
\proclaim{\the\sectioncount.\the\commoncount.\lable{#1} Corollary}
\expandafter\xdef\csname Corollary #1 \endcsname{%
\the\sectioncount.\the\commoncount}%
\global\advance\commoncount by1\it\relax}
\def\corollary#1{{\sl Corollary \csname Corollary #1 \endcsname}}
\def\corollaryno#1{\csname Corollary #1 \endcsname}

\def\newdefinition#1{
{\smallskip\noindent\bf\the\sectioncount.\the\commoncount.\lable{#1}
Definition.}
\expandafter\xdef\csname Definition #1 \endcsname{%
\the\sectioncount.\the\commoncount}%
\global\advance\commoncount by1\relax}
\def\definition#1{{\sl Definition \csname Definition #1 \endcsname}}
\def\definitionno#1{\csname Definition #1 \endcsname}
\let\dfntn=\definition

\def\makeeq#1{
\expandafter\xdef\csname Equation #1 \endcsname{%
\the\sectioncount.\the\eqcount}
\expandafter\xdef\csname Equationlable #1 \endcsname{{#1}}
\global\advance\eqcount by1\relax}
\def\eqqno#1{\csname Equation #1 \endcsname}

\def\eqlable#1{{\rom{[\csname Equationlable #1 \endcsname]}}}

\documentstyle{amsppt}
\let\item=\iteM 
{\catcode`\@=11\global\def\output@{\shipout\vbox{%
 \iffirstpage@ \global\firstpage@false
  \pagebody 
 \else \ifrunheads@ \makeheadline \pagebody
       \else \pagebody \makefootline \fi
 \fi}%
 \advancepageno \ifnum\outputpenalty>-\@MM\else\dosupereject\fi}
}

\def\example#1{{\sl Example \csname Example #1 \endcsname}}
\def\qed{\ \ \hfill\hbox{$\square$}\smallskip}

\let\definition=\dfntn

\def\rom{\ifmmode \fam\rmfam \else\rm \fi}

\input epsf.tex
\newdimen\xsize
\newdimen\oldbaselineskip
\newdimen\oldlineskiplimit
\xsize=.7\hsize
\def\nolineskip{\oldbaselineskip=\baselineskip\baselineskip=0pt%
\oldlineskiplimit=\lineskiplimit\lineskiplimit=0pt}
\def\restorelineskip{\baselineskip=\oldbaselineskip%
\lineskiplimit=\oldlineskiplimit}
\def\putm[#1][#2]#3{
\hbox{\vbox to 0pt{\parindent=0pt%
\vskip#2\xsize\hbox to0pt{\hskip#1\xsize $#3$\hss}\vss}}}%
%
\def\putt[#1][#2]#3{
\vbox to 0pt{\noindent\hskip#1\xsize\lower#2\xsize%
\vtop{\restorelineskip#3}\vss}}

\def\Month{\ifcase\month \or January\or February\or March\or April\or May\or
June\or July\or August\or September\or October\or November\or December\fi}

\font\bigbf=cmbx10 scaled 1200

\thinmuskip = 2mu
\medmuskip = 2.5mu plus 1.5mu minus 2.1mu  
\thickmuskip = 4mu plus 6mu
\def\longpoints{\leaders\hbox to 0.5em{\hss.\hss}\hfill \hskip0pt}
\font\teneusm=eusm10
\font\seveneusm=eusm7
\font\fiveeusm=eusm5
\newfam\eusmfam
\textfont\eusmfam=\teneusm
\scriptfont\eusmfam=\seveneusm
\scriptscriptfont\eusmfam=\fiveeusm
\def\scr#1{{\fam\eusmfam\relax#1}}
\font\tenmib=cmmib10
\font\sevenmib=cmmib7
\font\fivemib=cmmib5
\newfam\mibfam
\textfont\mibfam=\tenmib
\scriptfont\mibfam=\sevenmib
\scriptscriptfont\mibfam=\fivemib
\def\mib{\fam\mibfam}
\font\tensf=cmss10
\font\sevensf=cmss10 scaled 833
\font\fivesf=cmss10 scaled 694
\newfam\sffam
\textfont\sffam=\tensf
\scriptfont\sffam=\sevensf
\scriptscriptfont\sffam=\fivesf
\def\sf{\fam\sffam}
\font\sevensl=cmsl10 scaled 833
\font\fivesl=cmsl10 scaled 694
\scriptfont\slfam=\sevensl
\scriptscriptfont\slfam=\fivesl
\def\sl{\fam\slfam\tensl}
\font\mathnine=cmmi9
\font\rmnine=cmr9
\font\cmsynine=cmsy9
\font\cmexnine=cmex10 scaled 913
\font\ninesl=cmsl9
\font\ninesf=cmss9
\font\ninemib=cmmib9
\def\msmall#1{\hbox{$\displaystyle
\textfont0=\rmnine   \textfont1=\mathnine
\textfont2=\cmsynine \textfont3=\cmexnine
\textfont\slfam=\ninesl \textfont\sffam=\ninesf
\textfont\mibfam=\ninemib
{#1}$}}

\hyphenation{Lip-schit-zian Lip-schitz}
\def\cc{{\Bbb C}}
\def\hh{{\Bbb H}}

\def\ee{{\Bbb E}}
\def\rr{{\Bbb R}}

\def\pp{{\Bbb P}}
\def\ttt{{\Bbb T}} 
\def\zz{{\Bbb Z}}

\def\loc{{\rom loc}}
\def\st{_{\rom st}}

\def\area{{\sf area}\,}

\def\ch{{\sf ch}}

\def\const{{\sf const}}
\def\cotan{{\sf cotan}}
\def\cos{{\sf cos}}
\def\diam{{\sf diam}}
\def\dim{{\sf dim}\,}

\def\dimc{{\sf dim}_\cc}
\def\dist{{\sf dist}}

\def\End{{\sf End}}

\def\exp{{\sf exp}}
\def\ex{{\sf ex}\,}

\def\sfh{{\sf H}}

\def\id{{\sf Id}}

\def\inf{{\sf inf}\,}
\def\im{{\sf Im}\,}

\def\lim{\mathop{\sf lim}}
\def\log{{\sf log}}

\def\osc{{\sf osc}}
\def\pr{{\sf pr}}

\def\re{{\sf Re}\,}

\def\sh{{\sf sh}}

\def\tan{{\sf tan}}
\def\cotan{{\sf cotan}}

\def\eps{\varepsilon}
\def\epsi{\varepsilon}
\let\bs=\bss
\def\ogran{{\hskip0.7pt\vrule height8pt depth4pt\hskip0.7pt}}
\def\comp{\Subset}
\def\d{\partial}
\def\barr#1{\mskip1mu\overline{\mskip-1mu{#1}\mskip-1mu}\mskip1mu}
\def\dbar{{\overline\partial}}

\def\ddef{\mathrel{{=}\raise0.23pt\hbox{\rm:}}}
\def\deff{\mathrel{\raise0.23pt\hbox{\rm:}{=}}}
\def\ge{\geqslant}
\def\longto{\longrightarrow}
\def\inv{^{-1}}
\def\<{\langle}
\def\>{\rangle}
\def\fraction#1/#2{\mathchoice{{\msmall{ #1\over#2}}}%
{{ #1\over #2 }}{{#1/#2}}{{#1/#2}}}
\def\half{{\fraction1/2}}
\def\le{\leqslant}
\def\vph{^{\mathstrut}}
\def\lrar{\longrightarrow}


\def\mapdown#1|#2{\llap{$\vcenter{\hbox{$\scriptstyle #1$}}$}
 {{ \big\downarrow}}
  \rlap{$\vcenter{\hbox{$\scriptstyle #2$}}$}}

\def\emptyset{\varnothing}
\let\hook=\hookrightarrow

\let\wt=\widetilde
\let\ti=\tilde

\def\norm#1{\Vert #1\Vert}
\def\mapo{ \{ \hbox{\it marked}\hskip1ex\allowbreak \hbox{\it points}\} }
\def\scirc{\mathchoice{\mathop{\msmall{\circ}}}{\mathop{\msmall{\circ}}}%
{{\scriptscriptstyle\circ}}{{\scriptscriptstyle\circ}}}
\def\state#1. {\smallskip\noindent{\bf#1. }}
\def\qed{\ \hbox{ }\ \hbox{ }\ {\hfill$\square$}}
\def\Chi{\raise 2pt\hbox{$\chi$}}
\let\phI=\phi\let\phi=\varphi\let\varphi=\phI

%
%
\def\bfbeta{{\mathchar"0C0C}}%
\def\bfgamma{{\mathchar"0C0D}}%
\def\bfvartheta{{\mathchar"0C23}}%
\def\bflambda{{\mathchar"0C15}}%
\def\bfell{{\mathchar"0C60}}%

\def\eg{\hskip1pt plus1pt{\sl{e.g.\/\ \hskip1pt plus1pt}}}
\def\ie{{\sl i.e.\/}\ \hskip1pt plus1pt}
\def\wrt{w.r.t.\ }
\def\iff{\hskip1pt plus1pt{\sl iff\/\hskip2pt plus1pt }}
\def\.{\thinspace}
\def\3{\ss}
\def\isl{\mathchoice{\text{\sl i}}{\text{\sl i}}{\text{\sevensl i}}%
{\text{\fivesl i}}}
\def\sli{{\sl i)} }             
\def\slii{{\sl i$\!$i)} }       
\def\sliii{{\sl i$\!$i$\!$i)} } 
\def\sliv{{\sl i$\!$v)} }       
\def\slv{{\sl v)} }

\def\cal#1{{\scr{#1}}}
\def\cala{{\cal A}}
\def\alg(#1){\cala(\barr{#1})}
\def\calc{{\cal C}}

\def\calm{{\cal M}}

\def\calo{{\cal O}}

\def\cals{{\cal S}}

\def\calu{{\cal U}}
\def\calv{{\cal V}}

\def\bfv{{\mib v}}


\voffset=-.2in
\vsize=580pt
\baselineskip=12.33pt plus .2pt

\nolables

\long\def\comment#1\endcomment{}

\medskip
\centerline{\bigbf GROMOV COMPACTNESS THEOREM}
\medskip
\centerline{\bigbf FOR STABLE CURVES}

\bigskip
\centerline{\bf S. Ivashkovich , V. Shevchishin}
\def\version{\hbox{\font\fiverm=cmr5 \fiverm
Version of \the\day.\the\month.}}
\leftheadtext{\hss\vtop{%
\line{\hfil S.\.Ivashkovich\ \ V.\.Shevchishin \hfil
\llap{\version}}%
\vskip 4pt \hrule }\hss}

\rightheadtext{\hss\vtop{%
\line{\rlap{\version}%
\hfil Compactness for ps.-hol.\ curves\hfil}%
\vskip 4pt \hrule }\hss}


\bigskip
\bigskip
\firstsection0

\newsection{1}{Introduction}

\smallskip
The goal of this paper is to give a proof of the Gromov compactness theorem 
using the language of stable curves in the general situation, setting minimal 
assumption on almost complex structures and on pseudoholomorphic curves. In 
particular, we suppose that the almost complex structures on a target 
manifold are only {\sl continuous} (\ie of class $C^0$) and can vary.

\smallskip
More precisely, we consider a sequence $\{ J_n \}$ of continuous 
almost complex structures on a manifold $X$ which converges 
uniformly to a continuous structure $J_\infty$. Furthermore, let 
$\{ C_n \}$ be a sequence of Riemann surfaces with boundaries of fixed 
topological type. This means that all $C_n$ can be parametrized by the same
real surface $\Sigma$ (see \S\.1 for details). Denote by $\delta_n:\Sigma
\to C_n$ some parametrizations. Finally, let some sequence of
$J_n$-holomorphic maps $u_n:C_n\to X$ be given.

\state Theorem 1. {\it If the areas of $u_n(C_n)$ are uniformly bounded (with
respect to some fixed Riemannian metric on $X$) and the structures $j_{C_n}$ 
of the curves $C_n$ do
not degenerate at the boundary (see {\sl Definition 1.7}), then there exists
a subsequence, still denoted $(C_n, u_n)$, such that

\smallskip
\sl1) \it $C_n$ converge to some nodal curve $C_\infty$ in 
an appropriate completion of the moduli space of Riemann surfaces of given
topological type, \ie there exist a parametrization map $\sigma_\infty:
\Sigma \to C_\infty$ by the same real surface $\Sigma$;

\smallskip
\sl 2) \it one can chose new parametrizations $\sigma_n$ of $C_n$ in such a
way that each $\sigma_n$ coincides with the given parametrization $\delta_n:
\Sigma \to C_n$ outside some fixed compact subset $K\Subset \Sigma$ and the 
structures $\sigma_n^*j_{C_n}$ converge to $\sigma_\infty^*j_{C_\infty}$
in the $C^\infty$-topology on compact subsets outside of the finite set of
circles on $\Sigma$, which are pre-images of the nodal points of $C_\infty$
by $\sigma_\infty$;

{\sl 3)} maps $u_n\scirc \sigma_n$ converge, in the $C^0$-topology on the
whole $\Sigma$ and in the $L^{1,p}_\loc$-topology (for all $p<\infty$)
outside of the pre-images of the nodes of $C_\infty$, to map $u_\infty \scirc
\sigma_\infty$, such that $u_\infty$ is a $J_\infty$-holomorphic map 
$C_\infty\to X$.
}

\smallskip
For the definitions involved and the formal statement we refer to \S\.1 and 
{\sl Theorem 1.1}. Note, that this description of the convergence is precisely
the one given by Gromov in [G]. Our notion of a {\sl stable nodal curve} 
coincides, in fact, with the notion of a {\sl cusp-curve} of Gromov, and 
with the notion of a {\sl stable map} of Kontsevich and Manin [K-M]. The
choice of terminology is explained by the fact that we prefer to
consider our objects as curves rather than maps.

\smallskip
{\sl Theorem 1} generalizes the original result of Gromov in two directions. 
First, we note that the 
Gromov compactness theorem is still valid for continuous and continuously 
varying almost complex structures. This could have an interesting applications,
since now one can consider $C^0$-small perturbations of an almost complex 
structure on a manifold being insured that at least compactness theorem still 
holds true. 

\smallskip
Second, we consider not only the case of closed curves, but also the case 
when $C_n$ are open and of a fixed ``topological type'', so that 
the complex structures of $C_n$ can vary arbitrarily.  In \S\.2 we study 
moduli spaces of open nodal curves. In particular, we define a natural complex
structure for such moduli spaces and show that the condition 
of non-degeneration of complex structures of $C_n$ near boundary is equivalent
to boundedness of $C_n$ in an appropriate completion of the moduli space. 

\smallskip
Let us stop on this point, which is of independent interest. Fix a real 
oriented surface $\Sigma $ of genus $g$ with $m$ marked points and with 
boundary consisting from $b$ circles. Assume that $2g+m+b\ge 3$. Mark 
additionally a point on each boundary circle. Two complex structures $J_1$ 
and $J_2$ on $\Sigma $  are isomorphic if there exists a biholomorphism 
$\phi :(\Sigma ,J_1)\to (\Sigma ,J_2)$ isotopic to identity and preserving 
the marked points. 

Topological space of complex structures modulo this equivalence relation will 
be denoted by $\ttt_\Sigma$, more precise description and notations are 
given in \S\.2. 

For a non-closed Riemann surface $C=(\Sigma ,J)$, we construct the 
holomorphic double of $C$ which is a 
closed Riemann surface $C^d$ containing $C$, and the holomorphic involution 
$\tau$ of $C^d$ interchanching $C$ with $\tau (C)$. Put $D^d=D+\tau (D)$, 
where $D$ is the divisor of marked points, including boundary ones. We prove 
the following

\state
Theorem 2. {\it There is a natural structure of complex manifold on 
$\ttt_\Sigma$ of complex dimension $3g-3+m+2b$ with tangent space 
$T_C\ttt\Sigma$ at $C$ naturally isomorphic to the space  
$\sfh^1(C^d,{\cal O}(TC^d)\otimes {\cal O}(-D^d))^\tau$ of $\tau$-invariants
elements of $\sfh^1(C^d,{\cal O}(TC^d)\otimes {\cal O}(-D^d))$.}

\medskip
Another result of this paper, which we would like to mention in the 
introduction, is an apriori estimate for pseudoholomorphic maps of  ``long 
cylinders'', see {\sl Second Apriori Estimate} in \S\.3. This estimate gives
possibility to treat the degeneration of complex structure on the curves
$C_n$ and the ``bubbling'' phenomenon  in a uniform framework of ``long 
cylinders'' and to get a precise description of the convergence near 
``neck'' singularities where the usual ``strong'' convergence fails. In 
particular, this implies the Hausdorff convergence of the curves $C_n$ 
in {\sl Theorem 1}.

\smallskip
As an application the {\sl Second Apriori Estimate}, we prove in 
{\sl Corollary 3.6} the following generalization of removability theorem 
for the point singularity.

\state Theorem 3. {\it If the area of the image of $J$-holomorphic map
$u:(\check\Delta, J\st)\to (X,J)$ from the punctured disk into a compact almost
complex manifold has ``slow growth'' (``is not growing too fast"), \ie if 
$\area(u(R_k))\le \eps$ for all annuli $R_k \deff \{ z\in \cc :{1\over e^{k+1}
}\le | z| \le {1 \over e^k}\}$ with $k>\!>1$, then $u$ extends to origin.}

\smallskip
The positive constant $\eps $ here depends only on the Hermitian
structure $(J,h)$ of $X$. This theorem under stronger assumption $\sum_k
\area(u(R_k)) \equiv \norm{du}^2_{L^2(\check\Delta)} < \infty$ was proved by
Sacks and Uhlen\-beck [S-U] for harmonic maps, and by Gromov [G] for
$J$-holomorphic maps.
This fact (which is proved here for continuous $J$'s) is new even in the
integrable case. In fact, it measures the ``degree of non-hyperbolicity" (in
the sense of Kobayashi) of $(X,J,h)$.

\medskip
After the ``inner'' case considered in {\sl Theorem 1}, we prove in 
\S\.5 the compactness theorem for curves with boundary on totally real 
submanifolds. For this ``boundary'' case we give appropriate generalizations
of all ``inner'' constructions and estimates. In particular, in 
{\sl Corollary 5.7} we obtain a generalization of the Gromov's result
about removability of boundary point singularity, see [G]. An improvement is
that the statement remains valid also when one has {\sl different} boundary
conditions to the left and to the right from a singular point. Let us explain
this in more details.

\smallskip
Define the (punctured) half-disk by setting $\Delta^+ \deff \{ z\in \Delta :
\im(z)>0\}$ and $\check\Delta^+ \deff \Delta^+\bs \{0\}$. Define $I_- \deff
]-1, 0[ \subset \d\check\Delta^+$ and $I_+ \deff ]0,+1[ \subset \d\check
\Delta^+$. Let a $J$-holomorphic map $u:(\check\Delta^+, J\st)\to (X,J)$ is
given, where $J$ is again continuous. Suppose further that $u(I_+)\subset
W_+$ and $u(I_-) \subset W_-$, where $W_+,W_-$ are totally real submanifolds
of dimension $n= \half \dim_\rr X$ and intersect transversally. Note also, 
that transversality is understand here in a more general sense, see \S\.5
for details.

\smallskip
\state Theorem 4. {\it There is an $\eps^b>0$ such that if for all
half-annuli $R^+_k\deff\{ z\in \Delta^+ :e^{-(k+1)} \le | z| \le e^{-k} \}$
one has $\area(u(R^+_k))\le \eps^b$, then $u$ extends to origin $0\in
\Delta^+$ as an $L^{1,p}$-map for some $p>2$.
}

\smallskip
As in the ``inner'' case, the necessary condition is weaker than the
finiteness of energy. But unlike to ``inner" and smooth boundary cases, it is
possible that the map $u$ in the last statement is $L^{1,p}$-regular in the
neighborhood of ``corner point'' $0\in \Delta^+$ only for some $p>2$. For
example, the map $u(z) = z^\alpha$ with $0<\alpha <1$ satisfies totally real
boundary conditions $u(I_+) \subset \rr$, $u(I_-) \subset e^{\alpha\pi\isl}
\rr$ and is $L^{1,p}$-regular only for $p< p^* \deff {2\over1-\alpha}$. Note 
also, that by Sobolev imbedding $L^{1,p}\subset C^{0,\alpha }$ with $\alpha= 
1 - {2\over p}$, and thus $u$ extends to zero at least continuously. 
In particular, $u(0)\in W_+ \cap W_-$.

One can see such a point $x$ as a {\sl corner point} for a corresponding
pseudoholomorphic curve. Typical example appears in symplectic geometry when
one takes Lagrangian submanifolds as boundary conditions.

The compactness theorem for open stable curves, stated in {\sl Theorem 1},
was essentially used in [I-S1] and [I-S2] to describe envelopes of meromorphy
of 2-spheres in algebraic surfaces.

\smallskip
The organization of the paper is the following. In \S\S\.1 and 2 we present,
for the convenience of the reader, the basis notions concerning the topology
on the space of stable curves and complex structure on the Teichm\"uller
space of Riemann surfaces with boundary. In \S\.3 we give the necessary
apriori estimates for the inner case, and in \S\.4 the prove of {\sl Theorem
1}, related to curves with free boundary. This includes the case of closed
curves. In \S\.5 we consider curves with totally real boundary conditions,
obtain necessary apriori estimates at ``totally real boundary'', and prove
the compactness theorem for such curves. In particular, we prove {\sl Theorem
3} there.

\bigskip
\smallskip
\centerline{\bigbf Table of contents}

{\parindent=0pt
\smallskip
\sl 0. Introduction.\longpoints pp.1--3.

\smallskip
1. Stable curves and Gromov topology.\longpoints pp.4--11.

\smallskip
2. Complex structure on the space $\ttt_\Gamma$.\longpoints pp.12--20.

\smallskip
3. Apriori estimates.\longpoints pp.20--26.

\smallskip
4. Compactness  for the curves with free boundary.\longpoints pp.26--40.

\smallskip
5. Curves with boundary on totally real submanifolds.\longpoints pp.40--63.

\nobreak
\smallskip
References.\longpoints pp.63--64.
}

\bigskip\bigskip
\newsection{1}Stable curves and the Gromov topology

Before stating the Gromov compactness theorem, we need to introduce an
appropriate category of pseudoholomorphic curves. Since in the limit of a
sequence smooth curves one can obtain a singular one, a cusp-curve in the
Gromov's
terminology, we need to allow certain types of singularities of curves. On
the other hand, it is desirable to have singularities as simple as possible.

A similar problem appears in looking for a ``good'' compactification of
moduli spaces $\calm_{g,m}$ of abstract complex smooth closed curves of genus
$g$ with $m$ marked points. The Deligne-Mumford compactification $\barr \calm
_{g, m}$, obtained by adding the {\sl stable curves}, gives a satisfactory
solution of this problem and suggests a possible way of generalization to
other situations. In fact, the only singularity type one should allow are
nodes, or nodal points. An appropriate notion for curves in a complex
algebraic manifold $X$ was introduced by Kontsevich in [K]. Our definition of
stable curves over $(X,J)$ is simply a translation of this notion to almost
complex manifolds. The changee of  terminology from {\sl stable maps over
$(X,J)$} to {\sl stable curves} is motivated by the fact that we want to
consider our objects as curves rather than maps.

Recall that a {\sl standard node} is the complex analytic set $\cala_0 \deff
\{ (z_1,z_2)\in \Delta^2 : z_1\cdot z_2 =0\}$. A point on a complex curve is
called a {\sl nodal point}, if it has a neighborhood biholomorphic to the
standard node.

\state Definition 1.1. {\sl A {\it nodal curve} $C$ is a complex analytic
space of pure dimension 1 with only nodal points as singularities}.

In other terminology, nodal curves are called {\sl prestable}. We shall
always suppose that $C$ is connected and has a "finite topology", \ie $C$ has
finitely many irreducible components, finitely many nodal points, and that
$C$ has a smooth boundary $\d C$ consisting of finitely many smooth circles
$\gamma_i$, such that $\barr C \deff C \cup \d C$ is compact.

\state Definition 1.2. {\sl We say that a real oriented surface with boundary
$(\Sigma, \d\Sigma)$ {\it parameterizes} a complex nodal curve $C$ if there
is a continuous map $\sigma :\barr\Sigma \to \barr C$ such that:

\sli if $a\in C$ is a nodal point, then $\gamma_a = \sigma\inv(a)$ is a
smooth imbedded circle in $\Sigma \bs \d \Sigma $, and if $a\not= b$ then
$\gamma_a \cap \gamma_b= \emptyset$;

\slii $\sigma :\barr\Sigma \bs \bigcup_{i=1}^N\gamma_{a_i}\to \barr C \bs \{
a_1,\ldots ,a_N\} $ is a diffeomorphism, where $a_1,\ldots ,a_N$ are the
nodes of $C$.
}

\bigskip
\vbox{\xsize=.54\hsize\nolineskip\rm
\putm[.13][.01]{\gamma_1}%
\putm[.236][.19]{\gamma_2}%
\putm[.487][0.195]{\gamma_3}%
\putm[.62][.20]{\gamma_4}%
\putm[.84][.15]{\gamma_5}%
\putt[1.1][0]{\advance\hsize-1.1\xsize%
\centerline{Fig.~1}\smallskip
Circles $\gamma_1,..., \gamma_5$ are contracted by the parametrization
map $\sigma$ to nodal points  $a_1, \ldots a_5$.
}%
\putm[.56][.31]{\bigg\downarrow\sigma}%
\noindent
\epsfxsize=\xsize\epsfbox{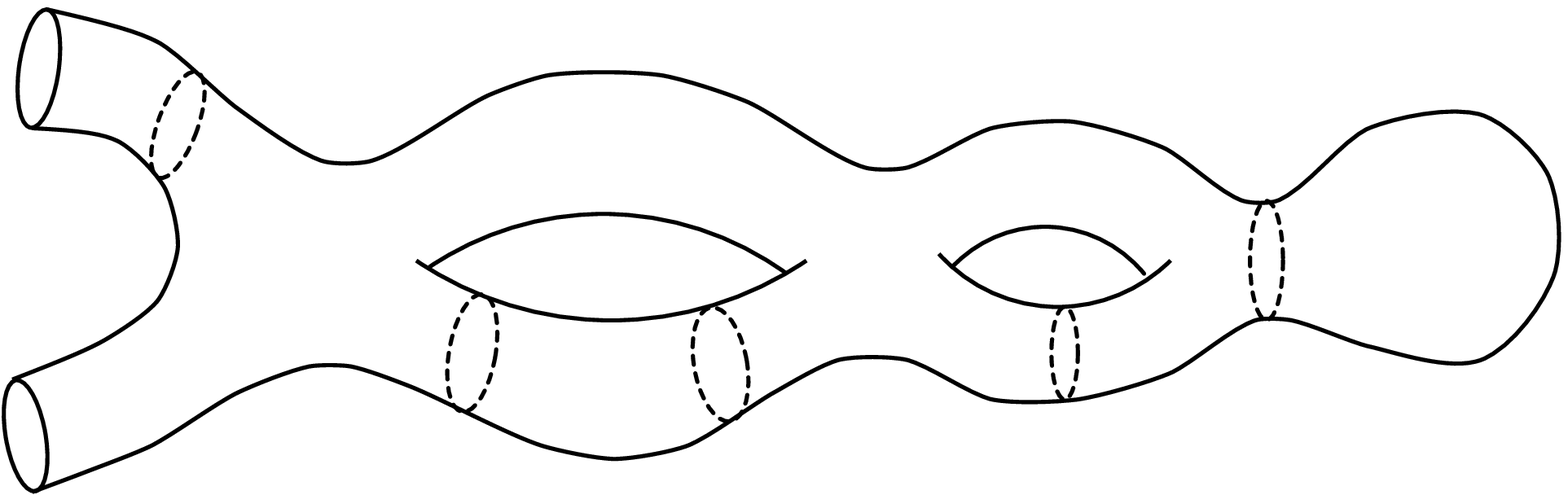}
\vskip15pt
\putm[.14][.045]{a_1}%
\putm[.31][.25]{a_2}%
\putm[.445][0.27]{a_3}%
\putm[.67][.23]{a_4}%
\putm[.805][.20]{a_5}%
\noindent
\epsfxsize=\xsize\epsfbox{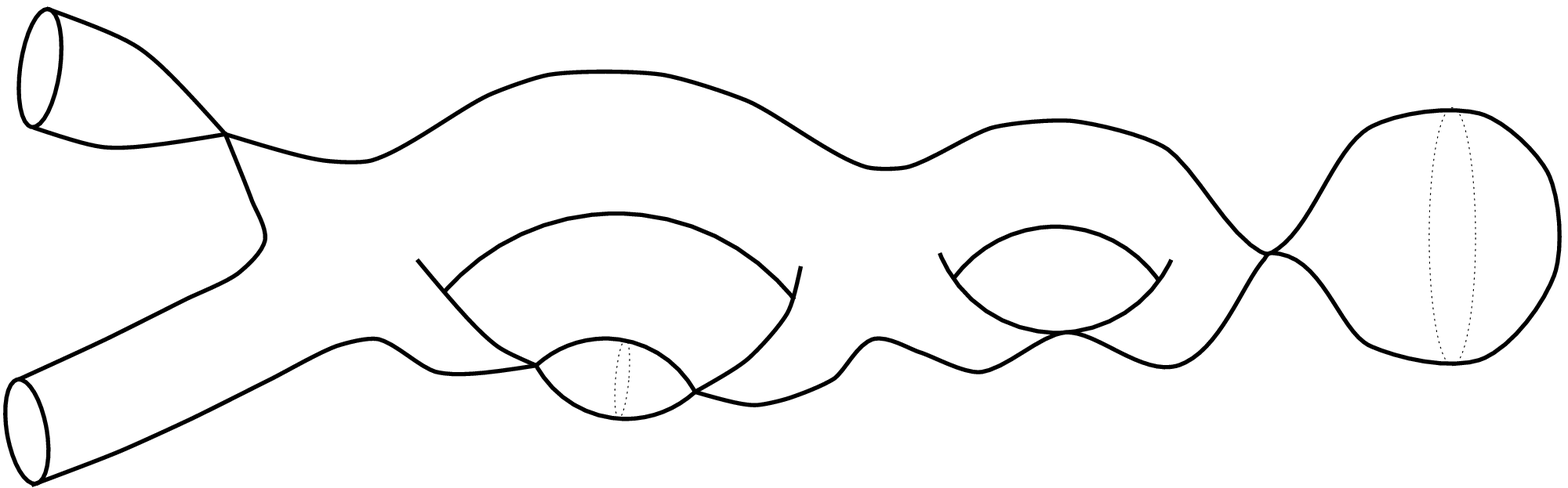}
}

\smallskip
Note that such a parametrization is not unique: if $g:\barr\Sigma \to
\barr\Sigma$ is any orientation preserving diffeomorphism then $\sigma \scirc
g: \barr\Sigma \to \barr C$ is again a parametrization.

\smallskip
A parametrization of a nodal curve $C$ by a real surface can be considered as
a method of ``smoothing'' of $C$. An alternative method of ``smoothing'' ---
the normalization --- is also useful for our purposes.

Consider the normalization $\hat C$ of $C$. Mark on each component of this
normalization the pre-images (under the normalization map $\pi_C: \hat C \to
C$) of nodal points of $C$. Let $\hat C_i$ be a component of $\hat C$. We can
also obtain $\hat C_i$ by taking an appropriate irreducible component $C_i$,
replacing nodes contained in $C_i$ by pairs of disks with marked points, and
marking remaining nodal points. Since it is convenient to consider components
in this form, we make the following

\state Definition 1.3. \it A component $C'$ \sl of a nodal curve $C$ is a
normalization of an irreducible component of $C$ with marked points selected
as above. \rm

\smallskip
This definition allows to introduce the Sobolev and H\"older spaces of
functions and (continuous) maps of nodal curves. For example, a continuous map
$u: C \to X$ is Sobolev $L_\loc ^{1,p}$-smooth if so are all its restrictions
on components of $C$. The  most interesting case is, of course, the one of
continuous $L_\loc ^{1,2}$-smooth maps. In this case the energy functional
$\norm{du}^2_{L^2(C)}$ is defined. The definition of the energy $\norm{du}^2
_{L^2(C)}$ involves Riemannian metrics on $X$ and $C$, which are supposed to
be fixed.

\smallskip
Let $C$ be a nodal curve and $(X,J)$ an almost complex manifold with
continuous almost complex structure $J$.

\nobreak
\state Definition 1.4. {\sl A continuous map $u:C \to X$ is {\it
$J$-holomorphic} if $u\in L^{1,2}_\loc(C,X)$ and
$$
du_x + J\scirc du_x\scirc j_C=0
\eqno(1.1)
$$
for almost all $x\in C$. Here $j_C$ denotes the complex structure on $C$.

The {\it area} of $J$-holomorphic map is defined as
$$
\area(u(C)) \deff \norm{du}^2_{L^2(C)}.
$$
}

\smallskip
We shall show later that every $J$-holomorphic $u$ is, in fact,
$L^{1,p}_\loc(C,X)$-smooth  for all $p< \infty$, see discussion after
{\sl Lemma 3.1}.

\state Remark. Our definition of the area uses the following fact. Let $g$ be
a Riemannian
metric on $C$ compatible with $j_C$, $h$ a Riemannian metric on $X$, and $u:C
\to X$ a $J$-holomorphic immersion. Then $\norm{du}^2_{L^2(C)}$ is independent
of the choice of $g$ and coincides with the area of the image $u(C)$ \wrt
the metric $h_J(\cdot, \cdot) \deff \half(h(\cdot, \cdot) + h(J\cdot, J\cdot)
)$. The metric $h_J$ here can be seen as a ``Hermitization'' of $h$ \wrt
$J$. The independence of $\norm{du}^2_{L^2(C)}$ of the choice of a metric $g$
on $C$ in the same conformal class is a well-known fact, see \eg [S-U]. Thus
we can use the flat metric $dx^2 + dy^2$ to compare the area and the energy.
For $J$-holomorphic map we get
$$
\norm{du}^2_{L^2(C)}= \int_C |\d_x u|_h^2 + |\d_y u|_h^2 =
\int_C |\d_x u|_h^2 + |J \d_x u|_h^2 = \int_C |du|_{h_J}^2 = \area_{h_J}(u(C)),
$$
where the last equality is another well-known result, see \eg [G]. Since we
consider changing almost complex structures on $X$, it is useful to know that
we can use any Riemannian metric on $X$ having reasonable notion of area.

\smallskip
\state Definition 1.5. {\sl A {\it stable curve over $(X,J)$} is a pair
$(C,u)$, where $C$ is a nodal curve and $u:C\to X$ is a $J$-holomorphic map,
satisfying the following condition: If $C'$ is a compact component of $C$, 
such that
$u$ is constant on $C'$, then there exist finitely many biholomorphisms
of $C'$ which preserve the marked points of $C$.
}

\state Remark. One can see that stability condition is nontrivial only in the
following cases:

\item{\sl 1)} some component $C'$ is biholomorphic to $\cc\pp^1$ with 1 or 2
marked points; in this case $u$ should be non-constant on any such component
$C'$;

\item{\sl 2)} some irreducible component $C'$ is $\cc\pp^1$ or a torus
without marked points.

\noindent
Since we consider only connected nodal curves, case {\sl 2)} can happen only
if $C$ irreducible, so that $C'=C$. In this case $u$ must be non-constant on
$C$.

\smallskip\rm
Now we are going to describe the Gromov topology on the space of stable
curves over $X$ introduced in [G]. Let a sequence $J_n$ of continuous almost
complex structures on $X$ be given, and suppose that  $J_n$ 
 converge to $J_\infty$ in the $C^0$-topology. Furthermore, let $(C_n,
u_n)$ be a sequence of stable curves over $(X, J_n)$, such that all $C_n$ are
parametrized by the same real surface $\Sigma$.

\state Definition 1.6. {\sl We say that $(C_n,u_n)$ {\it converges to a
stable $J_\infty$-holomorphic curve $(C_\infty,u_\infty)$ over $X$} if the
parametrizations $\sigma_n: \barr\Sigma \to \barr C_n$ and $\sigma_\infty:
\barr \Sigma \to \barr C_\infty$ can be chosen in such a way that the
following holds:

\sli $u_n\scirc \sigma_n$ converges to $u_\infty\scirc \sigma_\infty$ in the
$C^0( \Sigma, X)$-topology;

\slii if $\{ a_k \}$ is the set of nodes of $C_\infty$ and $\{\gamma_k\}$ are
the corresponding circles in $\Sigma$, then on any compact subset $K\comp
\Sigma \bs \cup_k \gamma_k$ the convergence $u_n\scirc \sigma_n \to u_\infty
\scirc \sigma_\infty$ is $L^{1,p}(K, X)$ for all $p< \infty$;

\sliii for any compact subset $K\comp \barr\Sigma \bs \cup_k\gamma_k$ there
exists $n_0=n_0(K)$ such that $ \sigma_n^{-1}(\{ a_k \}) \cap K= \emptyset$
for all $n\ge n_0$ and the complex structures $\sigma_n^*j_{C_n}$ converge
smoothly to $\sigma_0^*j_{C_0}$ on $K$;

\sliv the structures $\sigma_n^*j_{C_n}$ are constant in $n$ near the
boundary $\d\Sigma$.
}

\medskip
The reason for introducing the notion of a curve stable over $X$ is similar
to the one for the Gromov topology. We are looking for a completion of the
space of smooth imbedded pseudoholomorphic curves which has ``nice''
properties, namely: \.1) such a completion should contain the limit of a
subsequence of every sequence of smooth curves, bounded in an appropriate
sense; \.2) such a limit should exist also for a subsequence of every sequence 
in the
completed space; \.3) such a limit should be unique. The Gromov's compactness
theorem insures us that the space of curves stable over $X$ enjoys these nice
properties.

Condition \sliv is trivial if $\Sigma $ is closed, but it is important when
one considers the ``free boundary case'', \ie when $\Sigma$ (and thus all
$C_n$) are not closed and no boundary condition is imposed. However, we
would like to point out that in our approach the ``free boundary case'' is
essentially involved in the proof of compactness theorem also in the case
of closed curves. On the other hand, in the case of curves with boundary on
totally real submanifolds (see \S\.5) such a condition is unnecessary.

\smallskip
Recall that a complex annulus $A$ has a conformal radius $R>1$ if $A$ is
biholomorphic to $A(1,R) \deff \{ z\in \cc \,:\, 1 <|z| < R \}$.

\state Definition 1.7. {\sl Let $C_n$ be a sequence of nodal curves,
parametrized by the same real surface $\Sigma$. We say that the complex
structures on $C_n$ {\it do not degenerate near boundary}, if there exist
$R>1$, such that for any $n$ and any boundary circle $\gamma_{n, i}$ of $C_n$
there exist an annulus $A_{n,i} \subset C_n$ adjacent to $\gamma_{n, i}$,
such that all $A_{n,i}$ are mutually disjoint, do not contain nodal points of
$C_n$, and have the same conformal radius $R$.}

Since conformal radius of all $A_{n, i}$ is the same, we can identify them
with $A(1,R)$. This means that all changes of complex structures of $C_n$
take place away from boundary. The condition is trivial if $C_n$ and $\Sigma$
are closed, $\d\Sigma = \d C_n = \emptyset$.

\state Remark. Changing our parametrizations $\sigma_n: \Sigma \to C_n$, we
can suppose that for any $i$ the pre-image $\sigma_n\inv (A_{n,i} )$ is the
same annulus $A_i$ independent of $n$.

\medskip
Now we state our main result. Fix some Riemannian metric $h$ on $X$ and some
$h$-complete set $A\subset X$.

\state Theorem 1.1. {\it Let $\{(C_n,u_n)\}$ be a sequence of stable
$J_n$-holomorphic curves over $X$ with parametrizations $\delta_n: \Sigma \to
C_n$. Suppose that:

\item{\sl a)} $J_n$ are continuous almost complex structures on $X$,
$h$-uniformly converging to $J_\infty$ on $A$ and $u_n(C_n)\subset A$ for all
$n$;

\item{\sl b)} there is a  constant $M$ such that $\area [u_n (C_n)]\le M$
for all $n$;

\item{\sl c)} complex structures on $C_n$ do not degenerate near the boundary.

Then there is a subsequence $(C_{n_k},u_{n_k})$ and parametrizations $\sigma
_{n_k}: \Sigma \to C_{n_k}$, such that $(C_{n_k}, u_{n_k}, \sigma_{n_k})$
converges to a stable $J_\infty$-ho\-lo\-mor\-phic curve $(C_\infty, u_\infty,
\sigma_\infty)$ over $X$.

Moreover, if the structures $\delta_n^*j_{C_n}$ are constant on the fixed
annuli $A_i$, each adjacent to a boundary circle $\gamma_i$ of $\Sigma$, then
the new parametrizations $\sigma_{n_k}$ can be taken equal to $\delta_{n_k}$
on some subannuli $A'_i \subset A_i$, also adjacent to $\gamma_i$.
}

\state Remarks.~1. In the proof, we shall give a precise description of
convergence with estimates in neighborhoods of the contracted circles
$\gamma_i$. The convergence of curves with boundary on totally real
submanifolds will be studied in \S5.

\state 2. In applications, one uses a generalized version of the Gromov
compactness theorem for nodal curves with marked point. This version is an
immediate consequence of {\sl Theorem 1.1} due to the following construction.
Consider a nodal curve $C$ and a $J$-holomorphic map $u: C \to X$. Let $\bold
x \deff \{x_1,\ldots,x_m\}$ be the set of marked points on $C$ which are
supposed to be distinct from the nodal points of $C$. Define a new curve
$C^+$ as the union of $C$ with disks $\Delta_1,\ldots, \Delta_m$ such that $C
\cap \Delta_i = \{x_i\}$ and such that any $x_i$ becomes a nodal point of
$C^+$. Extend $f$ to a map $f^+: C^+ \to X$ by setting $f^+\ogran_{\Delta_i}$
to be constant and equal to $f(x_i)$. An appropriate definition of
stability, used for triples $(C, \bold x, f)$, is equivalent to stability of
$(C^+, f^+)$. Similarly, the Gromov convergence
$(C_n, \bold x_n, f_n) \to (C_\infty, \bold x_\infty, f_\infty)$
is equivalent to the Gromov convergence
$(C^+_n, f^+_n) \to (C^+ _\infty, f^+_\infty)$. Thus the Gromov compactness for
curves with marked points reduces to the case considered in our paper.
However, we shall consider curves with marked points as well.

\medskip
In the rest of this section we shall describe topology and conformal geometry
of nodal curves and compute the set of moduli parameterizing deformations of
complex structure. As a basic reference we use the book of Abikoff [Ab].

\smallskip
Let $C$ be a complex nodal curve parametrized by $\Sigma$.

\state Definition 1.8. {\sl A component $C'$ of $C$ is called {\it nonstable}
if one of the following two cases occurs:

\item{ 1)} $C'$ is $\cc\pp^1$ and has one or two marked points;

\item{ 2)} $C'$ is $\cc\pp^1$ or a torus and has no marked points.
}

\smallskip
This notion of stability of abstract closed curves is due to Deligne-Mumford,
see [D-M]. It was generalized by Kontsevich [K?] for the case of maps $f:C
\to X$, \ie for curves over $X$ in our terminology. As it was already noted,
the last case can happen only if $C=C'$. Strictly speaking, this case should
be considered separately. However since such considerations require only
obvious changes we just skip them and suppose that case {\sl 2)} does not
occur.

Our first aim is to analyse the behavior of complex structures in the sequence
$(C_n,u_n)$ of $J_n$-holomorphic curves stable over $X$ with uniformly bounded
area, which are parametrized by the same real surface $\Sigma$. At the moment,
uniform bound of area of $u_n(C_n)$ is needed only to show that the number of
components of $C_n$ is bounded. Passing to a subsequence, we can assume that
all $C_n$ are homeomorphic. This reduces the problem to a description of
complex structures on a fixed nodal curve $C$.

To obtain such a description, it is useful to cut the curve into pieces where
the behavior of complex structure is easy to understand. Such a procedure is
a {\sl partition into pants}. It is well known in the theory of moduli spaces
of complex structure on curves, see \eg [Ab], p.\.93. Making use of it, we
shall also do a slightly different procedure. Namely, we shall choose a
special covering of $\Sigma$ instead of its partition. Further, as blocks for
our construction we shall use not only pants, but also disks and annuli.
The reasons are that, firstly, the considered curves can have unstable
components and, second, it is convenient to use annuli for a description of
the deformation of the complex structure on curves. We start with

\state Definition 1.9. \it An annulus $A$ \sl on a real surface or on a
complex curve is a domain which is diffeomorphic (resp.\ biholomorphic) to
the standard annulus $A(r,R) \deff \{ z\in \cc \;:\; r<|z|<R \}$, such that
its boundary consists of smoothly imbedded circles. \it Pants \sl(also
called a \it pair of pants\sl) on a real surface or on a complex curve is a
domain which is diffeomorphic to a disk with 2 holes. \rm

The boundary of pants consists of three components, each of them being
either a smoothly imbedded circle or a point. This point can be considered as
a puncture of pants or as a marked point. An annulus or pants is {\sl
adjacent to a circle $\gamma$} if $\gamma$ is one of its boundary components.

\medskip
\vbox{
 \hbox{%
   \vtop{
\hsize=.47\hsize\noindent\epsfxsize=\hsize\epsfbox{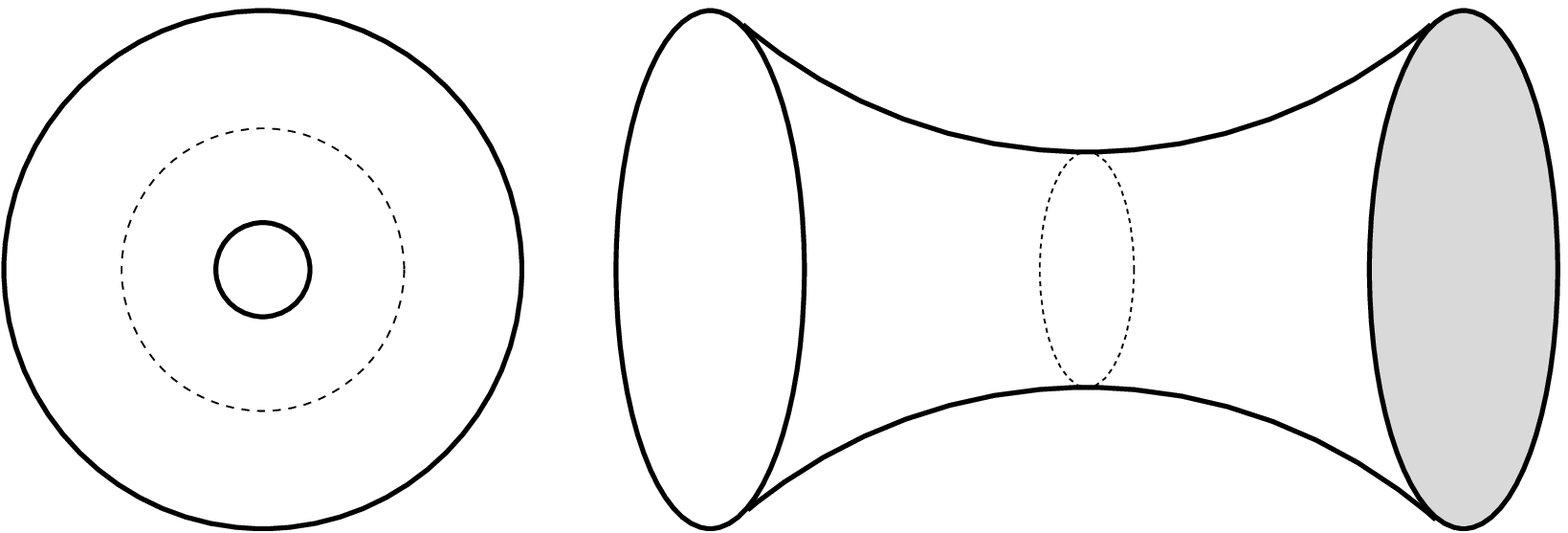}%
\smallskip
\centerline{Fig.~2. An annulus}
\smallskip
It is useful to imagine an annulus as a cylinder. After contracting the
middle circle of the annulus we get a node.
}
\hskip.05\hsize
   \vtop{
\hsize=.48\hsize\noindent\epsfxsize=\hsize\epsfbox{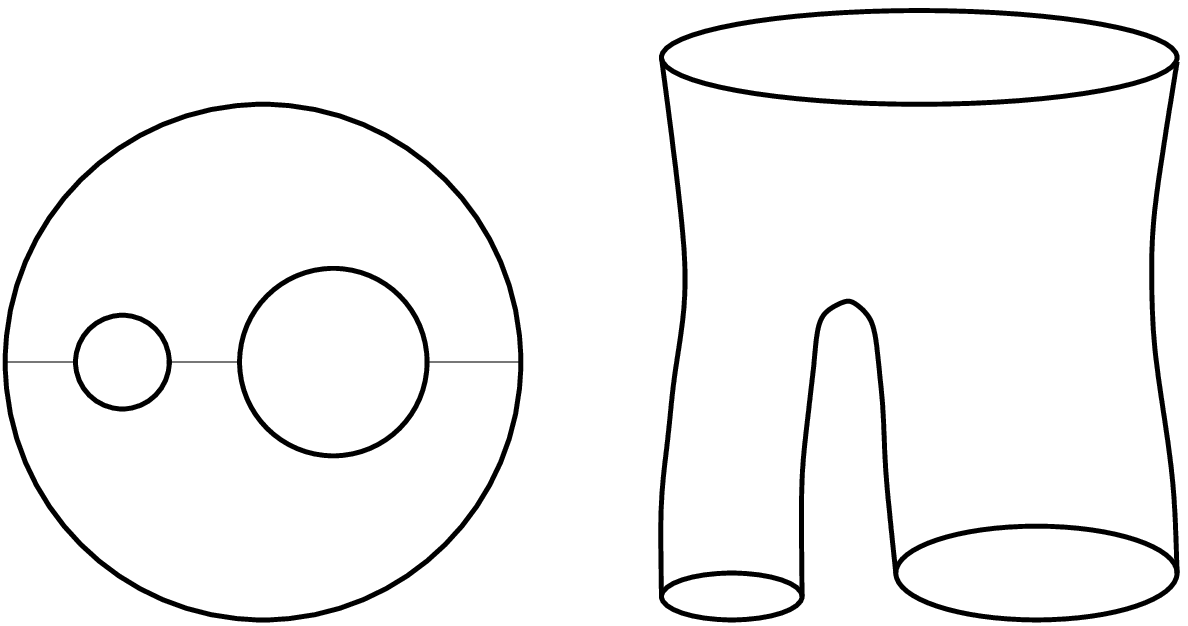}
\smallskip
\centerline{Fig.~3. Pants.}
\smallskip
One can consider pants as a disc with two holes or as a
sphere with three holes. }
}}

\medskip
Let $C$ be a nodal curve parametrized by a real surface $\Sigma$. We shall
associate with every such curve $C$ a certain graph $\Gamma_C$, which
determines $C$ topologically in a unique way. In fact, $\Gamma_C$ will also
determine a decomposition of some components of $C$ into pants.

\smallskip
By the definition, a compact component $C'$ is stable if it contains only a
finite number of automorphisms preserving marking points. In this case $C'
\bs\mapo$ possesses a unique so-called intrinsic metric.

\state Definition 1.10. The {\sl intrinsic metric} for a smooth curve $C$
with marked points $\{ x_i \}$ and with boundary $\d C$ is a metric $g$ on
$C\bs \mapo$ satisfying the following properties:

\sli $g$ induces the given complex structure $j_C$;

\slii the Gauss curvature of $g$ is constantly -1;

\sliii $g$ is complete in a neighborhood of every marked point $x_i$;

\sliv every boundary circle $\gamma$ of $C$ is geodesic \wrt $g$.

\smallskip
Note that such a metric, if exists, is unique, see \eg [Ab].

\smallskip
Now consider a component $C'$ of $C$ adjacent to some boundary circle of $C$.
Then $C' \bs \mapo$ is either

a)~a disk $\Delta$, or

b)~an annulus $A$, or

c)~a punctured disk $\check\Delta$, or else

d) $C'\bs \mapo$ admits the intrinsic metric.

\noindent
Note that if a component $C'$ is a disk or an annulus (both without marked
points), then $C'$ is the whole curve $C$. We shall considered cases a) and
b) later. Now we assume for simplicity that cases a) and b) do not occur.

\state Definition 1.11. \sl A component $C'$ of a nodal curve $C$ is {\it
non-exceptional} if $C' \bs \mapo$ admits the intrinsic metric. \rm

\smallskip
In particular, nonstable components are exceptional compact ones, and
exceptional non-compact components are those of following types a)--c) above.

\smallskip
Take some non-exceptional component $C'$ of $C$. There is a so called maximal
partition of $C'\bs \mapo$ into pants $\{ C_1,\ldots C_n\}$, such that all
boundary components of these pants are either simple geodesics circles in
intrinsic metric or marked points, see [Ab]. Let us fix such a partition and
mark the obtained geodesics circles on $C'$.

\medskip
Let now $\sigma: \Sigma \to C$ be some parametrization of $C$. This defines
the set $\bfgamma'$ of the circles on $\Sigma$ which correspond to the nodes
of $C$. Let $\bfgamma''$ be the set of $\sigma $-preimages of the geodesics,
chosen above. Then $\bfgamma\deff \bfgamma' \sqcup \bfgamma''$ forms a system
of disjoint ``marked'' circles on $\Sigma$, which encodes the topological
structure of $C$. Now the graph $\Gamma_C$ in question can be constructed as
follows.

Define the set $V_C$ of vertices of $\Gamma_C$ to be the set $\{ S_j \}$ of
connected components of $\Sigma \bs \cup_{\gamma \in \bfgamma} \gamma =
\sqcup_j S_j$. Any $\gamma \in \bfgamma$ lies between 2 components, say $S_j$
and $S_k$, and we draw an edge connecting the corresponding 2 vertices.
Further, any boundary circle $\gamma$ of $\Sigma$ has the uniquely defined
component $S_j$ adjacent to $\gamma$. For any such $\gamma$ we draw a~{\sl
tail}, \ie an edge with one end free, attached to vertex $S_j$. Finally, we
mark all edges which correspond to the circles $\bfgamma'$, \ie those coming
from from nodes.

\medskip\medskip
\vbox{\xsize=.5\hsize\nolineskip\rm
\putm[.20][-.01]{\gamma_1^*}%
\putm[.46][.13]{\gamma_2^*}%
\putm[.67][.12]{\gamma_3^*}%
\putm[.21][.58]{a_1^*}%
\putm[.38][.815]{a_2^*}%
\putm[.53][.83]{a_3^*}%
\putt[1.05][0]{\advance \hsize -1.05\xsize
\centerline{Fig.~4. Graph of a curve $C$.}
\smallskip
Graph $\Gamma_C$ determines the topology of the curve $C$ in a unique way.
Take as many oriented spheres as many vertices $\Gamma_C$ has. For each edge
take a handle and join the corresponding spheres by this handle. For each
tail make a hole (\ie remove a disk) in the corresponding sphere. Finally,
contract into points the circles on the handles corresponding to the marked
edges to get nodes. We obtain a topological space homeomorphic to $C$.
\smallskip
}
\noindent
\epsfxsize=\xsize\epsfbox{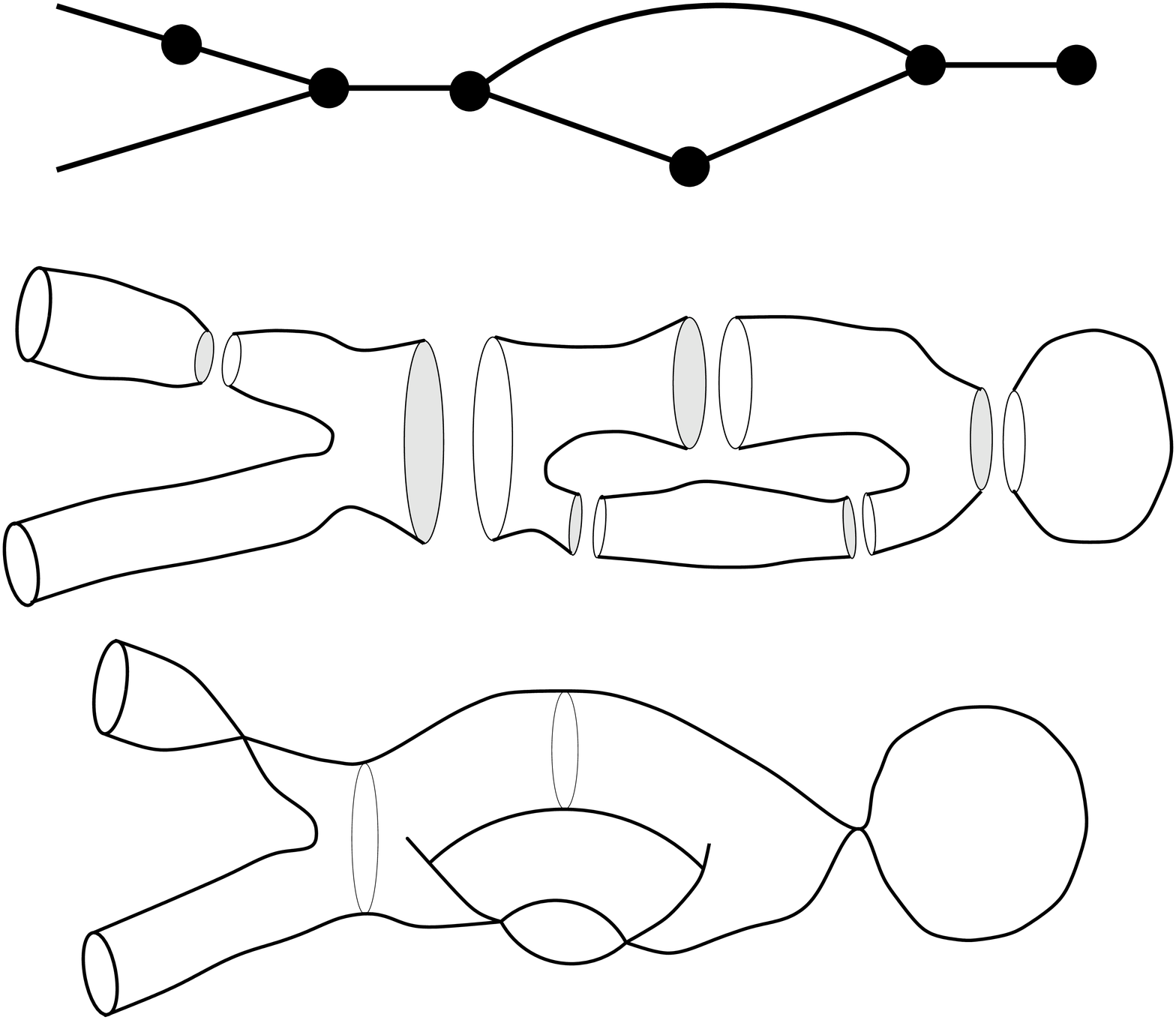}
}

\medskip
\bgroup\baselineskip=13.8pt
Having the graph $\Gamma$, which characterizes uniquely the topological
structure of $C$, we are now going to describe the set of parameters,
defining (uniquely) the complex structure of the curves $C$. This is
equivalent to determining the complex structure and marked points on all
components of $C$. If such a component $C'$ is a sphere with 1 or 2 marked
points or a disk with 1 marked point, then its structure is defined by its
topology uniquely up to diffeomorphism. Otherwise, the component $C'$ is
non-exceptional. In this case the complex structure and the marked points can
be restored by the so called {\sl Fenchel-Nielsen coordinates} on the
Teichm\"uller space $\ttt_{g,m,b}$. Recall that the space $\ttt_{g,m,b}$
parametrizes the complex structures on a Riemann surface $\Sigma$ of genus
$g$ with $m$ punctures (\ie marked points) and with boundary consisting of
$b$ circles, see [Ab].

\smallskip
Let $C$ be a smooth complex curve with marked points of non-exceptional type,
so that $C$ admits the intrinsic metric. Fix some parametrization $\sigma:
\Sigma \to C$. Consider the preimages of the marked points on $C$ as marked
points on $\Sigma$ or, equivalently, as punctures of $\Sigma$. Let $C \bs
\mapo = \cup_j C_j$ be a decomposition of $C$ into pants and $\Sigma\bs
\mapo = \cup_j S_j$ the induced decomposition of $\Sigma$.

Let $\{\gamma_i\}$ be the set of boundary circles of $\Sigma$. The boundary
of every pants $S_j$ has three components, each of them being either a marked
point of $\Sigma$ or a circle. In the last case this circle is either a
boundary component of $\Sigma$ or a boundary component of another pants, say
$S_k$. In this situation we denote by $\gamma _{jk}$ the circle lying between
the pants $S_j$ and $S_k$. Fix the orientation on $\gamma _{jk}$, induced
from $S_j$ if $j<k$ and from $S_k$ if $k<j$. For any such circle $\gamma
_{jk}$, fix a boundary component of $S_j$ different from $\gamma_{jk}$ and
denote it by $\d_k S_j$. In the same way fix a boundary component $\d_j S_k$.
Make similar notations on $C$ using primes to distinguish the circles on $C$
from those on $\Sigma$, \ie set $\gamma'_i \deff \sigma(\gamma_i)$ and
$\gamma'_{jk} \deff \sigma(\gamma_{jk})$.

By our construction, $\gamma'_{jk} = \sigma( \gamma_{jk})$ is a geodesic
\wrt intrinsic metric in $C$.

\egroup

\smallskip
\vbox{\xsize=.43\hsize\nolineskip\rm
\putm[.42][.12]{\gamma'_{jk}}%
\putm[.51][.21]{x^*_{j,k}}%
\putm[.33][.31]{x^*_{k,j}}%
\putm[.26][.24]{S_j}%
\putm[.70][.36]{S_k}%
\putm[-.05][.12]{\d_k S_j}%
\putm[.86][.29]{\d_j S_k}%
\putt[1.05][-.03]{\advance\hsize-1.05\xsize
\noindent
If the component $\d_k C_j$ is a marked point, we find on $C_j$ the (uniquely
defined) geodesic ray $\alpha_{j,k}$ starting at some point $x^*_{j,k} \in
\gamma'_{jk}$ and approaching $\d_k C_j$ at infinity, such that $\alpha
_{j,k}$ has no self-intersections and is orthogonal to $\gamma'_{jk}$ at $x^*
_{j,k}$. Otherwise, we find on $C_j$ the shortest geodesic $\alpha _{j,k}$
which connects $\d_k C_j$ with $\gamma'_{jk}$ and denote the point $\alpha
_{j,k} \cap \gamma'_{jk}$ by $x^*_{j,k}$. In both cases, this construction
detemies a distinguished point $x^* _{j,k} \in \gamma' _{jk}$.
}
\noindent
\epsfxsize=\xsize\epsfbox{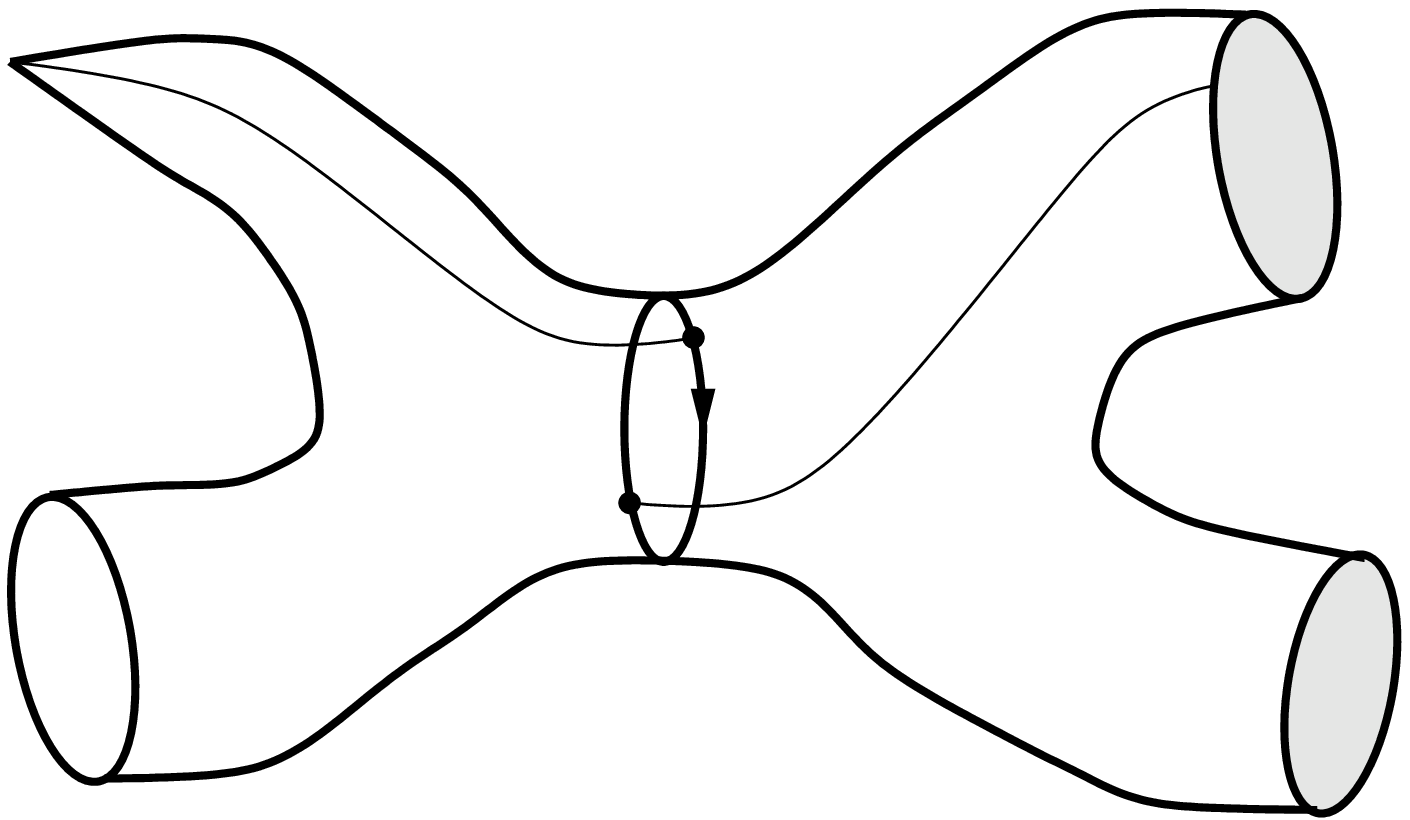}
\medskip
\vbox{\hsize=\xsize\centerline{Fig.~5. Marked points}
\centerline{on the circle $\gamma'_{jk}$.}
}
\smallskip
}

\medskip
Doing the same procedure in $C_k$, we obtain another point $x^*_{k,j} \in
\gamma' _{jk}$. Denote by $\ell_{jk}$ (resp.\ by $\ell_i$) the intrinsic
length of $\gamma'_{jk}$ (resp.\ of $\gamma'_i \deff \sigma( \gamma_i)$) in
$C$. For $j<k$ define $\lambda_{jk}$ as the intrinsic length of the arc on
$\gamma' _{jk}$, which starts at $x^* _{j,k}$ and goes to $x^* _{k,j}$ in the
direction determined by the orientation of $\gamma _{jk}$. Set $\vartheta
_{jk} \deff {2\pi \lambda _{jk} \over \ell_{jk} }$. We shall consider
$\vartheta _{jk}$ as a function of the complex structure $j_C$ on $C$ with
values in $S^1 \cong \rr/ 2\pi \zz$.

The parameters $\bfell \deff (\ell_i, \ell_{jk} )$ and
$\bfvartheta \deff(\vartheta_{jk})$ are called
{\sl Fenchel-Nielsen coordinates} of the complex
structure $j_C$. The reason is that these parameters determines up to
isomorphism the complex structure $j_C$ on the smooth complex curve with
marked points parametrized by a real surface $\Sigma$. In other words,
$(\bfell, \bfvartheta)$ can be considered as coordinates on $\ttt_{g,m,b}$.
More precisely, one has the following

\state Proposition 1.2. \it Let $\Sigma$ be a real surface of genus $g$ with
$m$ marked points and with the boundary consisting of $b$ circles, so that
$2g+m +b\ge3$. Let $\Sigma \bs \mapo = \cup_j S_j$ be its decomposition into
pants. Then

\item{\sl i)} for any given tuples $\bfell= (\ell_i, \ell_{jk})$ and
$\bfvartheta = (\vartheta_{jk})$ with $\ell_i, \ell_{jk} >0$ and $\vartheta_{jk} \in S^1$
there exists a complex structure $j_C$ on $\Sigma$, such that boundary circles
of all $S_j$ are geodesic \wrt the intrinsic metric on $\Sigma\bs \mapo$
defined by $j_C$, and such that the given $(\bfell, \bfvartheta)$
are Fenchel-Nielsen coordinates of $j$; moreover, such a structure $j_C$
is unique up to a diffeomorphism preserving the pants $S_j$ and the marked
points;

\item{\sl ii)} let $C$ be a smooth complex curve with parametrization $\sigma:
\Sigma \to C$ which has $m$ marked points; then there exists a parametrization
$\sigma_1: \Sigma \to C$ isotopic to $\sigma$, which maps boundary components
and marked points of $\Sigma$ onto the ones of $C$ in prescribed order, and
such that the boundary circles of $\sigma(S_j)$ are geodesic \wrt the
intrinsic metric on $C\bs \mapo$. \rm

\nobreak
\state Proof. See [Ab]. \qed

\bigskip\bigskip
\newsection{2}{Complex structure on the space $\ttt_\Gamma$}

Let $\Sigma$ be a real oriented surface of genus $g$ with $m$ marked points 
and with
the boundary consisting of $b$ circles. Assume that $2g+m +b\ge3$. Then there
exists a decomposition of $\Sigma \bs \mapo$ into pants, which is in general
not unique. The topological type of such a decomposition can be encoded in
graph $\Gamma$, associated with the decomposition. It is constructed in
similar way as above, but this time we must draw a tail also for every marked
point, and then mark all those tails on the graph.

Let such a graph $\Gamma$ be fixed. We call two complex structures $J_1$ and
$J_2$ on $\Sigma$ isomorphic if there exists a biholomorphism $\phi: (\Sigma,
J_1) \cong (\Sigma, J_1)$ preserving the marked points of $\Sigma$ and the
decomposition of $\Sigma$ into pants given by graph $\Gamma$. Denote by
$\ttt' _\Gamma$ the space of isomorphism classes of complex structures on
$\Sigma$. By {\sl Proposition 1.2}, Fenchel-Nielsen coordinates identify
$\ttt'_\Gamma$ with the real manifold
$\rr_+^{3g-3+m+2b} \times (S^1)^{3g-3+m+b}$.

It is desirable to equip $\ttt'_\Gamma$ with some natural complex structure.
In doing so, the main difficulty is that the real dimension of $\ttt'_\Gamma$
can be odd. A possible explanation of this fact is that not all relevant
information (\ie parameters) about a complex structure has been taken into
consideration. Note that for any ``inner circle'' $\gamma_{jk}$ which appears
after the decomposition into pants we have obtained a pair of coordinates,
manely the length $\ell_{jk}$ and the angle $\vartheta_{jk}$. On the other
hand, for any boundary circle $\gamma_i$ of $\Sigma$ we have got only the
length $\ell_i$. An obvious way to produce additional angle coordinates is to
introduce an additional marking of every boundary circle.

\state Definition 2.1. \sl A real surface $\Sigma$ or a nodal complex curve
$C$ is said to have a {\it marked boundary} if on every boundary circle of
$\Sigma$ (resp.\ $C$) a point is fixed. \rm

\state Remark. Later in \S\.5 we shall consider complex curves with several
marked points on boundary circles. But now we shall assume that on every
boundary circle exactly one point is marked.

\smallskip
``Missed'' angle coordinates $\vartheta_i$ can be now introduced similarly to
$\vartheta_{jk}$. For a boundary circle $\gamma_i$ we consider the adjacent
pants $S_j$. Fix a boundary component $\d_i S_j$ different from $\gamma_i$.
Let $J$ be a complex structure on $\Sigma$ such that boundary circles of all
pants $S_k$ are geodesic \wrt the intrinsic metric defined by $J$. Using
constructions from above, find a geodesic (resp.\ a ray) $\alpha_i$ starting
at point $x^*_i \in \gamma_i$ and ending at boundary circle $\d_i S_j$
(resp.\ approaching marked point $\d_i S_j$ of $\Sigma$). Take the marked
boundary point $\zeta_i$ on $\gamma_i$ and consider the length $\lambda_i$ of
the geodesic ark on $\gamma_i$, starting at $x^*_i$ and going to $\zeta_i$ in
the direction defined by the orientation of $\gamma_i$. Define $\vartheta_i
\deff {2\pi \lambda_i \over \ell_i }$, $\vartheta_i \in S^1\cong \rr/ 2\pi\zz$.
We include the coordinates $\vartheta_i$ into the system of angle coordinates
$\bfvartheta$. Denote by $\ttt_\Gamma$ the set of isomorphism classes of
complex structures on $\Sigma$ with marked boundary and with a given
decomposition into pants.

\smallskip
Let $C$ be a smooth complex curve with marked points and a marked boundary, $C
\bs \mapo = \cup_j C_j$ its decomposition into pants, $\sigma: \Sigma \to C$
a parametrization, and $\Sigma \bs \mapo = \cup_j S_j$ the induced
decomposition of $\Sigma$. To define a complex structure on $\ttt_\Gamma$, we
introduce special local holomorphic coordinates in a neighborhood of
boundary of pants on $C$. Consider some pants $S_j$ and its boundary circle
$\gamma^*$. It can be a boundary circle of $\Sigma$, $\gamma_i$ in our
previous notation, or a circle $\gamma_{jk}$ separating $S_j$ from another
pants $S_k$. Let $\ell^*$ be the intrinsic length of $\gamma^*$. Fix some
small $a>0$ and consider the annulus $A$ consisting of those $x\in S_j$ for
which the intrinsic distance $\dist(x,\gamma^*) <a$. The universal cover
$\ti A$ can be imbedded into hyperbolic plane $\hh$ as an infinite strip
$\Theta$ of constant width $a$, such that one of its borders is geodesic line
$L$. The action of a generator of $\pi_1(A) \cong \zz$ on $\ti A$ is defined
by the shift of $\Theta$ along $L$ by distance $\ell^*$.

Now consider the annulus $A' \deff [0, {\pi^2\over \ell^*}[ \times S^1$ with
coordinates $\rho, \theta$, $0\le\rho < {\pi^2\over \ell^*}$, $0\le \theta \le
2\pi$ and with the metric $({\ell^* \over 2\pi} / \cos {\ell^* \rho\over 2\pi}
)^2 (d\rho^2 + d\theta^2)$. A direct computation shows that this metric is of
constant curvature -1 and that boundary circle $\d_0A_1 \deff S^1 \times \{0\}$
is geodesic of length $\ell^*$, whereas $A'$ is complete in a neighborhood of
the other boundary circle. Consequently, the universal cover $\ti A{}'$ of
$A'$ can be imbedded in the hyperbolic plane $\hh$ as a hyperbolic half-plane
$\hh^+_L$ with a boundary line $L$, such that $\Theta \subset \hh^+_L$.
Moreover, the action of $\pi_1(A') \cong \zz$ on $\ti A{}' \cong \hh^+_L$ is
the same as for $\ti A \cong \Theta$. This shows that there exists an
{\sl isometric} imbedding of $A$ into $A'$ which maps $\gamma^*$ onto
$\d_0 A'$. Moreover, such an imbedding is unique up to rotations in the
coordinate $\theta$. This leads us to the following

\state Proposition 2.1. \it Let $C_j$ be pants with a complex structure
and $\gamma^*$ its boundary circle of the intrinsic length $\ell^*$. Let
$x^*$ be a point on $\gamma$. Then some collar annulus $A$ of $\gamma^*$
possesses the uniquely defined conformal coordinates $\theta\in S^1\cong \rr/
2\pi\zz$ and $\rho$, such that the intrinsic metric has the form $({\ell^*
\over 2\pi} / \cos {\ell^* \rho\over 2\pi})^2 (d\rho^2+ d\theta^2)$, $\rho|_{
\gamma^*} \equiv0$, $\theta(x^*) =0$, and such that the orientation on $S_j$
is given by $d\theta\wedge d\rho$. \rm

\medskip
We shall represent $\rho$ and $\theta$ also in the complex form $\zeta \deff
e^{-\rho + \isl \theta}$ and call $\zeta$ the {\sl intrinsic coordinate} of
the pants $C_j$ at $\gamma^*$. An important corollary of the description of
the intrinsic metric in a neighborhood of boundary circle is the following
statement about non-degenerating complex structures in pants, see {\sl
Definition 1.7}.

\state Lemma 2.2. \it Let $C$ be a smooth complex curve with marked points
admitting the intrinsic metric and let $\gamma^*$ be a boundary circle of $C$
of length $\ell^*$.

\sli If there exists an annulus $A \subset C$ of conformal radius $R$ (\ie $A
\cong \{\,z \in \cc \;:\; 1<|z| < R \;\})$, adjacent to $\gamma^*$ and
containing no marked points, then $\log R \le {\pi^2 \over \ell^*}$

\slii There exists a universal constant $a^*$ such that the condition $\ell^*
\le 1$ implies that there exist an annulus $A \subset C$ of conformal radius
$R$ with $\log R \ge { \pi^2\over \ell^*} -{2\pi\over a^*}$, which is
adjacent to $\gamma^*$, has area $a^*$ and contains no marked points of $C$.

\sliii Let $\gamma\subset C$ be a simple geodesic circle of the length $\ell$
and $A \subset C\bs\mapo$ annulus of conformal radius $R$ homotopy equivalent
to $\gamma$. Then $\log R \le {2\pi^2 \over \ell}$.
\rm

\state Proof. Let $\Omega$ be the universal cover of $C\bs \mapo$ equipped
with the intrinsic metric lifted from $C$. Then $\Omega$ can be isometrically
imbedded into the hyperbolic plane $\hh$ as a domain bounded by geodesic lines,
such that each of these lines  covers some boundary circle $\gamma_i$. Take
some (not unique!) line $L$ covering the circle $\gamma^*$ and fix a hyperbolic
half-plane $\hh^+_L$ with a boundary line $L$, so that $\Omega \subset
\hh^+_L$.

Now consider the universal cover $\ti A$ of the annulus $A$ and provide it
with the metric induced from $C$. Then we can isometrically imbed $\ti A$
in $\hh^+_L$ in such a way that the line covering $\gamma^* \subset \d A$
will be mapped onto $L$. The action of a generator of $\pi_1(A) \cong \zz$ on
$\ti A$ is defined by the shift of $\hh$ along $L$ onto distance $\ell^*$.
Consequently, $A$ can be isometrically imbedded into $\hh^+_L/ \pi_1(A)$, which
is the annulus $A'= [0, {\pi^2\over \ell^*}[ \times S^1$ with coordinates
$\rho, \theta$, $0\le\rho < {\pi^2\over \ell^*}$, $0\le \theta \le 2\pi$ and
with metric $({\ell^* \over 2\pi} / \cos {\ell^* \rho\over 2\pi})^2 (d\rho^2 +
d\theta^2)$. Note that the conformal radius of $A'$ is $e^{\pi^2/\ell^*}$. The
monotonicity of the conformal radius of annuli (see \eg [Ab], Ch.II, \S1.3)
yields the inequality $R\le e^{\pi^2/\ell^*}$ which is equivalent to first
assertion of the lemma.

\smallskip
Part \sliii of the lemma can be proved by same argument. More precisely,
under the hypothesis of part \sliii we imbed the annulus $A$ into the
annulus $A''= ]-{\pi^2\over \ell}, {\pi^2\over \ell}[ \times S^1$ with 
coordinates $\rho, \theta$, $- {\pi^2\over \ell} < \rho < {\pi^2\over \ell}$, 
$0\le \theta \le 2\pi$ and with metric 
$({\ell \over 2\pi} / \cos {\ell \rho \over 2\pi})^2 (d\rho^2 +d\theta^2)$.
The conformal radius of $A$ is now estimated by the conformal radius of $A''$,
which is equal to $e^{2\pi^2\over \ell}$.

\medskip
The second part of our lemma follows from results of Ch.II, \S\.3.3 of [Ab].
{\sl Lemma 2} there says that there exists a universal constant $a^*$ with the
following property: If $\ell^* \le 1$, then there exists a collar neighborhood
$A$ of constant width $\rho^*$ and of area $a^*$, which is an annulus imbedded
in $C$ and contains no marked points of $C$. In particular, we can extend the
intrinsic coordinates $\rho$ and $\theta$ in $A$. Using these coordinates, we
present $A$ in the form $\{ (\rho,\theta): 0\le \rho \le \rho^*\}$ and compute
the area,
$$
a^*=\area A = 2\pi \int_{\rho=0}^{\rho^*} \left(\msmall{\ell^*/2\pi \over
\cos(\ell^*\rho/2\pi) }\right)^2 d\rho =
\ell^* \tan\left(\msmall{\ell^*\rho^* \over 2\pi}\right).
$$
Consequently, $\tan\left({\pi\over2}-{\ell^*\rho^* \over 2\pi}\right) =
\cotan\left({\ell^*\rho^* \over 2\pi}\right) = {\ell^*\over a^*}$. This
implies ${\pi\over2}-{\ell^*\rho^* \over 2\pi} \le {\ell^*\over a^*}$, which
is equivalent to $\rho^* \ge { \pi^2\over \ell^*} -{2\pi\over a^*}$. To finish
the proof we note that the conformal radius $R$ of $A$ is equal $e^{\rho^*}$.
\qed

\medskip
Let $C$ be a smooth complex curve with marked points, $C_j$ a piece of a
decomposition of $C\bs \mapo$ into pants and $\gamma^*$ its boundary circle.
Then as a ``base point'' $x^*=\{\theta=0=\rho\}$ for the definition of the
intrinsic coordinate we shall use the point $x^*_{j,k}$ if $\sigma( \gamma^*
)$ is the geodesic separating $C_j$ from another pants $C_k$, or respectively
the point $x^*_i$ if $\gamma^*$ is a boundary circle of $C$. We denote these
coordinates $\zeta_{j,k}= e^{-\rho_{j,k} + \isl \theta_{j,k}}$ and $\zeta_i=
e^{-\rho_i + \isl \theta_i}$. Note that $\vartheta_i$ is exactly the
$\theta$-coordinate of the marked boundary point $x_i\in \gamma_i$ with
respect to $x^*_i$ and $\vartheta_{jk}$ is the $\theta$-coordinate of
$x^*_{k,j}$ with respect to $x^*_{j,k}$.

Note also that any intrinsic coordinate of a pair $(\zeta_{j,k}, \zeta_{k,j}
)$ extends canonically from one collar neighborhood of $\gamma_{jk}$ to
another side in such a way that the formula for the intrinsic metric
remains valid. This extension possesses the property $\zeta_{j,k}\cdot \zeta
_{k,j} \equiv e^{\isl \vartheta_{jk}}$ where $\vartheta_{jk}$ a constant
function. We can view this relation as the transition function from
$\zeta_{j,k}$ to $\zeta_{k,j}$.

A similar construction is possible in the case of a boundary circle
$\gamma_i$.
Namely, allowing $\rho_i$ to change also in the interval $]-{\pi^2 \over
\ell_i}, 0]$ and maintaining the formula $({\ell_i \over 2\pi} / \cos {\ell_i
\rho_i\over 2\pi})^2 (d\rho_i^2+ d\theta_i^2)$ for the metric we can glue to
$\Sigma$ an annulus $]-{\pi^2\over \ell_i}, 0] \times S^1$ and extend there
the coordinate $\zeta_i = e^{-\rho_i + \isl \theta_i}$.

Making such a construction with every boundary circle $\gamma_i$ we obtain a
complex curve $C^{(N)}$ with the following properties. $C$ is relatively
compact in $C^{(N)}$ and the intrinsic metric of $C$ extends to a complete
Riemannian metric on $C^{(N)}$ with constant curvature -1. Such the
extension and the metric are unique. $C^{(N)}$ is called the {\sl Nielsen
extension of $C$}, see [Ab]. Note that the complex coordinate $\zeta_i$ can
be extended further to the unit disk $\{|\zeta_i|<1\}$.

\smallskip
Using the introduced complex coordinates $\zeta_i$ and $\zeta_{j,k}$, we
define a deformation family of complex structures on the curve $C$ with
marked boundary. Let $\lambda_i$ and $\lambda_{jk}$ be complex parameters
changing in small neighborhoods of $e^{\isl\theta_i}$ and $e^{\isl \theta
_{jk}}$ respectively. Having these data $\bflambda =(\lambda_i, \lambda_{jk}
)$, construct a complex curve $C_\bflambda$ in the following way. Take the
pants $\{C_j\}$ of the given decomposition of $C$ and extend all the
complex coordinates $\zeta_i$ and $\zeta_{jk}$ outside the pants.
Glue the pairs of coordinates $(\zeta_{j,k},\zeta_{k,j})$ with new transition
relations $\zeta_{j,k}\cdot \zeta_{k,j}= \lambda_{jk}$ (constant functions).
Move original boundary circles $\gamma_i= \{|\zeta_i|=1\}$ of $C$ to new
positions defined by the equations $|\zeta_i|= |\lambda_i|$ and mark the
points $\zeta_i = \lambda_i$ on them.

\state Theorem 2.3. \it The natural map $F:\bflambda \to(\bfell,\bfvartheta)$
is non-degenerated. In particular, $\bflambda$ can be considered as the set
of local complex coordinates on $\ttt_\Gamma$ and $\calc\deff \{ C_\bflambda
\}$ as a (local) universal holomorphic family of curves over $\ttt_\Gamma$.
\rm

\state Proof. Write the functions $\bflambda = (\lambda_i, \lambda_{jk})$
in the form $\lambda_i= e^{-r_i + \isl\phi_i}$, $\lambda_{jk}= e^{-r_{jk} +
\isl\phi_{jk}}$. From the definition of the map
$$
\textstyle
F: (e^{-r_i + \isl\phi_i}, e^{-r_{jk} + \isl\phi_{jk}})
\mapsto (\ell_i, \ell_{jk}; \vartheta_i, \vartheta_{jk})
$$
it is easy to see
that ${\d (\vartheta_i, \vartheta_{jk}) \over \d(\phi_i, \phi_{jk})}$ is the
identity matrix, whereas ${\d (\ell_i, \ell_{jk}) \over \d(\phi_i, \phi_{jk}
)}$ is equal to 0. So it remains to show that the matrix ${\d (\ell_i,
\ell_{jk}) \over \d(r_i, r_{jk})}$ is non-degenerate.

\smallskip
Consider a special case when $C$ is pants with the boundary circles
$\gamma_i$ (at least one) and, possibly, with marked points $x_j$. We shall
consider such points as punctures of $C$. Let $J$ denote the complex structure
on $C$ and let $\mu_0$ be the intrinsic metric. Extend the coordinates
$\zeta_i$ and the metric $\mu_0$ outside of $\gamma_i$ to some bigger complex
curve $\wt C$ with $C \comp \wt C$.

Fix real numbers $v_i$ and consider the domains $C_t$ in $\wt C$ defined in
local coordinates $\zeta_i= e^{-\rho_i + \isl\theta_i}$ by inequalities the
$\rho_i \ge v_it$. This defines a family of deformations of $C= C_0$ parametrized
by a real parameter $t$, corresponding to a real curve in the parameter space
$\{ \bflambda \}$ given by $\lambda_i(t) = e^{v_it}$. Note that the
deformation is made in such a way that original complex structure $J$ and
local holomorphic coordinates are preserved. Thus we can use them as
``invariable basis'' in our calculations.

Let $\mu_t$ be the intrinsic metric of $C_t$. Without loss of generality we
may assume that $\mu_t$ extends to $\wt C$ as a metric with constant curvature
-1, which induces the original complex structure $J$ on $\wt C$. Since
the intrinsic metric depends smoothly on the operator $J$ of complex
structure, $\mu_t$ are smooth in $t$. In a local holomorphic coordinate $z=x
+\isl y$ we can present $\mu_t$ in the form $e^{2\psi(t,z)}(dx^2 + dy^2)$.
The condition ${\sf Curv}(\mu_t) \equiv -1$ is equivalent to the differential
equation
$$
\d^2_{xx}\psi(t,\cdot) + \d^2_{yy} \psi(t, \cdot)= e^{2\psi(t, \cdot)}
$$
where $\d_x$ denotes the partial derivation ${\d \over \d x}$ and so on.
Differentiating it in $t$ we get $e^{-2\psi(t,\cdot)}(\d^2_{xx}+ \d^2_{yy})
\dot \psi(t,\cdot) = \dot\psi(t,\cdot)$, where $\dot\psi(t,\cdot)$
denotes the derivative of $\psi(t,\cdot)$ in $t$.

Note that $\d_t\mu_t = 2 \dot \psi(t,\cdot)\cdot \mu_t $, so $\dot \psi_t
(\cdot)= \dot\psi(t, \cdot)$ is independent of the choice of a local
holomorphic coordinate $z= x+ \isl y$ and is defined globally. The equation
$e^{-2\psi_z(t,\cdot)} (\d^2_{xx}+ \d^2_{yy})\dot\psi_z(t,\cdot) = \dot\psi_z
(t,\cdot)$ can be rewritten in the form $\Delta_t \dot\psi_t = 2 \dot \psi_t
$, with $\Delta_t$ denoting the Laplace operator for the metric $\mu_t$.

The condition that the circle $\gamma_i(t) \deff \{ \rho_i = v_it\}$ is
$\mu_t $-geodesic means that the covariant derivative $\nabla_{\theta_i}
(\d_{\theta_i})$ of the vector field $\d_{\theta_i}$, the tangent vector field
to $\gamma_i(t)$, must be parallel to $\d_{\theta_i}$ along $\gamma_i(t)$.
Expressing this relation in local coordinates $\rho_i$ and $\theta_i$, we get
$\d_{\rho_i}\psi_i(t; v_it, \theta_i)=0$, where $\mu_t =e^{2 \psi_i (t;
\rho_i, \theta_i)}(\d\rho_i^2 + \d\theta_i^2)$ is a local representation of
the metric $\mu_t$. Deriving in $t$, we get $\d_{\rho_i}\dot \psi_i(t; v_it,
\theta_i) + v_i \d^2_{\rho_i\rho_i}\psi_i(t; v_it, \theta_i) =0$.

In the case $t=0$ we have $\psi_i(0; 0, \theta_i) \equiv \log{\ell_i \over
2\pi}$, a constant. Hence $\d^2_{\theta_i\theta_i}\psi_i(0; 0, \theta_i) =
e^{2 \psi_i(0; 0, \theta_i)} $ and $\d_{\rho_i}\dot \psi_i(0; 0, \theta_i) =
- v_i e^{2 \psi_i(0; 0, \theta_i)}= - v_i \left({\ell_i \over2\pi} \right)^2
$. On the other hand, $\d_{\rho_i} = - {\ell_i \over 2\pi} \d_\nu$ on
$\gamma_i(0) = \gamma_i$ where $\nu$ denotes the unit outer normal field to
$C_0=C$. Consider the integral $\int_C |d\dot\psi {}_0|^2 + 2 \dot\psi{}_0^2
d\mu_0$. Integrating by parts, we get
$$
\eqalignno{
&\int_C |d\dot\psi{}_0|^2 + 2 \dot\psi{}_0^2 d\mu_0 =
\int_C \dot\psi_0(2 \dot\psi_0 - \Delta_0 \dot\psi_0) d\mu_0
+ \int_{\d C} \dot\psi_0 \d_\nu\dot\psi_0 dl=
\cr
\noalign{\vskip0pt\allowbreak}
=&
\sum_i \int_{\gamma_i} \dot\psi_0 \d_\nu\dot\psi_0
\msmall{\ell_i \over 2\pi}
d\theta_i =
\sum_i - \int_{\gamma_i}  \dot\psi_0 \d_{\rho_i}\dot\psi_0
d\theta_i =
\sum_i \int_{\gamma_i}  \dot\psi_0 v_i \left(\msmall
{\ell_i \over 2\pi}\right)^2 d\theta_i =
\cr
\noalign{\vskip0pt\allowbreak}
=&
\sum_i \msmall{v_i \ell_i \over 2\pi} \int_{\gamma_i}  \dot\psi_0
e^{\psi_i(0; 0, \theta_i)} d\theta_i =
\sum_i \msmall{v_i \ell_i \over 2\pi} \int_{\gamma_i}  (\dot\psi_0 +
 v_i \d_{\rho_i} \psi_i(0; 0, \theta_i))e^{\psi_i(0; 0, \theta_i)} d\theta_i =
\cr
\noalign{\vskip0pt\allowbreak}
=&
\sum_i \msmall{v_i \ell_i \over 2\pi}
\left.\msmall{\d \over \d t}\right|_{t=0} \int_{\gamma_i}
e^{\psi_i(t; v_it, \theta_i)} d\theta_i =
\sum_i \msmall{v_i \ell_i \over 2\pi}
\left.\msmall{\d \over \d t}\right|_{t=0}\ell_i(t) =
\sum_i \msmall{ \ell_i \over 2\pi} v_i \dot\ell_i.
}
$$
Here $\dot\ell_i$ denotes the derivative of the length parameter $\ell_i$ for
the curve $C_t$ at $t=0$, so that
$(\dot\ell_1,\dot\ell_2,\dot\ell_3)= dF (v_1, v_2, v_3)$.
The obtained relation shows that the Jacobi matrix
$dF= {\d(\ell_i, \ell_2, \ell_3) \over \d(r_1, r_2, r_3)}$ is non-degenerate.
Otherwise there would exist a nonzero vector $(v_1, v_2, v_3)$ such that for
the deformation constructed above we get $\dot\ell_i=0$. But then $\dot\psi_0
\equiv 0$, which is a contradiction.

\smallskip
Now consider a general situation. Let $\Sigma$ be a real surface with marked
boundary, $C$ a smooth curve with marked points, $\sigma:\Sigma \to C$ a
parametrization, and $C \bs \mapo = \cup_j C_j$ a decomposition into pants
with a given graph $\Gamma$. Let $\{\gamma_i\}$ be the set of boundary circles
and $\{ \gamma_{jk} \}$ the set of circles lying between the pants $C_j$ and
$C_k$ respectively. Consider these pants separately. Then for any circle $\gamma
_{jk} = \gamma_{kj}$ we obtain 2 distinguished ones, $\gamma_{j,k}$ considered
as a boundary circle of $C_j$, and $\gamma_{k,j}$ considered as a boundary circle
of $C_k$. Take real numbers $\bfv\deff (v_i, v_{j,k}, v_{k,j})$ where $v_i$ is
associated with the circle $\gamma_i$, $v_{j,k}$ with $\gamma_{j,k}$, and
$v_{k,j}$ with $\gamma_{k,j}$ respectively. Let $C_j(t\bfv)$ denote the pants
obtained from $C_j$ by the above construction using the corresponding parameters
$v_i$ and $v_{j,k}$. For $\bfv$ lying in a small ball $B = \{ |\bfv| < \eps \}$
all such families $C_j(t\bfv)$ can be extended for all $t \in [-1,1]$. Thus over
$B$ we obtain a collection of deformation families $C_j(\bfv)$ of complex structure
on pants $C_j$.

Let $\ell_i(\bfv)$, $\ell_{j,k}(\bfv)$, and $\ell_{k,j}(\bfv)$ denote the
lengths of circles $\gamma_i$, $\gamma_{j,k}$, and $\gamma_{k,j}$ \wrt
obtained intrinsic metrics $\mu_j(\bfv)$ on $C_j(\bfv)$. Denote by $\dot\ell_i$
a linear functional $\d_t|_{t=0} \ell_i(t\bfv)$, and define $\dot\ell_{j, k}$
similarly. The explicit formula for an intrinsic metric near a
boundary circle showss that $C_j(\bfv)$ can be glued to $C_k(\bfv)$ along
$\gamma_{jk}$ exactly when $\ell_{j,k}(\bfv) = \ell_{k,j}(\bfv)$. Since the
Jacobian $\d\bfell(\bfv) \over \d \bfv$ is non-degenerate, the conditions
$\ell_{j,k}(\bfv) = \ell_{k,j}(\bfv)$ define a submanifold $V \subset B$
whose tangent space $T_0V$ is given by relations $\dot\ell_{j, k}= \dot\ell
_{k,j}$. Note that this defines a deformation family of complex structures on
$C$ over the base $V$ such that the map $\bfv \in V \mapsto \bfell(\bfv)$ is a
diffeomorphism.

We state that the set $(v_i, v_{j,k} + v_{k,j})$ is a system of coordinates
on $V$ in the neighborhood of $0\in V$. To prove this it is sufficient to
show that the linear map $\bfv=(v_i, v_{j,k}, v_{k,j}) \in T_0V \mapsto (v_i,
v_{j,k} + v_{k,j})$ is non-degenerate. If it would be not true, then there
would exist a nontrivial $\bfv=(v_i, v_{j,k}, v_{k,j})\in T_0V$ with $v_i=0$
and $v_{j,k} + v_{k,j}=0$. Let $\dot\ell_i=\dot\ell_i(\bfv)$, $\dot\ell
_{j,k}= \dot\ell_{j,k}(\bfv)$ and $\dot \ell_{k,j} =\dot \ell_{k,j} (\bfv)$
be the corresponding derivatives of length. Then $\dot\ell_{j,k} = \dot
\ell_{k,j}$ and
$$
0< \sum_i \msmall{ \ell_i \over 2\pi} v_i \dot\ell_i +
\sum_{j<k} \msmall{ \ell_{jk} \over 2\pi} v_{j,k} \dot\ell_{j,k} +
\sum_{j<k} \msmall{ \ell_{jk} \over 2\pi} v_{k,j} \dot\ell_{k,j}=
\sum_{j<k} \msmall{ \ell_{jk} \over 2\pi} (v_{j,k}+ v_{k,j}) \dot\ell_{j,k}
=0.
$$
The obtained contradiction leads us to the following conclusion: The
functions $v_i$ and $v_{j,k}+ v_{k,j}$ define a coordinate system
on $V$ equivalent to $\bfell=(\ell_i, \ell_{jk})$.

\smallskip
Let us return to the holomorphic deformation family of complex structures on
$C$, defined by complex parameters $\lambda_i= e^{-r_i + \isl\phi_i}$ and
$\lambda_{jk}= e^{-r_{jk} + \isl\phi_{jk}}$. It is easy to see that the
Jacobian $\d(v_i, v_{j,k}+ v_{k,j}) \over \d( r_i , r_{jk})$  at the point
$(r_i , r_{jk})=0$ is the identity matrix. This fact proves the statement of
the lemma. \qed

\state Remark. At this point we give a possible reason why the complex
(\ie holomorphic) structure introduced by the complex coordinates
$\bflambda$ can be regarded as natural. Let $C$ be a complex curve with
marked points and nonempty marked boundary. In the case when $C$ is a disk or
an annulus assume additionally that at least one inner point of $C$ is
marked. The in a neighborhood of every boundary circle $\gamma_i$ of $C$ we
can construct the intrinsic coordinate $\zeta_i$. Take 2 copies $C^+$ and
$C^-$ of $C$ and denote by $\tau$ the natural holomorphic map $\tau: C^\pm
\to C^\mp$ interchanging the copies. Denote by $\zeta_i^\pm$ the local
intrinsic coordinate on $C^\pm$ at boundary circles $\gamma_i^\pm$, both
corresponding to $\gamma_i$. Now we can glue $C^+$ and $C^-$ together along
every pair of circles $(\gamma_i^+, \gamma_i^-)$ by setting $\zeta_i^+ \cdot
\zeta_i^- =1$ as transition relations. We obtain a closed complex curve $C^d$
which admits a natural holomorphic involution $\tau: C^d \to C^d$. For the
constructed family $\{C_\bflambda \}$ the corresponding family $\{C^d
_\bflambda \}$ will be holomorphic. In fact, the statement of {\sl Theorem
2.3} means that $\{ C^d _\bflambda \}$ is a minimal complete family of
deformation of $C^d$ in the class of curves with holomorphic involution. This
construction of doubling should not be confused with another construction of
the {\sl Schottky double} $C^{Sch}$ of $C$ which provides an {\sl
antiholomorphic} involution $\tau^{Sch}: C^{Sch} \to C^{Sch}$. We shall use
the Schottky double $C^{Sch}$ in {\sl Section 5} considering curves with
totally real boundary conditions.

\smallskip
The construction of (holomorphic) double $C^d$ shows how to give an invariant
description of holomorphic structure on $\ttt_\Gamma$. Let $C$ be a smooth
complex curve with marked points and marked boundary, and $x\in \ttt_\Gamma$
the corresponding point on moduli space. Denote by $D$ the divisor of marked
points. If curve the $C$ is not closed and $C^d$ is its double with holomorphic
involution $\tau$, we denote by $D^d \deff D + \tau(D)$ the double of $D$.

\state Lemma 2.4. \it If $C$ is closed, then the tangent space $T_x \ttt
_\Gamma$ is naturally isomorphic to $\sfh^1(C, \calo(TC)\otimes \calo(-D) )$.

If $C$ is not closed, then the space $T_x\ttt_\Gamma$ is naturally isomorphic
to the space $\sfh^1(C^d, \calo(TC^d)\otimes \calo(-D^d))^\tau$ of
$\tau$-invariant elements in $\sfh^1(C^d, \calo(TC^d)\otimes \calo(-D^d))$.

In both cases the complex structure on $T_x\ttt_\Gamma$ induced by local
complex coordinates $\bflambda$ coincides with those from $\sfh^1(C,
\calo(TC) \otimes \calo(-D))$ (resp.\ $\sfh^1(C^d, \calo(TC^d) \otimes
\calo(-D^d) )^{(\tau)}$). In particular, this defines a global complex
structure on space $\ttt_\Gamma$. \rm

\state Proof. The part concerning closed curves is well-known. In fact, the
natural isomorphism $\psi: T_x\ttt_\Gamma \to \sfh^1(C, \calo(TC)\otimes
\calo(-D))$ is a Kodaira-Spencer map. Its description is very simple in the
introduced local coordinates $\zeta_{j,k}$ on $C$ and $\bflambda =( \lambda
_{jk})$ on $\ttt_\Gamma$. Let $C\bs \mapo = \cup C_j$ be the decomposition of
$C$ into pants with the graph $\Gamma$. For every pants $C_j$ choose an open
set $\ti C_j$, containing a closure $\barr C_j= C_j \cup \d C_j$. Without
loss of generality we may assume that $\ti C_j$ are chosen not too big, so
that the covering $\calu\deff \{ \ti C_j \}$ is acyclic for the sheaf
$\calo(TC)$ and that the local coordinates $\zeta_{j,k}$ are well-defined in
the intersections $\ti C_j \cap \ti C_k$. Then vector $v \in T_x\ttt_\Gamma$
with local representation $v= \sum_{j<k} v_{jk} {\d \over \d \lambda_{jk}}$
is mapped by Kodaira-Spencer map $\psi$ to the \v{C}ech 1-cohomology class
$$
\psi(v) \in \sfh^1(C, \calo(TC)\otimes \calo(-D)) \cong
\check \sfh^1(\calu, \calo(TC) \otimes \calo(-D)),
$$
represented by the 1-cocycle
$$
\left(v_{jk} \zeta_{j,k} \msmall{\d \over \d \zeta_{j,k}} \right)
\in \prod_{j<k} \Gamma(\ti C_j \cap \ti C_k,  \calo(TC)\otimes \calo(-D)).
$$
For more details see [D-G].

\smallskip
Using this description of Kodaira-Spencer map for closed curves with marked
points, it is easy to handle the case of curves with boundary. Let $C$ be
a non-compact curve with marked points and with decomposition $C \bs \mapo
= \cup_j C_j$. Take its double $C^d$ with the involution $\tau$. Then the
decomposition of $C$ induces a $\tau$-invariant decomposition
$C^d =\bigcap_j(C_j  \cap \tau C_j)$. The corresponding covering $\calu^d$
of $C^d$ can be also chosen to be $\tau$-invariant.

The local coordinates $\zeta_i$, corresponding to boundary circles $\gamma_i$
of $C$, can be now extended to a both side neighborhood of $\gamma_i$ in
$C^d$. The coordinates $\zeta_{j,k}$, corresponding to inner circles
$\gamma_{jk}$, induce local complex coordinates $\zeta^\tau_{j,k}\deff
\zeta_{j,k} \scirc \tau$ in $\tau( \ti C_j \cap \ti C_k)$.

Any deformation of the complex structure on $C$ induces a deformation of
the complex structure on $C^d$. This defines a map $\phi: \ttt_\Gamma \to
\ttt_{\Gamma^d}$, with $\Gamma^d$ denoting the graph corresponding to the
$\tau$-invariant decomposition of $C^d$ into pants. Using introduced
coordinates, we present a tangent vector $v \in T_x \ttt_\Gamma$ in the form
$$
v= \sum_i v_i \msmall{\d \over \d \lambda_i} +
\sum_{j<k} v_{jk} \msmall{\d \over \d \lambda_{jk}}.
$$
Then the composition of the Kodaira-Spencer map $\psi^d$ of $C^d$ with
the differential of $\phi: \ttt_\Gamma \to \ttt_{\Gamma^d}$ maps $v$ to
$$
\psi^d \scirc d\phi(v) \in \sfh^1(C^d, \calo(TC^d)\otimes \calo(-D^d)) \cong
\check \sfh^1(\calu^d, \calo(TC^d) \otimes \calo(-D^d)),
$$
represented by the \v{C}ech 1-cocycle
$$
\check v \deff \left(v_i \zeta_i \msmall{\d \over \d \zeta_i},
v_{jk} \zeta_{j,k} \msmall{\d \over \d \zeta_{j,k}},
v_{jk} \zeta^\tau_{j,k} \msmall{\d \over \d \zeta^\tau_{j,k}} \right)
\in \prod_i \Gamma(\ti C(i) \cap \tau C(i),  \calo(TC^d)\otimes \calo(-D^d))
\times
$$
$$
\prod_{j<k} \Gamma(\ti C_j \cap \ti C_k,  \calo(TC^d)\otimes \calo(-D^d))
\times \prod_{j<k} \Gamma(\tau(\ti C_j \cap \ti C_k),
\calo(TC^d)\otimes \calo(-D^d)),
$$
where $C(i)$ denotes the pants of $C$ adjacent to circle $\gamma_i$. It is
obvious that if all $v_i$ vanish, then this \v{C}ech 1-cocycle is $\tau
$-invariant. On the other hand, the relation $\zeta_i \cdot (\zeta_i \scirc
\tau) \equiv \lambda_i= \const$ implies that $\tau_*(\zeta_i {\d\over \d
\zeta_i}) = -\zeta_i {\d \over \d \zeta_i}$. The additional change of sign of
the corresponding part of cocycle $\check v$ comes from the fact that $\tau$
interchange $C(i)$ with $\tau(C(i))$. This shows that $\check v$ is $\tau
$-invariant and the statement of the lemma follows. \qed

\medskip
Now we study the connection between the geometry of $\ttt_\Gamma$ and the
degeneration of complex structures on a real surface $\Sigma$ with marked
points and marked boundary. Let $\Sigma \bs \mapo= \cup_j S_j$ be a
decomposition into pants with graph $\Gamma$. The Fenchel-Nielsen coordinates
on $\ttt_\Gamma$ define a map $(\bflambda, \bfvartheta): \ttt_\Gamma \to (\rr
\times S^1)^{3g-3+ m+2b}$, which is a diffeomorphism by {\sl Proposition
1.2}. So, if $\{ j_n\}$ is a sequence of complex structures on $\Sigma$, its
degeneration means that the sequence of Fenchel-Nielsen coordinates of $\{
j_n\}$ is not bounded in $(\rr \times S^1)^{3g-3+ m+2b}$.

One can see that, in fact, we have two types of the degeneration. The first
one occurs when the maximum of the length coordinates $\ell_i$ and $\ell_{jk}$
of $j_n$ increases infinitely, and the second one is present
when minimum of the length coordinates of $j_n$ vanishes. It should be pointed
out that for an appropriate sequence one can have both types of degeneration.

Note that by {\sl Proposition 1.2} the Fenchel-Nielsen coordinates of a
complex structure $j$ on $\Sigma$ are defined by a choice of a topological
type of decomposition of $\Sigma$ into pants, encoded in graph the $\Gamma$.
Thus the introduced notion of degeneration also depends on the choice of
$\Gamma$. Possibly, the best choice of such decomposition is
established by the following statement, proved in [Ab], Ch.II, \S\.3.3.

\state Proposition 2.5. \it Let $C$ be a complex curve with parametrization
$\sigma: \Sigma \to C$. Then

{\sl a)} there exists a universal constant $l^*>0$ such that any two geodesic
circles $\gamma'$ and $\gamma''$ on $C$ satisfying $\ell( \gamma') < l^*$ and
$\ell(\gamma'') < l^*$ are either disjoint or coincide;

{\sl b)} there exists pants decomposition $C \bs \mapo= \cup_j S_j$ such
that the lengths of inner boundary circles $\gamma_{jk} = \barr S_j \cap \barr
S_k$ are bounded from above by a constant $L$ which depends only on the topology
of\/ $\Sigma$ and the maximum $M$ of the lengths of the boundary circles of $C$;
moreover, any simple geodesic circle $\gamma$ on $C$ with $\ell(\gamma) <l^*$
occurs as a boundary circle of some $C_j$.

\smallskip
\state Corollary 2.6. \it Let $C_n$ be a sequence of nodal curves parametrized
by a real surface $\Sigma$ with uniformly bounded number of components.
Suppose that the complex structures of $C_n$ do not degenerate near boundary.
Then, passing to a subsequence, one can find a decomposition $\Sigma =\cup_j
S_j$ and new parametrizations $\sigma'_n : \Sigma \to C_n$, such that:

\sli the decomposition $\Sigma =\cup_j S_j$ induces a decomposition of every
non-exceptional component of $C_n$ into pants, whose boundary circles are
geodesics;

\slii the intrinsic length of these geodesics are bounded uniformly in $n$.
\rm

\state Proof. By {\sl Lemma 2.2}, the intrinsic lengths of boundary circles
of non-exceptional components of $C_n$ are bounded uniformly in $n$. Find a
decomposition into pants of every non-exceptional component of $C_n$
satisfying the conditions of part {\sl b)} of {\sl Proposition 2.5}. 
Let $\Gamma_n$
denotes the obtained graph of the decomposition of $C_n$. Since the number of
components of $C_n$ is uniformly bounded, we can a subsequence $C_{n_k}$ with
the same graph $\Gamma$ for all $n_k$.

It follows from the proof in [Ab], that the constant $L$ from part {\sl b)} of
{\sl Proposition 2.5} depends continuously on the maximum $M$ of the lengths 
of the boundary circles of $C_n$. This implies condition {\sl i$\!$i)}.
Applying {\sl Proposition 1.2}, we complete the proof. \qed

\bigskip\bigskip
\newsection{}{Apriori estimates.}

\medskip\noindent
\rm
Let $(X,J)$ be an almost complex manifold. In what follows the  tensor $J$ is
supposed to be only continuous, \ie of class $C^0$. Fix some Riemannian metric
$h$ on $X$. All norms and distances will be taken with respect to $h$.
In particular, we have the following

\state Definition 3.1. {\sl A continuous almost complex structure $J$ is
called {\it uniformly continuous on $A \subset X$ with respect to $h$}, if
$\norm{J}_{L^\infty(A)} < \infty$ and for any $\epsi >0$ there exists
$\delta=\delta(J, A, h) >0$ such that for any $x\in A$ one can find a
$C^1$-diffeomorphism $\phi: B(x, \delta) \to B(0, \delta)$ from the ball
$B(x, \delta)\deff \{y\in X \;:\; \dist_h(x,y) <\delta \}$ onto the standard
ball in $\cc^n$ with the standard metric $h\st$, such that
$$
\norm{J -\phi^*J\st}_{L^\infty(B(x,\delta)\cap A)} +
\norm{h -\phi^*h\st}_{L^\infty(B
(x,\delta)\cap A)}
\le\epsi.
$$
}

Roughly speaking, this means that on the set $A$ we can $\norm\cdot _{L^\infty}
$-approximate $J$ by an integrable structure in $h$-metric balls of radius
independent of $x\in A$. The function $\mu(J,A,h)$ whose value at $\epsi>0$ is
the biggest possible $\delta\le1$ with the above property is called the {\it
modulus of uniform continuity of $J$ on $A$}. Note that every continuous almost
complex structure $J$ is always uniformly continuous on {\it relatively
compact} subsets $K \Subset X$.

\smallskip
Let $J^*$ be a continuous almost complex structure on $X$ and $A\subset X$
a subset. Assume that $J^*$ is uniformly continuous on $A$ and denote
by $\mu_{J^*}= \mu( J^*,A, h)$ the modulus of uniform continuity of
$J^*$ on $A$.

\state Lemma 3.1. {\sl (First Apriori Estimate). \it For every $p$ with
$2< p<\infty$ there exists an $\eps_1 =\eps_1(\mu_{J^*}, A, h)$ (independent
of $p$) and $C_p=C(p,\mu_{J^*}, A, h)$, such that for any continuous almost
complex structure $J$ with $\norm{ J- J^*} _{L^\infty(A)}<\eps_1 $ and
for every $J$-holomorphic map $u\in C^0\cap L^{1,2}(\Delta ,X)$, satisfying
$u(\Delta )\subset A$ and $\norm{du} _{L^2( \Delta )} <\eps_1 $ one has
the estimate
$$
 \norm{du}_{L^p({1\over 2}\Delta )}\le C_p\cdot
  \norm{du}_{L^2(\Delta )}.\eqno(3.1)
$$
}

\smallskip
\state Proof.
\sl Step 1. First, we prove first inequality $(4.1)$ for the case when
$A\subset U \subset \cc^n$, $h$ is the Euclidean metric, and $J^*$ is the
standard complex structure in $\cc^n=\rr^{2n}$.\rm

In the Schwarz spaces $\cals(\cc)= \cals(\cc,\cc^n)$ and $\cals'(\cc)=
\cals'(\cc, \cc^n)$ we consider the Cauchy-Green operators $\d= {\d \over \d
z}$, $\dbar= {\d\over \d \bar z}$, $T=T_{CG}={1 \over 2\pi\isl z}*(\cdot )$
and $\barr T={1 \over 2\pi\isl\bar z} * (\cdot)$, where the star $*$ denotes
the convolution of distributions. Note that operators $T$ and $\barr T$ map
$\cals$ only to $\cals'$ and not in $\cals$. Nevertheless one has the
following identities in the spaces  $\cals$ and $L^p(\cc )$:
$$
\dbar
\scirc T=T\scirc \dbar = \id
\qquad \text{and} \qquad
\d \scirc \barr T=\barr T\scirc \d = \id.
$$
Recall also that the Calderon-Zygmund inequality states that for any $p$,
$1<p< \infty$, there exists a constant $C_p$ such that for any $f\in L^p(\cc )$
one has
$$
\norm{ (\d \scirc T)(f)} _{L^p(\cc )}\le C_p\cdot \norm{f}_{L^p(\cc )}
\qquad \text{and} \qquad
\norm{ (\dbar \scirc \barr T)(f)}_{L^p(\cc )}\le C_p\cdot
\norm{f}_{L^p(\cc )}.
$$
This implies that taking any $f\in L^p(\cc)$ and setting $g \deff Tf$ (or
$g \deff \barr Tf$) we get the regularity property $g\in L^{1,p}_\loc(\cc)$
with the estimate $\norm{dg}_{L^p(\cc )}\le (1+C_p)\norm{f }_{L^p(\cc )}$.

Consider now a continuous linear complex structure $J^*(z)$ in the trivial
bundle $\cc \times \rr^{2n}\to \cc $. This meant that $J^*(z)$ is a continuous
family of endomorphisms $\rr^{2n}\to \rr^{2n}$ with $J^*(z)^2=-\id$. Define
an operator $\dbar_{J^*}:\cals' (\cc ,\rr^{2n})\to \cals'(\cc ,\rr^{2n})$ by
formula
$$
(\dbar_{J^*}f)(z)={1\over 2}[\d_x f(z) + J^*(z)\d_y f(z)].
$$
If $J$ is another continuous complex structure in the bundle $\cc \times
\rr^{2n}$ then for $f\in L^p(\cc ,\rr^{2n})$ one has the estimate
$$
\norm{ (\dbar_J\scirc T - \dbar_{J^*}\scirc T)f}_{L^p(\cc )}\le
\norm{J - J^*}_{L^\infty(\cc )}\cdot \norm{d(Tf)}_{L^p(\cc)}
\le
$$
$$
\le \norm{J - J^*}_{L^\infty(\cc )}(1+C_p)\norm{f}_{L^p(\cc )}.
\eqno(3.2)
$$
If we take $J^*(z)\equiv J\st$, the standard structure in $\cc^n$, then,
according to the above remark, $\dbar_{J^*}\scirc T: L^p(\cc ,\cc^n)\to L^p(\cc,
\cc^n)$ is the identity. From (3.2) we see that if $\norm{J - J^*}<\eps_p
\deff{1\over 1+C_p}$ then $\dbar_{J}\scirc T: L^p(\cc, \cc^n) \to L^p( \cc,
\cc^n)$ is an isomorphism. Moreover, since $\dbar_J \scirc T = \dbar_{J^*}
\scirc T + (\dbar_J -\dbar_{J^*})\scirc T$, we have
$$
(\dbar_J\scirc T)\inv = (\id + (\dbar_J-\dbar_{J^*})\scirc T)^{-1} =
\Sigma_{n=0}^\infty(-1)^n[(\dbar_J -
\dbar_{J^*})\scirc T]^n.\eqno(3.3)
$$
This shows, in particular, that the operator $(\dbar_J\scirc T)\inv $ doesn't
depend on the choice of $p>1$. Now we shall prove the following

\state Lemma 3.2. {\it For any $u\in L^{1,2}(\cc ,\rr^{2n})$ with compact
support and any continuous $J$ with $\norm{J-J\st}_{L^\infty(\cc, \End(
\rr^{2n}))}< \eps_p$ the condition $\dbar_Ju\in L^p(\cc ,\rr^{2n})$ implies
$$
\norm{du }_{L^p(\cc )}\le C\cdot \norm{ \dbar_Ju}_{L^p(\cc )}.
\eqno(3.4)
$$
for some $C=C(p,\norm{J-J\st}_{L^\infty(\cc )})$.
}

\smallskip
\state Proof. Put $v=u-T\scirc \dbar_{J\st} u$. Then $\dbar_{J\st}v=0$.
So $v$ is holomorphic and decreases at infinity. Thus $v=0$, which implies
$u=(T\scirc \dbar_{J\st})u$. By Calderon-Zygmund inequality, in order to
estimate
$\norm{ du}_{L^p(\cc )}$ it is sufficient to estimate $\norm{ \dbar_{J\st}u}
_{L^p(\cc )}$. Using $(\dbar_J\scirc T)\scirc \dbar_{J\st}u =\dbar_J u
\subset L^p(\cc )\cap L^2(\cc )$ and (3.3) we get that $\dbar_{J\st}u\in
L^p(\cc )\cap L^2(\cc )$ with the estimate
$$
\norm{ \dbar_{J\st}u}_{L^p(\cc )}\le \Sigma_{n=0}^\infty\norm{
(\dbar_J-\dbar_{J\st})\scirc T}_p^n\cdot \norm{ \dbar_Ju}_p\le
C\cdot \norm{ \dbar_Ju}_p,
$$
which yields (3.4). 

\smallskip
To finish {\sl Step 1}, consider a $J$-holomorphic map $u:\Delta \to (\rr^n,J)$,
with $u(\Delta )\subset A$ and $\norm{J-J\st}<\eps_p$. Define a linear complex
structure in the bundle $\Delta \times \rr^{2n}$ setting $J(z) \deff J(u(z))$.
Then $u$ is a $J$-holomorphic section of $(\Delta \times \rr^{2n}, J)$ and
$\norm{J-J\st}_{L^\infty (\Delta )}< \eps_p$. Extend $J$ to $\cc \times
\rr^{2n}$ with the same estimate.

Let $\psi$ be a non-negative cut-off function supported in
$\Delta(0, {3\over4})$ which is identically 1 on $\Delta(0,\half)$.
Put $u_1 \deff u\psi $. Then $u_1\in L^{1,2}(\Delta )$ and $\dbar_Ju_1=
u\dbar_J\psi \in L^p(\cc )$ with
$\norm{ \dbar_Ju_1}_{L^p(\Delta )}=\norm{ u\dbar_J\psi }_{L^p(\Delta )} \le
C\norm{du}_{L^2(\Delta )}$. Here we use the Sobolev imbedding $L^{1,2}(
\Delta ,\cc )\to L^p(\Delta ,\cc )$,  $p<\infty $. Now (3.4) applies to
get the estimate of {\sl Step 1}.

Using the Sobolev imbedding $L^{1,p}\subset C^{1-{2\over p}}$ and obvious
properties of $L^p$-norms by dilations, one derives easily from {\sl Step 1}
the following property

\medskip\noindent\sl
Step 2. Fix $2<p<\infty $. There exists $\eps _2=\eps_2(\mu_{J^*}, A,h)> 0$
such that for any $u$ and $J$ as in {\sl Lemma 3.1} with 
$\diam(u(\Delta (x,r))) <\eps_2$  one has the estimate
$$
\diam(u(\Delta (x,{\textstyle{r\over 2}})))\le C_1 r^{1-{2\over p}} \cdot
\norm{ du}_{L^p(\Delta (x,r/2))}\le C_2r^{1-{2\over p}}\cdot
\norm{ du}_{L^2(\Delta (x,r))}.
\eqno(3.5)
$$
for any disc $\Delta (x,r)\subset \Delta$. 
\rm

\medskip
Now consider the function
$$
  \alpha(r) \deff \left\{ \matrix \format \l&\ \c\ &\r&\c&\l\\
1     & \text{ if }&        &r& \le 1/2 \cr
3-4r  & \text{ if }&1/2 \le &r& \le 3/4 \cr
0     & \text{ if }&3/4 \le &r&         \endmatrix
\right.
$$
For $x\in \Delta$ set
$$
f(x)\deff\max \left\{ t \in [0, {\textstyle {1\over8}}] \;:\;
\diam(u(\barr\Delta (x,t\cdot \alpha (| x| ))))\le \eps_2\right\}.
$$
Clearly, $f$ is continuous and $f\equiv {1\over8}$ if ${3\over4}\le |x| <1$.

\smallskip\noindent\sl
Step 3. $f(x)\equiv {1\over 8}$.

\smallskip\rm
Suppose that there is an $x_0$ with $f(x_0)=\min \{ f(x)\;:\; x\in \Delta \}
< {1\over8}$. It is clear that $f(x_0) >0$.

Take the disk $\Delta (x_0,a)$ with $a\deff f(x_0)\alpha (| x_0| )$. Note that
$$
\diam(u(\Delta (x_0,a)))=\eps_2.
\eqno(3.6)
$$
Using the Sobolev embedding $L^{1,4}(\Delta )\subset C^{0, {1\over 2}}(\Delta )$,
estimate (3.5), and relation (3.6), we obtain that $\diam(u(\Delta (x_0,{a\over
2}))) \le C\cdot \norm{du}_{L^2(\Delta(x_0,a))}$. Take a point $x_1\in \Delta (x
_0, a)$ with $|x_1 -x_0| \le {3\over 4}a$. Since $f(x_0)={a\over \alpha (| x_0|)}$
is the minimum of $f$, we have that $f(x_1)\ge {a\over \alpha (| x_0| )}$ and,
thus, $f(x_1) \alpha (| x_0| )\ge a$. At the same time $\alpha (| x_1| )\ge
\alpha (| x_0| )- 3a$, so $f(x_1)\alpha (| x_1 | )\ge a-3a^2\ge {a\over 2}$
because $a\le {1\over 8}$. This means that $\diam(u(\Delta(x_1,{a\over 2})))\le
\eps_2$ and so $\diam(u(\Delta(x_1,{a\over 4})))\le C\cdot \norm{du}_{L^2 (\Delta
 )}$. So $\diam(u(\Delta(x_0,a)))\le 4C\cdot | | du| |_{L^2(\Delta )}$. If
$\eps $ is taken smaller than ${\eps_2\over 4C}$ then we obtain a contradiction
with (3.4). {\sl Step 3} is completed.

This means that $\diam(u(\Delta(x,{1\over 8})))\le \eps_2$ for any $x\in
\Delta(0, {1\over2})$. So {\sl Step 2} with $r={1\over 8}$ gives us the
assertion of {\sl Lemma 3.1}. \qed

\medskip
This lemma can be used to prove that a $J$-holomorphic map $u:\Delta \to(X,
J)$ is $L^{1,p}$-smooth for any $p<\infty $, provided $J$ is continuous.
To show this, we note that $u_{\eps }(z)\deff u(\eps z)$ is also
$J$-holomorphic and $\norm{ du_{\eps }}_{L^2(\Delta (0,1))} = \norm{du}
_{L^2(\Delta (0,\eps ))}$. For $\eps $ small enough we obtain $\norm{
du_\eps}_{ L^2( \Delta(0,1) )} = \norm{du}_{L^2(\Delta (0,\eps ))}< \eps_1$,
where $\eps_1$ is defined in from {\sl Lemma 3.1}.
Now, estimate (3.1) gives us the $L^{1,p}$-continuity of
$u_\eps$ and thus of $u$ in the neighborhood of zero. Another immediate
consequence of the First Apriori Estimate (3.1) is the following

\smallskip
\state Corollary 3.3. {\it Let $X$ be a manifold, $h$ some metric, and $\{ J_n
\}$ a sequence of continuous almost complex structures on $X$ such that
$J_n\longto J$ in $C^0$-topology on $X$. Let $A\subset X$ be a closed
$h$-complete subset, such that $J$ is uniformly continuous on $A$
\wrt $h$.

Let $u_n\in C^0\cap L^{1,2}_\loc(\Delta ,X)$ be a sequence of
$J_n$-holomorphic maps such that $u_n(\Delta)\subset A$, $\norm{du_n}_{L^2(
\Delta )}\le \eps_1$, and $u_n(0)$ is bounded in $X$. Then there exists a
subsequence $\{ u_{n_k} \}$ which $L^{1,p}_\loc$-converges to a $J$-holomorphic
map $u_\infty$ for all $p< \infty$.

In particular, for any $K \Subset \Delta$ norms $\norm{du_{n_k}}_{L^2(K)}$
tend to $\norm{du_\infty}_{L^2(K)}$.
}

\state Proof. First Apriori Estimate of {\sl Lemma 3.1} with the  Sobolev
imbedding $L^{1,p} (\Delta) \subset C^{1-{2\over p}}(\Delta)$ implies that
for every $r<1$ the sequence of the sets $\{u_n(\Delta(0,r))\}$ is uniformly
bounded in $X$. Consequently, $u_n(\Delta (0,r)) \subset K_r$ for some
relatively compact subset $K_r \Subset X$ independent of $n$. This implies
the existence of a subsequence $\{ u_{n_k} \} $ which converges to $u_\infty$
in $C^\alpha (\Delta, X)$ for any $\alpha <1$.

To show the (strong) $L^{1,p}_\loc(\Delta)$-convergence of $\{ u_{n_k}\}$,
take any $x\in\Delta$. Then we can find $r>0$, such that all images
$u_{n_k} (\barr\Delta(x,r))$ lie in some chart $U \subset X$. Moreover, we may
assume that $U$ is a domain in $\cc^n$, such that $\norm{J -J\st}_{L^\infty(U)}
\le \epsi_1$. Take $\phi \in C_0^\infty(\Delta(x,r) )$ with
$\phi\ogran_{\Delta (x, r/2)} \equiv 1$. Then
$$
\dbar_{ J_{n_k}} (\phi u_{n_k})=\d_x (\phi
u_{n_k}) +J_{n_k} (u_{n_k}) \d_y( \phi u_{n_k})= (\d_x\phi +J_{n_k} (u_{n_k})
\d_y\phi) u_{n_k},
$$
which is $C^0$-bounded and thus $L^p$-convergent for any  $p< \infty$.
The estimate (3.4) for $\dbar_{ J_{n_k}}$-operator gives us
the $L^{1,p}$-convergence of $\{ u_{n_k}\}$ on $\Delta(x,r/2)$. \qed

\smallskip\state Definition 3.2. {\sl Define a {\it cylinder} $Z(a,b)$ by
$Z (a,b) \deff S^1 \times [a,b]$, equiping it with coordinates
$\theta \in [0,2\pi]$, $t\in[a,b]$, with metric $ds^2= d\theta^2 + dt^2$ and
the complex structure $J\st({\d \over \d\theta}) ={\d \over \d t}$. Denote
$Z_i\deff Z(i-1,1)= S^1 \times [i-1,i]$.}

\smallskip
Let $J^*$ be some continuous almost complex structure on $X$ and $A$ a
subset of $X$, such that $J^*$ is uniformly continuous on $A$. Let
$\mu_{J^*}$ denote the modulus of uniform continuity of $J^*$ on $A$.

\smallskip
\state Lemma 3.4. {\sl (Second Apriori Estimate). \it There exist constants
$\gamma \in ]\,0,1\,[$ and $\eps_2=\eps_2(\mu_{J^*}, A,h)>0$ such that for any
$J$ with $\norm{ J- J^*} <\eps_2$ and every $J$-holomorphic map $u:Z(0,5)\to X$
with $u(Z(0,5)) \subset A$ the condition $\norm{du}_{ L^2(Z_i)}<\eps_2$ for
$i=1,\ldots,5$ implies
$$
\norm{du}^2_{L^2(Z_3)}\le {\gamma \over 2}\bigl(\norm{du}^2_{L^2(Z_2)} +
\norm{du}^2_{L^2(Z_4)}\bigr).
\eqno(3.7)
$$
}

\state Proof. Take $\eps_2>0$ small enough, such that $\mu_{J^*}(\epsi_2) <
\epsi_1$, where $\eps_1$ is the constant from {\sl Lemma 3.1}. Then for any
$A' \subset A$ the condition $\diam(A') \le \epsi_2$ implies that $\osc(J^*,
A') \le \eps_1$. Due to {\sl Lemma 3.1}, we may assume that $u(Z_i)\subset B$
for $i=2,3,4$, where $B$ is a small ball in $\rr^{2n} =\cc^2$ with the
structure $J\st$. Moreover, we may assume that $\norm{ J^* - J\st}_{
L^\infty(B)} \le \eps_1$.

Find $v\in C^0\cap L^{1,2}(Z(1,4),\cc^n)$ such that $\dbar_{J\st}v=0$ and
$\norm{du - dv}_{L^2(Z(1,4))}$ is minimal. We have
$$
\norm{  \dbar_{J\st}(u-v)}_{L^2(Z_i)} = \norm{(J\st-J(u))\d_y u
}_{L^2(Z_i)}
\le \norm{J\st - J}_{L^\infty(B)}\norm{du}_{L^2(Z_i)}.
$$
So for $i=2,3,4$ we get
$$
\norm{du - dv}_{L^2(Z_i)}\le C\norm{J\st -J}_{L^\infty(B)}
\norm{ du
}_{L^2(Z(1,4))}.\eqno(3.8)
$$
Now let us check the inequality (3.7) for $v$. Write $v(z)= \Sigma_{k=
-\infty} ^\infty v_k e^{k(t+i\theta)}$. Then $\norm{dv}^2_{L^2 (S\times \{
t\} )}=4\pi \Sigma_{k=-\infty}^\infty k^2| v_k|^2 e^{2kt}$. Since obviously
$$
\int_2^3e^{2kt}\le {\gamma_1\over 2}
\left( \int_1^2e^{2kt}dt+\int_3^4e^{2kt}dt \right)
$$
for all $k \not=0$ with $\gamma_1={2\over e^2}$, one gets the required
estimate for all holomorphic $v$.

Using (3.8) with $\norm{J\st - J}_{L^\infty}$ sufficiently small, we conclude
that the estimate (3.7) holds for $u$ with appropriate $\gamma > \gamma_1$.
\qed

\smallskip
\state Corollary 3.5. {\it Let $X$, $h$, $J^*$, $A$, and the constants
$\eps_2$ and $\gamma$ be as in {\sl Lemma 3.4}. Suppose that $J$ is a
continuous almost complex structure on $X$ with $\norm{J-J^*}_{L^\infty(A)}
<\eps_2$ and $u \in C^0\cap L^{1,2}(Z(0,l),X)$ a $J$-holomorphic map, such
that $u(Z)\subset A$ and $\norm{du}_{ L^2(Z_i)}<\eps_2$ for any $i=1,\ldots
,l$. Let $\lambda>1$ be (the uniquely defined) real number with
$\lambda = {\gamma \over 2} (\lambda^2+ 1)$.

Then for $2\le k\le l-1$ one has
$$
\norm{du}^2_{L^2(Z_k)} \le \lambda^{-(k-2)} \cdot \norm{du}^2_{L^2(Z_2)}
+ \lambda^{-(l-1-k)} \cdot \norm{du}^2_{L^2(Z_{n-1})}.
\eqno(3.9)
$$
}

\state Proof. The definition of $\lambda$ implies that for any $a_+$ and
$a_-$ the sequence $y_k \deff a_+ \lambda^k + a_- \lambda^{-k}$ satisfies the
recurrent relation $y_k = {\gamma\over2}(y_{k-1} + y_{k+1})$. In particular,
so does the sequence
$$
A_k \deff
\msmall{ \lambda^{-(k-2)} - \lambda^{6-2l+ k-2}
\over 1- \lambda^{6-2l} } \norm{du}^2_{L^2(Z_2)}
+ \msmall{ \lambda^{-(l-1-k)} - \lambda^{6-2l+ l-1-k}
\over 1- \lambda^{6-2l} } \norm{du}^2_{L^2(Z_{l-1})},
$$
which is detemined by the values $A_2 = \norm{du}^2_{L^2(Z_2)}$ and $A_{l-1}
= \norm{du}^2_{L^2(Z_{l-1})}$.

We claim that for $2\le k\le l-1$ one has the estimate $\norm{du}^2 _{L^2(
Z_k)} \le A_k$, which is obviously stronger than (3.9). Suppose that there
exists a $k_0$, such that $2\le k_0\le l-1$ and
$\norm{du}^2 _{L^2(Z_{k_0})}>A_{k_0}$.
Choose $k_0$ so that the difference $\norm{du }^2_{L^2(Z_{k_0})}-A_{k_0}$
is maximal. By {\sl Lemma 3.3} and by our recurrent definition of $A_k$ we
have that $2< k_0 < l-1$ and
$$
\norm{du }^2_{L^2(Z_{k_0})}-A_{k_0}\le {\gamma \over 2}(
\norm{du }^2_{L^2(Z_{k_0+1})}-A_{k_0+1}+
\norm{du }^2_{L^2(Z_{k_0-1})}-A_{k_0-1})\le
$$
$$
\le {\gamma \over 2}2(\norm{du }^2_{L^2(Z_{k_0})}-A_{k_0})
$$
The second inequality follows from the fact that $\norm{du}^2_{L^2(
Z_{k_0})} - A_{k_0}$ is maximal. This gives a contradiction. \qed

\smallskip
An immediate corollary of this estimate is the following improvement of
Sacks-Uhlenbeck theorem about removability of a point singularity, see [S-U]
and [G].

\smallskip
\state Corollary 3.6. {\sl (Removal of point singularities). \it Let $X$ be a
manifold with a Riemannian metric $h$, $J$ a continuous almost complex
structure, and $u:(\check\Delta, J\st)\to (X,J)$ a pseudoholomorphic map
from the punctured disk. Suppose that

\item\sli $J$ is uniformly continuous on $A \deff u(\check\Delta)$ \wrt
$h$ and the closure of $A$ is $h$-complete;

\item \slii there exists $i_0$, such that for all annuli $R_i\deff\{ z\in \cc
:{1\over e^{i+1}}\le | z| \le {1\over e^i}\} $ with $i\ge i_0$ one has
$\norm{du}^2_{L^2(R_i)}\le \eps_2$, where $\eps_2$ is defined in 
{\sl Lemma 3.3}.

\smallskip\noindent
Then $u$ extends to the origin.
}

\smallskip
Condition \sli is automatically satisfied if $A= u(\check\Delta)$ is
relatively compact in $X$. Condition \slii of ``slow growth'' is clearly weaker
than just the boundedness of the area, see e.g.\ [S-U], [G]. It is sufficient
to have $\lim_{i\lrar\infty} \norm{du}^2_{L^2(R_i)} =0$, whereas boundedness
of the area means $\sum_{i=1}^\infty \norm{du}^2_{L^2(R_i)} <\infty$.

\state Proof. The exponential map $\exp(t,\theta) \deff e^{-t+i\theta}$
defines a biholomorphism between the infinite cylinder $Z(0,\infty)$ and
the punctured disk $\check\Delta$, identifying every annulus $R_i$ with
the cylinder $Z(i,i+1)$. Applying {\sl Corollary 3.5} to the map
$u \scirc \exp$ on cylinders $Z(i_0,l)$ and setting $l\lrar \infty$, we get
the estimate
$$
\norm{du}^2_{L^2(R_i)} \le \lambda^{-(i-i_0)} \cdot
\norm{du}^2_{L^2(R_{i_0})}, \qquad i>i_0.
$$
Using this and {\sl Lemma 3.1} we conclude that $\diam(u(R_i)) \le C\cdot
\lambda^{-i/2}$ for $i>i_0$. Since $\sum \lambda^{-i/2} <\infty$, $u$ extends
continuously into $0\in \Delta$. \qed

\medskip
In the proof of the compactness theorem we shall use the following corollary
from {\sl Lemma 3.3}. Let $X$ be a manifold with a Riemannian metric $h$, $J$
a continuous almost complex structure on $X$, $A\subset X$ a closed 
$h$-complete
subset, such that $J$ is $h$-uniformly continuous on $A$. Furthermore, let
$\{J_n\}$ be a sequence of almost complex structures uniformly converging to
$J$, $\{l_n\}$ a sequence of integers with $l_n\to \infty$, and $u_n:Z(0,l_n)
\to X$ a sequence of $J_n$-holomorphic maps.

\medskip
\state Lemma 3.7. {\it Suppose that $u_n(Z(0,l_n))\subset A$ and $\norm{du_n}
_{L^2 (Z_i)}\le \eps_2$ for all $n$ and $i\le l_n$. Take a sequence $k_n\to
\infty$ such that $k_n<l_n-k_n\to \infty$. Then:

\noindent
{\sl1)} $\norm{du_n}_{L^2(Z(k_n,l_n-k_n))}\to 0$ and $\diam\bigl(u_n
(Z(k_n,l_n-k_n))\bigr)\to 0$;

\noindent
{\sl2)} if, in addition, all images $u_n(Z(0,l_n))$ are contained in some bounded
subset of $X$, then there is a subsequence $\{ u_n \}$, still denoted $\{
u_n \}$, such that both $u_n |_{Z(0,k_n)}$ and $u_n |_{Z(k_n,l_n)}$ converge
in $L^{1,p}$-topology on compact subsets in $ \check \Delta \cong Z(0, +\infty )$
to $J$-holomorphic maps $u^+_\infty: \check \Delta \to X$ and $u^-_\infty:
\check \Delta \to X$. Moreover, both $u^+_\infty$ and $u^-_\infty$ extend to
the origin and $u^+_\infty(0)= u^-_\infty(0)$.
}

\state Remarks.~\bf1. The punctured disk $\check\Delta$ with the standard
structure $J_{\Delta } {\d \over \d r}={1\over r}{\d \over \d \theta }$ is
isomorphic to $Z(0,\infty )$ with the structure $J_Z{\d \over \d t}=-{\d \over
\d \theta }$ under a biholomorphism $(\theta ,t)\mapsto e^{-t+ \isl\theta}$.
Thus statement (2) of this corollary is meaningful.

\noindent
{\bf 2. \sl Lemma 3.6} describes explicitly how the sequence of
$J_n$-holomorphic maps of the cylinders of growing conformal radii
converges to a $J$-holomorphic map of the standard node.

\smallskip
\state Lemma 3.8. {\it There is an $\eps_3 =\eps_3(\mu_{J_\infty},A,h)$ such
that for any continuous almost-complex structure $J$ on $X$ with $\norm{J -
J_\infty}_{L^\infty}\le \eps_3 $ and any non-constant $J$-ho\-lo\-mor\-phic
sphere $u: \cc\pp^1 \to X$, $u(\cc\pp^1)\subset A$ one has the inequalities
$$
\area(u(\cc\pp^1 ))\ge \eps_3
\qquad\text{and}\qquad
 \diam(u(\cc\pp^1))\ge \eps_3.
$$}

\state Proof. Let $\eps_1$ be the constant from {\sl Lemma 3.1}. Suppose that
$\area u(\cc\pp^1)=\norm{du}^2_{L^2(\cc\pp^2)}\le \eps_1^2$. Cover $\cc\pp^1$
by two disks $\Delta_1$ and $\Delta_2$. By (3.1) and Sobolev imbedding
$L^{1,p} \subset C^{0,1-{2\over p}}$ we obtain that $\diam(u(\Delta_1))$ and
$\diam(u (\Delta_2) )$ are smaller than $const\cdot \eps_1$. Thus the diameter
of the image of the sphere is smaller than $const\cdot \eps_1$.

So we can suppose that the image $u(S^2)$ is contained in the coordinate
chart, \ie in a subdomain in $\cc^n$, and the structures $J$ and $J_\infty$ are
$L^\infty$-close to a standard one. Consider now $u:S^2\to U\subset \cc^n$
as a solution of the linear equation
$$
\d_xv(z) + J(u(z))\cdot\d_yv(z) = 0\eqno(3.10)
$$
on the sphere. The operator $\dbar _J(v) = d_xv(z) + J(u(z))\cdot
\d_yv(z)$ acts from $L^{1,p}(S^2, \cc^n)$ to $L^p(S^2, \cc^n)$ and is a small
perturbation of the standard $\dbar$-operator. Note that the standard $\dbar$
is surjective and Fredholm. Thus small perturbations are also surjective
and Fredholm, having the kernel of the same dimension. But the kernel of
$\dbar$ consists of constant functions. Since all constants are in the
kernel of (3.10), our $u$ should be a constant map.

We have proved that if the area or a diameter of $J$-holomorphic map is
sufficiently small then this map is constant.\qed

\smallskip
\state Remark. The same statement is true for the curves of arbitrary genus
$g$. In that case, in addition to the estimate (3.1), one should also use the
estimate (3.7). This yields the existence of an $\eps $ which depends on $g$
(and, of course, on $X$, $J$, and $K$), but not on the complex structure on the
parameterizing surface.

\bigskip\bigskip
\newsection{}{Compactness for curves with free boundary}
 
\baselineskip=12.3pt plus .7pt

\medskip
In this section we give a proof of the Gromov compactness theorem for the 
curves with boundaries of fixed finite topological type and without boundary
conditions on maps. The case of closed curves is obviously included in this
one.

Throughout this section we assume that the following setting holds.

\sl Let $X$ be a manifold with a Riemannian metric $h$, $J_\infty$ a continuous
almost complex structure on $X$, $A\subset X$ an $h$-complete subset, $\{ C_n
\}$ a sequence of nodal curves parametrized by a real surface $\Sigma$ with
parametrizations $\delta_n:\Sigma \to C_n$, and $u_n:(C_n,j_n) \to (X, J_n)$
a sequence of pseudoholomorphic maps. Further, $J_\infty$ is $h$-uniformly
continuous on $A$, $J_n$ are also continuous and converge to $J_\infty$,
$h$-uniformly on $A$, $u_n(C_n)\subset A$ for all $n$. \rm

\smallskip
Let us explain the main idea of the proof of {\sl Theorem 1.1}. The Gromov
topology on the space of stable curves over $X$ is introduced in order to
recover breaking of ``strong" (\ie $L^{1,p}$-type) convergence of a sequence
$(C_n,u_n)$ of pseudoholomorphic curves of bounded area. The are two reasons
for this. The first one is that a sequence of (say, smooth) curves $C_n$ could
diverge in an appropriate moduli space and the second one is a phenomenon of
``bubbling". In both cases one has to do with appearance of new nodes, \ie
with certain degeneration of complex structure on curves. The ``model''
situation of {\sl Lemma 3.6} describes a convergence of ``long cylinders''
$u_n :Z(0, l_n) \to X$, $l_n \lrar \infty$, to a node $u_\infty: \cala_0 \to
X$. In our proof we cover curves $C_n$ by pieces, which are either ``long
cylinders'' converging to nodes or have the property that complex structures
and maps ``strongly'' converge. Here the ``strong'' convergence means the
usual one, \ie \wrt the $C^\infty$-topology for complex structures, and \wrt
the $L^{1,p}$-topology with some $p>2$ for maps. In fact, the strong
convergence of maps is equivalent to the uniform one, \ie \wrt the
$C^0$-topology, and implies further regularity in the case when $J_n$
and $J_\infty$ have more smoothness. One consequence of this is that we remain
is the category of nodal curves. Another one is that we treat degeneration
of complex structure on $C_n$ and the
``bubbling'' phenomenon  in a uniform framework of ``long cylinders''.

\smallskip
For the proof we need some additional results.

\state Lemma 4.1. \it For any $R>1$ there exists an $a^+= a^+(R)>0$ with the
following property. For any cylinder $Z = Z(0,l)$ with $0 <l\le +\infty$
and any annulus $A \subset Z(0,l)$, which is adjacent to $\d_0 Z = S^1\times
\{0\}$ and has conformal radius $R$, one has $Z(0,a^+) \subset A$.
\rm

\state Proof. Without loss of generality we may assume that $l=+\infty$ and
identify $Z$ with the punctured disk $\check\Delta$ via the exponential map
$(-t+ \isl \theta) \mapsto e^{ -t+ \isl \theta}$, such that $\d_0 Z$ is
mapped onto $S^1= \d\Delta$.

Suppose that the statement is false. Then there would exist holomorphic
imbeddings $f_n : A(1, R) \to \check\Delta$ and points $a_n \in \Delta \bs
f_n(A(1, R))$, such that $f_n(A(1, R))$ are adjacent to $\d\Delta$ and $a_n
\lrar a\in \d\Delta$. Passing to a subsequence, we may assume that $\{f_n\}$
converges uniformly on compact subsets in $A(1, R)$ to a holomorphic map
$f: A(1, R) \to\Delta$.

If $f$ is not constant, then $f(A(1, R))$ must contain some annulus $\{ b
<|z|<1 \}$ with $b<1$. But then $\{ \sqrt{b} <|z|<1\} \subset f_n(A(1, R))$
for $n>\!>1$, which is a contradiction.

If $f$ is constant, then the diameter of images of the middle circle
$\gamma\deff \{ |z|= \sqrt{R} \} \subset A(1, R)$ must converge to $0$. But
$\diam(f_n(\gamma)) \ge \dist(0, a_n) \lrar1$. The obtained contradiction
finishes the proof. \qed

\medskip
For the proof of {\sl Theorem 1.1} we need a special covering of\/ $\Sigma$
which will be constructed in the following theorem.

\state Theorem 4.2. {\it Under the conditions of {\sl Theorem 1.1}, after
passing to a subsequence, there exist a finite covering $\calv$ of\/ $\Sigma$
by open sets $V_\alpha$ and parametrizations $\sigma_n:\Sigma \to C_n$
such that:

{\sl(a)} all $V_\alpha$ are either disks, or annuli, or pants;

{\sl(b)} for any boundary circle $\gamma_i$ of\/ $\Sigma$ there is some
annulus $V_\alpha$ adjacent to  $\gamma_i$;

{\sl(c)} $\sigma_n^*j_n\ogran_{V_\alpha}$ doesn't depend on $n$ if
$V_\alpha$ is a disk, pants, or an annulus adjacent to a boundary circle
of\/ $\Sigma$;

{\sl(d)} all non-empty intersections $V_\alpha \cap V_\beta$ are annuli,
where the structures $\sigma^*j_n$ are independent of $n$;

{\sl(e)} if $a$ is a node of $C_n$ and $\gamma^n_a=\sigma_n\inv (a)$ the
corresponding circle, then $\gamma_a^n=\gamma_a$ doesn't depend on $n$, is
contained in some annulus $V_\alpha$, containing only
one such ``contracting" circle for any $n$; moreover, the structures
$\sigma_n^*j_n\ogran_{V_\alpha\bs \gamma_a}$ are independent of $n$;

{\sl(f)} if $V_\alpha$ is an annulus and $\sigma_n( V_\alpha)$ are not
nodes, then the conformal radii of $\sigma_n( V_\alpha)$ converge to some
positive $R_\alpha^\infty>1$ or to $+\infty$.

{\sl(g)} if for initial parametrizations $\delta_n$ and fixed annuli
$A_i$, each adjacent to the boundary circle $\gamma_i$ of\/ $\Sigma$,
the structures $\delta_n^* j_n\ogran_{A_i}$ do not depend on $n$ $A_i$, then
the new parametrizations $\sigma_{n_k}$ can be taken equal to $\delta_{n_k}$
on some subannuli $A'_i \subset A_i$ also adjacent to $\gamma_i$.  }

\state Proof. We shall prove the properties {\sl(a)--(f)}. The property
{\sl(g)} will follow from {\sl Lemma 4.3} below.

There are 4 cases where the existence of such a covering is obvious. If
all $C_n$ are disks or annuli without nodal points, there is nothing to
prove. In the third case each $C_n$ is a sphere, and we cover it by 2
disks.

In the forth case each $C_n$ is a torus without marked points. Then
Any complex torus can be represented in the form $\cc{\bigm/}(\zz+
\tau \zz)$ with $|\re\tau| \le {1\over 2}$ and $\im\tau > {1\over 2}$.
Considering the map $z\in \cc \mapsto e^{2\pi\isl z} \in \check\cc \deff
\cc\bs\{0\}$, we represent $(T^2,j)$ as the quotient $\cc\bigm/ \{ z \sim
\lambda^2 z\}$ with $\lambda = e^{\pi\isl\tau}$, so that $|\lambda|<
e^{-\pi/2}< {1\over 3}$. The annuli $\{ {|\lambda| \over2} <|z| <1\}$ and
$\{{|\lambda|^2 \over2} <|z| <|\lambda| \}$ form the needed covering.

\smallskip
In all remaining cases we start with constructing of appropriate graphs
$\Gamma_n$ associated with some decomposition of $C_n$ into pants. {\sl Lemma
3.4} and non-degeneration of the complex structure $j_n$ on $C_n$ shows that
lengths of all boundary circles of all non-exceptional components $C_{n,i}$ of
$C_n$ are uniformly bounded from above. At this point we make the following

\state Remark. The Collar Lemma from [Ab], Ch.II, \S\.3.3 yields the existence
of the universal constant $l^*$ such that for any simple geodesic circles
$\gamma'$ and $\gamma''$ on $C_{n,i}$ the conditions $\ell(\gamma') <l^*$ and
$\ell( \gamma'') <l^*$ imply $\gamma' \cap \gamma'' = \emptyset$. We shall
call geodesic circles $\gamma$ with $\ell(\gamma) <l^*$ {\sl short geodesics}.

\smallskip
The fact that $(C_n,u_n)$ are pseudoholomorphic and of bounded area shows
that $C_n$ have uniformly bounded number of components. Indeed, the number of
exceptional components, which  are spheres and disks, is bounded by the energy
(see {\sl Lemma 3.8}), whereas the  number of boundary circles of $\Sigma$.
Further, the operation of contracting a circle on $\Sigma$ to a nodal point
either dininish genus of some component of $C_n$, or increase the number of
components. Thus, the number of nodal points on $C_n$ and the total number of
marked points on its components are also uniformly bounded. This implies that
the number of possible topological types of components is  finite.

In this situation the Teichm\"uller theory (see [Ab], Ch.II, \S\.3.3) states
the existence of decomposition of every non-exceptional component $C_{n,i}\bs
\mapo$ into pants with the following properties:

\sli every short geodesic is a boundary circle of some pants of the
decomposition;

\slii the intrinsic length of every boundary circle is bounded from above by
a (uniform) constant depending only on an upper bound of lengths of boundary
circles and possible topological types of $C_{n,i} \bs \mapo$.

Having decomposed all $C_{n,i} \bs \mapo$ into pants, we associate with every
curve $C_n$ its graph $\Gamma_n$. As it was noted above, the number of
vertices and edges of $\Gamma_n$ is uniformly bounded. Thus after passing to
a subsequence, we can assume that all $\Gamma_n$ are isomorphic to each other
(as marked graphs). Denote this graph by $\Gamma$. Now, the parametrizations
$\sigma_n: \Sigma \to C_n$ can be found in such a way that the decompositions
of $C_{n,i} \bs\mapo$ into pants define the same set $\bfgamma =\{ \gamma
_\alpha \}$ of circles on $\Sigma$ and induce the same decomposition $S\bs
\cup_\alpha \gamma_\alpha = \cup_j S_j$ with the graph $\Gamma$.

\smallskip
By our construction of the graph $\Gamma$, each edge of $\Gamma$ corresponds
either to a circle in $\Sigma$ contracted by every parametrization $\sigma_n$
to a nodal point, or to a circle mapped by every $\sigma_n$ onto a geodesic
circle separating two pants. Furthermore, each tail of $\Gamma$ corresponds
to a boundary circle of $\Sigma$. Thus we shall use the same notation $\gamma
_\alpha$ for an edge or a tail of $\Gamma$ and for the corresponding circle
on $\Sigma$. If $\gamma_\alpha$ is a boundary circle of some pants $S_j$, 
then the intrinsic length $\ell_{n,\alpha}= \ell_n (\gamma_\alpha)$ of
$\sigma_n( \gamma_\alpha)$ is well defined. This happens in the following two
cases:

{\sl \.a)}~$\gamma_\alpha$ separates two pants, or else

{\sl \.b)}~$\gamma _\alpha$ is a boundary circle $\Sigma$ and for any $n$
the irreducible component of $C_n$ adjacent to $\sigma_n(\gamma_\alpha)$ is 
not a disk with a single nodal point. Note that appearance of these two cases
is independent of $n$.

By our choice of $\gamma_\alpha$, the lengths $\ell_n (\gamma_\alpha)$ are
uniformly bounded from above. Passing to a subsequence, we may assume that
for any fixed $\alpha$ the sequence $\{ \ell_{n,\alpha} \}$ converges to
$\ell_{\infty,\alpha}$.

\smallskip
As one can expect, the condition $\ell_{n, \alpha} \lrar 0$ means that the
circle $\gamma_\alpha$ is shrunk to a nodal point on the limit curve. We
shall prove the statement of the theorem by induction in the number $N$ of
those circles $\gamma_\alpha$ for which $\ell_{\infty,\alpha} =0$.

\smallskip
The case $N=0$, when there are no such circles, is easy. Passing to a
subsequence, we may assume that the Fenchel-Nielsen coordinates $(\bfell_n,
\bfvartheta_n)$ of any non-exceptional component $C_{n,i}$ of $C$ converge to
the Fenchel-Nielsen coordinates $(\bfell_\infty, \bfvartheta_\infty)$ of some
smooth curve $C_{\infty, i}$ with marked points. Gluing together appropriate
pairs of marked points we obtain a nodal curve $C_\infty$, which admits a
suitable parametrization $\sigma_\infty: \Sigma \to C_\infty$ and has the
same graph $\Gamma$. {\sl Lemma 3.5} shows that for $n>\! >1$ the curves
$C_n$ can be obtained from $C_\infty$ by deformation of the transition
functions for the intrinsic local coordinates on non-exceptional components
of $C_\infty$. Note that such a deformation can be realized as a deformation
of operator $j_\infty$ of complex structure on $C_\infty$, localized in small
neighborhoods of circles $\sigma_\infty(\gamma_\alpha)$, see Fig.~6. In the
case when $\gamma_\alpha$ is a boundary circle, we may additionally assume
that the annulus, where $j_n$ changes, lies away from $\gamma_\alpha$. Now
the existence of the covering with desired properties is obvious.

\medskip
\vbox{\nolineskip\xsize.54\hsize%
\putm[.01][.38]{\underbrace{\hskip.35\xsize}_{V_\beta}}%
\putm[.54][.44]{\underbrace{\hskip.46\xsize}_{V_\gamma}}%
\putm[.27][.06]{\overbrace{\hskip.36\xsize}^{V_\alpha}}%
\putm[.27][.3]{\underbrace{\hskip.09\xsize}_{W_\beta}}%
\putm[.53][.3]{\underbrace{\hskip.09\xsize}_{W_\gamma}}%
\putt[1.05][0]{\advance\hsize-1.05\xsize\parindent=0pt%
\centerline{Fig.~6. }
\smallskip
$V_\alpha$ and $V_\beta$ represent elements of the covering where the complex
structure is constant. The change of complex structure is done in the painted
part of  $V_{\gamma}$. $W_\beta\deff V_\alpha \cap V_\beta$ and $W_\gamma
\allowbreak \deff V_\alpha \cap V_\gamma$ represent the annuli with the
constant complex structure.
}
\noindent
\epsfxsize=\xsize\epsfbox{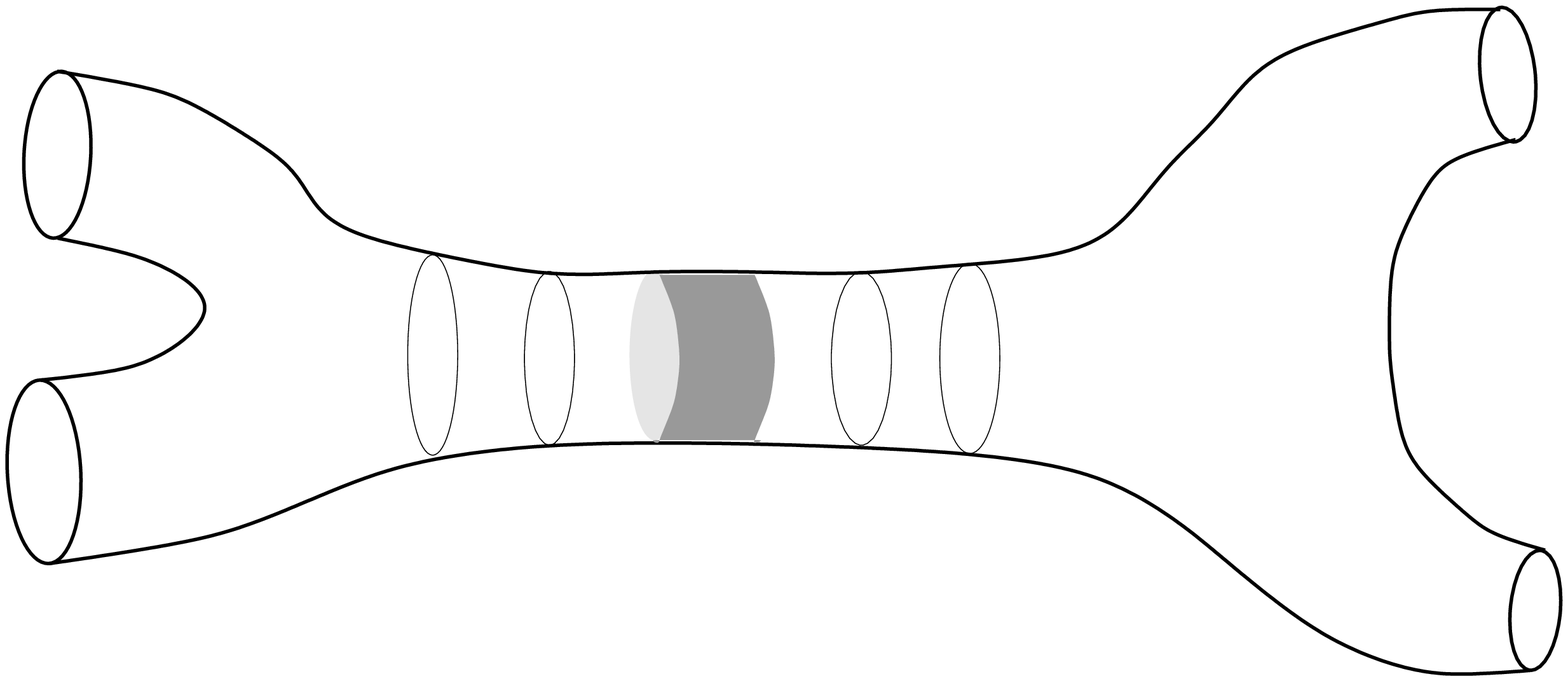}
}
\bigskip\smallskip

\baselineskip=12.3pt plus .2pt

\medskip
Let us consider now the general case when the number $N$ of ``shrinking
circles'' is not zero. Take a circle $\gamma_\alpha$ with 
$\ell_\infty (\gamma _\alpha) =0$. Let $S_j$ be pants adjacent to 
$\gamma_\alpha$. Consider the intrinsic coordinates $\rho_\alpha$ and 
$\theta_\alpha$ at $\sigma_n(\gamma_\alpha)$ and the annuli
$$
\eqalign{
A_{n,\alpha,j} &
\deff \left\{\, (\rho _\alpha, \theta_\alpha) \in \sigma_n(S_j) 
\;:\; 0\le \rho_\alpha \le 
\msmall {\pi^2 \over \ell_{n,\alpha} } - \msmall {2\pi \over a^*} 
\;\right\}
\cr
A^-_{n,\alpha,j} &\deff
\left\{\, (\rho_\alpha, \theta_\alpha) \in \sigma_n(S_j) \;:\;
0\le \rho_\alpha \le \msmall {\pi^2 \over \ell_{n,\alpha} } - 
\msmall {2\pi \over a^*}-1
\;\right\},
}
$$
adjacent to $\sigma_n( \gamma_\alpha)$. Note that ${\pi^2 \over \ell_{n,
\alpha} }- {2\pi \over a^*}$ (resp.\ ${\pi^2 \over \ell_{n,\alpha} }- {2\pi
\over a^*}-1$) is the logarithm of the conformal radius of $A_{n, \alpha, j}$
(resp.\ of $A^-_{n, \alpha,j}$). Consequently, we can use {\sl Lemma 2.2} to
shows that these annuli are well-defined.

If $\gamma_\alpha$ is a boundary circle we set $C^-_n \deff C_n \bs A^-_{n,
\alpha,j}$. Otherwise $\gamma_\alpha$ separates two pants, say $S_j$ and
$S_k$. Then we define in a similar way the annuli $A_{n,\alpha, k} \subset 
\sigma_n(S_k)$ and $A^- _{n,\alpha,k}\subset \sigma_n(S_k)$, 
set $A_{n,\alpha} \deff A_{n, \alpha,j} \cup A_{n, \alpha, k}$ and 
$A^-_{n,\alpha} \deff A^-_{n,\alpha,j} \cup A^-_{n, \alpha, k}$, and
put $C^-_n \deff C_n \bs A^-_{n,\alpha}$.

The parametrizations $\sigma_n: \Sigma \to C_n$ can be chosen in such a way
that the annuli $\sigma_n\inv(A^-_{n,\alpha,j})$ (resp.\ $\sigma_n\inv(A^-
_{n, \alpha,k})$) define the same annulus $A^-_{ \alpha,j}$ (resp.\ $A^-
_{\alpha, k}$) on $\Sigma$. Let $\gamma^-_{ \alpha,j}$ (resp.\ $\gamma^-
_{\alpha,k}$) denote its boundary circles different from $\gamma_\alpha$.
Thus, the curves $C^-_n$ are parametrized by a real surface $\Sigma^- \deff
\Sigma \bs A^-_{ \alpha,j}$ (resp.\ $\Sigma^- \deff \Sigma \bs (A^-_{\alpha,
j} \cup A^- _{\alpha, k}$), and the restrictions of $\sigma_n$ can be chosen
as parametrization maps. Thus the decompositions of components of $C_n$ into
pants define the combinatorial type of decompositions of components of
$C^-_n$ into pants. Moreover, the corresponding graph $\Gamma^-$ will be the
same for all $C^-_n$. It coincides with $\Gamma$ if $\gamma_\alpha$ is a
boundary circle. Otherwise $\Gamma^-$ can be obtained from $\Gamma$ by
replacing the edge corresponding to $\gamma_\alpha$ by 2 new tails for 2 new
boundary components.

Let $C_{n,i}$ be the component of $C_n$ adjacent to $\sigma_n(\gamma _\alpha)$
and $C^-_{n,i}$ the corresponding componet of $C^-_n$. Decompose $C^-_{n,i}$ 
into pants according to graph $\Gamma^-$ in the canonical way, so that the
boundary circles of the obtained pants are geodesic. Note that even if the
constructed pants are in combinatorial one-to-one correspondence with the
pants of $C_n$, the intrinsic metric on $C^-_{n,i}$ and the obtained geodesics
circle $\gamma^-_\beta$ differs from the corresponding objects on $C_{n,i}$

Nevertheless, we claim that for the obtained decomposition of $C^-_{n,i}$ the
intrinsic lengths are uniformly bounded from above (possibly by a new
constant) and that the sequence $\{ C^-_n \}$ has less ``shrinking circles'' 
than $\{ C_n \}$. This explains the  meaning of the above construction, when 
the curves $C^-_n$ are obtained from $C^-_n$ by cutting off the annuli $A^-
_{n,\alpha,j}$ (and resp.\ the annuli $A^-_{n,\alpha,j}$). The annuli we choose
are sufficiently long so that one ``shrinking circle'' disappears, but not
too long so that the complex structures of the curves $C^-_n$ remain
non-degenerating near the boundary.

\smallskip
Indeed, the complex structures on $C^-_n$ do not degenerate at the boundary
circle $\gamma^-_{ \alpha,j}$ (resp.\ at $\gamma^-_{\alpha,k}$), since
$C^-_n$ contain annuli $A_{n,\alpha,j} \bs A^-_{n,\alpha,j}$ (resp.\
$A_{n,\alpha,k} \bs A^- _{n,\alpha,k}$) of the constant conformal radius
$R=e>1$. This implies that the lengths $\ell^-_n(\gamma^-_{\alpha,j})$ of
$\sigma_n( \gamma^- _{\alpha,j})$ (resp.\ $\ell^-_n(\gamma^-_{\alpha,k})$ of
$\sigma_n(\gamma^- _{\alpha,k})$) with respect to the intrinsic metrics on
$C^-_n$ are uniformly bounded.

On the other hand, the lengths $\ell^-_n(\gamma^-_{\alpha,j})$ (resp.\
$\ell^-_n(\gamma^-_{\alpha,k})$) are also uniformly bounded from below by a
positive constant. Otherwise, by {\sl Lemma 3.4}, after passing to a
subsequence, there would exist annuli $A_n\subset C^-_n$ of infinitely
increasing radii $R_n$, adjacent to $\sigma_n(\gamma^-_{\alpha,j})$ (resp.\
to $\sigma_n(\gamma^-_{ \alpha,k})$). The superadditivity of the logarithm of
the conformal radius of annuli, see [Ab], Ch.II, \S\.1.3, shows that the
conformal radius $R^+_n$ of the annulus $A^-_{n,\alpha,j} \cup A_n$ satisfies
the inequality $\log R^+_n \ge {\pi^2 \over \ell_{n,\alpha} }- {2\pi \over
a^*}-1 + \log R_n$, which contradicts {\sl Lemma 2.2}, part {\sl i)}.

\smallskip
Now we estimate the intrinsic lengths of boundary circles and the number of
``shrinking circles'' on $C^-_{n,i}$. 
Compute the width $L_n$ of $A_{n,\alpha,j} \bs A^- _{n, \alpha,j}$ \wrt the 
intrinsic metric on $C_n$. Using $\ell_{n, \alpha} \lrar0$, we get
$$
\eqalign{
L_n &= \int_{\rho={\pi^2 \over \ell_{n, \alpha} }- {2\pi \over a^*}-1}
^{ {\pi^2 \over \ell_{n, \alpha} }- {2\pi \over a^*}}
\left(\msmall{ {\ell_{n, \alpha} \over 2\pi} \over \cos 
{\ell_{n, \alpha} \rho \over 2\pi}} \right) d\rho 
=
\left[ \log\;\cotan \left( \msmall{\pi \over 4} 
-\msmall{\ell_{n, \alpha} \rho \over 4\pi} \right)
\right]_{\rho={\pi^2 \over \ell_{n, \alpha} }- {2\pi \over a^*}-1}
^{{\pi^2 \over \ell_{n, \alpha} }- {2\pi \over a^*}}
\cr
&= \log \msmall{ \cotan{\ell_{n, \alpha}  \over2a^*} \over 
\cotan({\ell_{n, \alpha} \over 2a^*} + {\ell_{n, \alpha} \over4\pi}) }
\approx \log \msmall{ {\ell_{n, \alpha}  \over2a^*} + 
{\ell_{n, \alpha} \over4\pi} \over
{\ell_{n,\alpha} \over2a^* } } =
\log \left( 1 + \msmall{ a^* \over 2\pi } \right) >0.
}
$$
Let $\gamma_\beta\not= \gamma_\alpha$ be another boundary circle of $\Sigma$,
such that the circles $\sigma_n(\gamma_\beta)$ bound $C_{n,i}$.
Denote by $\ell_{n, \beta}$ (resp.\ by $\ell^-_{n, \beta}$) the length of 
$\sigma_n(\gamma_\beta)$ \wrt the intrinsic length of $C_{n,i}$ (resp.\ 
of $C^-_{n,i}$). By {\sl Lemma 2.2} for any $n$ we can find an annulus 
$A_{n, \beta} \subset C_n$ of the constant width which is adjacent to 
$\sigma_n(\gamma_\beta)$ and has the conformal radius $R_{n, \beta}$ with 
$\log R_{n, \beta}= {\pi^2\over  \ell_{n, \beta}}-{2\pi \over a^*}$. The width
of $A_{n, \beta}$ is then 
$$
L_{n, \beta} = \log\, \cotan \left( \msmall{\pi \over 4} -
\msmall{\ell_{n, \beta}\over 4\pi}
\left( \msmall{\pi^2\over  \ell_{n, \beta}}- {2\pi \over a^*}\right)
 \right)=  
\log\, \cotan \left( \msmall{\ell_{n, \beta}\over 2a^*}\right).
$$
Since $\ell_{n, \beta}$ are uniformly bounded from above, $L_{n, \beta}$ are 
uniformly bounded from below. This implies that we can find a subannuli $A^-
_{n, \beta} \subset C^-_{n,i}$, which are ajacent to $\sigma_n(\gamma_\beta)$ 
and have the constant width $L^*>0$ \wrt the intrinsic metric of $C_{n,i}$.
Computing the conformal radius $R_{n, \beta}$ of such an annulus $A^- _{n, 
\beta}$ we get the relation 
$$
L^* = \log\,\cotan\left (\msmall {\pi\over4} -
 \msmall {R_{n, \beta} \ell_{n, \beta} \over 4\pi} 
\right).
$$
Conseqently, $R_{n, \beta}={M \over \ell_{n, \beta}}$ for some constant $M>0$.
Since $\ell_{n, \beta}$ are uniformly bounded from above, $R_{n, \beta}$
are uniformly bounded from below. This means that the complex structures of 
$C^-_n$ do not degenerate near boundary.

\smallskip
Now we estimate the number of ``shrinking circles'' on $C^-_n$. Take a circle
$\gamma\subset \Sigma^-$ such that for every $n$ there exists a simple closed 
geodesic $\gamma^-_n \subset C^-_n$ homotopic to $\sigma^-_n(\gamma)$. Assume 
that the intrinsic lengths $\ell^-(\gamma^-_n)$ vanish. Using {\sl Lemma 2.2},
we can construct annuli $A_n \subset C^-_n$ homotopic to $\gamma^-_n$ whose
conformal radii $R_n$ increase infinitely. But then there exist simple closed 
geodesics $\gamma_n \subset C_n$ homotopic to $\sigma_n(\gamma)$, and the
$A_n \subset C_n$ are homotopic to $\gamma_n$. {\sl Lemma 2.2} implies that
the intrinsic lengths $\ell(\gamma_n)$ also vanish. Thus every ``shrinking 
circle'' on $C^-_n$ appears from some ``shrinking circle'' on $C_n$.

\smallskip
Thus we have shown, that $\{ C^-_n \}$ with the intrinsic metric and the 
defined by decomposition $\Gamma^-$ have uniformly bounded lengths of marked 
circles and less ``shrinking circles'' than $\{ C_n \}$. By induction, we 
may assume that there exist the covering of $\Sigma^-$ and parametrizations of 
$C^-_n$ by $\Sigma^-$ with the desired properties. Since $C_n \bs C^-_n$ are 
annuli of increasing conformal radii, the statement of the theorem is valid 
for $C_n$. \qed

\state Lemma 4.3. \it Let $C_n$ be a sequence of annuli with complex 
structures $j_n$, $\Sigma$ some fixed annulus, and $\delta_n: \Sigma \to C_n$
some parametrizations. Suppose that for two fixed annuli $A_1, A_2 \subset
\Sigma$ adjacent to the boundary circles of $\Sigma$ the restrictions
$\delta_n^*j_n|_{A_i}$ do not depend on $n$.

Then one can find parametrizations $\sigma_n: \Sigma \to C_n$ such that
$\sigma_n$ coincide with $\delta_n$ on some (possibly smaller) annuli $A'_i$,
also adjacent to the boundary circles of $\Sigma$ restrictions, and such that:

\smallskip
\sli if conformal radii $R_n$ of $C_n$ converge to $R_\infty <\infty$, then
$\sigma_n^*j_n$ converge to some complex structure;

\smallskip
\slii if conformal radii $R_n$ of $C_n$ converge to $\infty$, then for some
circle $\gamma \subset \Sigma$ structures $\sigma_n^*j_n$ converge on
compact subsets $K \Subset \Sigma \bs \gamma$ to the complex structure of 
the disjoint union of two punctured disks. Moreover, as such $\gamma$ an 
arbitrary imbedded circle generating $\pi_1(\Sigma)$ can be chosen. \rm

\state Proof. Without loss of generality we may assume that $\Sigma = A(1,
10)$, $A_1=A(7,10)$, $A_2=A(1,4)$ and that the given circle $\gamma$ lies in
$A(3,7)$. Let $\delta_n : \Sigma \to C_n$ be the given parametrizations.
There exist biholomorphisms $\phi_n: C_n \to A(r_n,1)$ with $r_n\inv=R_n$
being the conformal radii of $C_n$, such that $\phi_n (\delta_n(A_1))$ is
adjacent to $\{ |z|=1\} = \d\Delta$ and $\phi_n (\delta_n(A_2))$ is adjacent
to $\{ |z|=r_n\}$. Define $\phi'_n(z) \deff {r_n \over \phi_n(z)}$.

Recall that structures $\delta_n^*j_n|_{A_i}$ are the same for all $n$. We call
this structure $j$. Consider the maps $\phi_n \scirc \delta_n : (A_1, j) \to
A(r_n,1) \subset \Delta$ and $\phi'_n \scirc \delta_n : (A_2, j) \to A(r_n,1)
\subset \Delta$. Passing to a subsequence we can suppose that $r_n \lrar
r_\infty <1$ and that the maps $\phi_n \scirc \delta_n$, $\phi'_n \scirc 
\delta_n$
converge on $A_1\cup A_2$ to holomorphic maps $\psi: (A_1\cup A_2, j) \to
\Delta$ and $\psi': (A_1\cup A_2, j) \to \Delta$ respectively. This means that
the maps $(\phi_n \scirc \delta_n, \phi'_n \scirc \delta_n): (A_1 \cup A_2)
\to \Delta^2$ take values in $\{(z,z') \in \Delta^2 \;:\; z{\cdot} z'= r_n
\}$ and converge to the map $(\psi, \psi'): (A_1 \cup A_2) \to \Delta^2$ with
values in $\{(z,z') \in \Delta^2 \;:\; z{\cdot}z' = r_\infty \}$.

The arguments from the proof of {\sl Lemma 4.1} show that the annuli $\psi(
A_1)$ and $\psi'(A_2)$ are adjacent to $\d\Delta$. This implies that for $n>
\!>1$ there exist diffeomorphisms $(\psi_n, \psi'_n): \Sigma \to \{(z,z') \in
\Delta^2 \;:\; z{\cdot}z' = r_n \}$ such that $\psi_n \equiv \phi_n \scirc
\delta_n$ on $A(9,10)$, $\psi'_n(z) \equiv \phi'_n \scirc \delta_n(z)$ for
$z\in A(1,2)$, and $(\psi_n, \psi'_n)$ converge to $(\psi, \psi')$ on
$\Sigma$. Moreover, we may assume that $|\psi_n(t)|=|\psi'_n(t)|= \sqrt{r_n}$
for any $t\in \gamma$. This means that $(\psi_n, \psi'_n)(\gamma)$ lies on
the middle circle $\{ (z,z'): |z|= \sqrt{r_n},\; z'= {r_n\over z} \}$ of
$\{(z,z') \in \Delta^2 \;:\; z{\cdot}z' = r_n \}$.

Set $\sigma_n \deff \phi_n\inv \scirc \psi_n: \Sigma \to C_n$. Then, obviously,
$\sigma_n \equiv \delta_n$ on $A(1,2)$ and on $A(9,10)$, and $\sigma_n^*j_n=
(\psi_n,\psi')_n^* J\st \lrar (\psi,\psi')^* J\st$, where $J\st$
denotes the standard complex structure on $\Delta^2$.
\qed

\medskip
\state Proof of Theorem 1.1. Let $\{(C_n, u_n)\}$ be the sequence from the
hypothesis of the theorem. Then the condition {\sl c)\/} of the theorem and
{\sl Lemma 4.1} yields the existence of parametrizations $\delta_n: \Sigma 
\to C_n$ and annuli $A_i$, adjacent to each boundary circle $\gamma_i$, such
that $\delta_n^* j_{C_n}$ are constant in every $A_i$. Thus we may assume that
$\delta_n$ with these properties are given.

Take a covering $\calv= \{V_\alpha\}$ and parametrizations $\sigma_n$ as in
{\sl Theorem 4.3}. With every such covering we can associate the curves
$C_{\alpha,n} \deff \sigma_n(V_\alpha)$, the parametrizations $\sigma_{\alpha,
n} \deff \sigma_n\ogran_{V_\alpha}: V_\alpha \to C_{\alpha,n}$, and the maps
$u_{\alpha,n} \deff u_n \ogran_{C_{\alpha, n}} : C_{\alpha, n} \to X$. Consider
the following type of convergence of sequences $\{ (C_{\alpha, n}, u_{\alpha,
n}, \sigma_{\alpha, n} )\}$ with $\alpha$ fixed:

\sl
\item{A)} $C_{\alpha, n}$ are annuli of infinitely growing conformal radii
$l_n$ and the conclusions of {\sl Lemma 3.7} hold;

\item{B)} every $C_{\alpha, n}$ is isomorphic to the standard node $\cala_0=
\Delta \cup_{ \{0\} } \Delta$, such that the compositions $V_\alpha \buildrel
\sigma _{\alpha, n} \over \lrar C_{\alpha, n} \buildrel \cong \over \lrar
\cala_0$ define the same parametrizations of $\cala_0$ for all $n$;
furthermore, the induced maps $\ti u_{\alpha, n}: \cala_0 \to X$ strongly
converge;

\item{C)} the structures $\sigma_n^*j\vph_n \ogran_{V_\alpha}$ and the maps
$u_{\alpha, n}\scirc \sigma _{\alpha, n}: V_\alpha \to X$ strongly converge.

\rm\noindent
Here the strong convergence of maps is the one in the $L^{1,p}$-topology on
compact subsets for some $p>2$ (and hence for all $p<\infty$), and the
convergence of structures means the usual $C^\infty $-convergence.

\smallskip
Suppose that there is a subsequence, still indexed by $n\to \infty$, such
that for any $V_\alpha$ we have one of the convergence types {\sl A)--C)}.
Then the sequence of global maps $\{ (C_n , u_n, \sigma _n )\}$ converges in
the Gromov topology which gives us the proof, and also a precise description 
of the convergence.

Otherwise, we want to find a refinement of our covering $\calv$ and
parametrizations $\sigma_n$ which have the needed properties. We shall proceed
by induction estimating the area of pieces of coverings of $\Sigma$. To do so,
we fix $\eps>0$ satisfying $\eps \le {\eps_1^2 \over2 }$ with $\eps_1$ from 
{\sl Lemma 3.1}, $\eps \le {\eps_2^2 \over2 }$ with $\eps_2$ from {\sl Lemma 
3.3}, and $\eps \le {\eps_3\over 3}$ with $\eps_3$ from {\sl Lemma 3.8}. 
Consider first

\smallskip\noindent{\sl
Special case: $\area(u_n(\sigma_n(V_\alpha))) \le \eps$ for any $n$ and any
$V_\alpha \in \calv$}. We can consider every $V_\alpha$ separately. If the
structures $\sigma_n^* j_n |_{V_\alpha}$ are constant, then, due to {\sl 
Corollary 3.3}, some subsequence of $u_n \scirc \sigma_n$ strongly converges.

If structures $\sigma_n^* j_n |_{V_\alpha}$ are not constant, then $V_\alpha$
must be an annulus. Fix biholomorphisms $\phi_n: Z(0,l_n) \buildrel \cong
\over \lrar \sigma_n(V_\alpha)$. If $l_n \lrar \infty$, then
{\sl Lemma 3.7} shows
that (and describes how!) an appropriate subsequence of $u_n \scirc \phi_n$
converges to a $J_\infty$-holomorphic map of a standard node. Otherwise we can
find a subsequence, still denoted $(C_n, u_n)$, for which $l_n \lrar l_\infty
<\infty$ and $u_n \scirc \phi_n$ converge to a $J_\infty$-holomorphic map of
$Z(0,l_\infty)$ in $L^{1,p}$-topology on compact subsets $K \Subset Z(0,l_\infty)$
for any $p<\infty$. To construct refined parametrizations $\ti\sigma _{\alpha,
n} : V_\alpha \to C_{\alpha, n}= \sigma_n(V_\alpha)$, we use property
{\sl(d)} of {\sl Theorem 4.2} and apply {\sl Lemma 4.3}.

Thus we get one of the convergence types {\sl A)--C)} which completes the 
proof in {\sl Special case}.

\smallskip\noindent
{\sl General case.}
Suppose that the theorem is proved for all sequences of $J_n$-holo\-mor\-phic
curves $\{(C_n, u_n)\}$ with parametrizations $\delta_n: \Sigma \to C_n$
which satisfy the additional condition $\area(u_n(C_n)) \le (N-1) \eps$ for
all $n$. We consider this as the induction hypothesis in $N$, so that our
{\sl Special case} is the base of the induction.

\smallskip
Assume that there exists a subsequence, still indexed by $n\to \infty$, such
that for every $V_\alpha$ and for the curves $C_{\alpha, n} = \sigma_n(
V_\alpha )$ the statement of the theorem holds. This means the existence of
refined coverings $V_\alpha= \cup_i V _{\alpha, i}$ and new parametrizations
$\ti\sigma _{\alpha, n}: V_\alpha \to C_{\alpha, n}$, such that $\ti\sigma_n$
coincide with $\sigma_n$ near the boundary of every $V_\alpha$ and such that
for curves $C_{\alpha, i,n} \deff \ti\sigma_n (V _{\alpha ,i})$ we have the
convergence of one of the types {\sl A)--C)}. Then we can glue $\ti\sigma
_{\alpha, n}$ together to global parametrizations $\ti\sigma_n: \Sigma \to
C_n$ and set $\wt\calv \deff \{ V _{\alpha, i} \}$, getting the proof.

In particular, it is so for any $V_\alpha$, such that $\area(u_n(\sigma_n(
V_\alpha))) \le (N-1) \eps$ for all $n$ due to inductive hypothesis.

This implies that it is sufficient to consider only those $V_\alpha$,
for which $(N-1) \eps\le \area(u_n(\sigma_n( V_\alpha ))) \le N \eps$
for all $n$. Obviously, it is sufficient to show the desired property only
for one such piece of covering, say for $V_1$. To construct the refined
parametrizations $\ti\sigma_{1,n}$ and the covering $V_1 = \cup_i V _{1, i}$,
we consider 4 cases:

\smallskip
{\sl Case 1):  The structures $\sigma_n^*j_n|_{V_1}$ do not change and
$C_{1,n}$ are not isomorphic to the standard node $\cala_0$}.
Then we can realize $(V_1, \sigma_n^*j_n)$ as a constant bounded
domain $D$ in $\cc$. Hence we can consider $u_n \scirc \sigma_n: V_\alpha \to
X$ as pseudoholomorphic maps $u_n :D \to (X,J_n)$. Now we use the ``patching
construction'' of Sacks-Uhlenbeck [S-U].

Fix some $a>0$. Denote $D_{-a} \deff \{ z\in D\;:\; \Delta(z,a) \subset D\}$.
Find a covering of $D_{-a}$ by open sets $U_i \subset D$ with $\diam(U_i) <a$,
such that any $z\in D$ lies in at most 3 pieces $U_i$. Then for any $n$ there
exists at most $3N$ pieces $U_i$ with $\area( u_n(U_i)) > \eps$. Taking a
subsequence, we may assume that the set of such ``bad'' pieces $U_i$ is the
same for all $n$. Repeat successively the same procedure for ${a\over2}$,
${a\over4}$, and so on, and then take the diagonal subsequence. We obtain at
most $3N$ ``bad'' points $y^*_1, \ldots, y^*_l$, such that a subsequence of
$u_n$ converges in $D\bs \{y^*_1, \ldots,y^*_l\}$ strongly, \ie in the 
$L^{1,p}$-topology on compact subsets $K \Subset D\bs \{y^*_1, \ldots,y^*_l\}$.
These ``bad'' points $y^*_1, \ldots, y^*_l$ are characterized by the following
property:
$$
\text{\ \ for any $r>0$\ \ }
\area(u_n(\Delta(y^*_i, r)) >\eps \text{\ \ for $n$ all sufficiently big}.
\eqno(4.1)
$$

\state Remark. As we shall see now, every such point is a place where the
``bubbling'' occurs. Therefore we shall call $y^*_i$ {\sl bubbling
points}. The characterization property of a bubbling point is (4.1).

\smallskip
If there are no bubbling points, \ie $l=0$, then the chosen subsequence
$u_n$ converges strongly and we can finish the proof by induction.

Otherwise we consider the first point $y^*_1 \in D$. Take a disk $\Delta(
y^*_1,\varrho)$ which doesn't contain any other bubbling points $y^*_i$,
$i>1$.

Then for any $n$ we can find the unique $r_n$ such that

(1) $r_n \le {\varrho \over 2}$ and $\area(u_n(\Delta(x, r_n)))\le\eps$ for
any $x\in \barr\Delta (y^*_1,{\varrho \over 2})$;

(2) $r_n$ is maximal \wrt (1).

\smallskip
Then $r_n\lrar 0$, since otherwise for $r^+ \deff {\sf lim\,sup}\; r_n >0$
and for some subsequence $n_k \lrar \infty$ with $r_{n_k} \lrar r^+$
we would have
$$
\area\bigl(u_{n_k}(\Delta(y^*_1, r^+))\bigr) \le \eps,
$$
contradicting to (4.1).

\state Lemma 4.4. {\it For every $n>\!>1$ there exists 
$x_n\in \barr \Delta (y^*_1, {\varrho \over2})$, such that $x_n\to y^*_1$ and 
$\area( u_n( \Delta(x_n, r_n))) =\eps$.}

\state Proof. If not, then for some subsequence $n_k\lrar \infty$ and every
$x\in \barr\Delta (y^*_1, {\varrho \over2})$ we would have $\area(u_{n_k}(
\Delta (x, r_{n_k})))< \eps $. Since $r_n \lrar0$, this would contradict
with the maximality of $r_n$. In particular, there exists a sequence 
$\{ x_n \}$ with the desired properties.

If $x_n$ do not converge to $y^*_1$, then after going to a subsequence we
would find $y'=\lim_{n\to \infty } x_n\not= y^*_1$. By our construction, $y'$
does not coincide with any other bubbling point $y^*_i$. Take $a>0$ such that
$\Delta(y', a)$ contains no bubbling point. Then by {\sl Corollary 3.3} some
subsequence $u_{n_k}$ would converge to some $u' \in L^{1,p}_\loc (\Delta(y', 
a), X)$ in strong $L^{1,p}(K)$-topology for any compact subset $K \Subset
\Delta(y', a)$ and any $p<\infty$. In particular, for sufficiently small
$b<a$ we would have $\area(u_{n_k}(\Delta(y', b))) \lrar \area(u'( \Delta(
y',b))) < \eps$, which would contradict the choice of $r_n$ and $x_n$.
\qed

\smallskip
Using $r_n$ and $x_n$ constructed above, define the maps $v_n:\Delta (0, {
\varrho  \over2r_n}) \to (X,J_n)$ by $v_n(z)\deff u_n(x_n +r_n z)$. By the 
definition of $r_n$ we have
$$
\area(v_n(\Delta (x ,1))\le \eps \ \ \text{for all $x\in \Delta
(0,{\textstyle{\varrho \over 2r_n}}-1)$}.
\eqno(4.2)
$$
On the other hand, $\area(v_n(\Delta(0,1)) = \area(u_n( \Delta (x_n, r_n)) =
\eps$ by {\sl Lemma 4.4}. Thus $v_n$ converge (after going to a subsequence)
on compact subsets in $\cc$ to a nonconstant $J_\infty $-holomorphic map
$v_\infty$ with finite energy. Consequently, $v_\infty$ extends to $S^2$ by
the removable singularity theorem of {\sl Corollary 3.6}.

Since $v_\infty$ is nonconstant, $\norm{ dv_\infty }^2_ {L^2( S^2)}= \area(
v_\infty(S^2))\ge 3\eps$ by {\sl Lemma 3.8} and the choice of $\eps$.
Choose $b>0$ in such a way that
$$
\area(v_\infty (\Delta (0,b)) =
\norm{ dv_\infty}_{L^2(\Delta (0,b))}^2 \ge 2\eps.
\eqno(4.3)
$$
By {\sl Corollary 3.3} this implies that
$$
\norm{du_n}_{L^2(\Delta (x_n,br_n))}^2=
\norm{ dv_n}_{L^2(\Delta (0,b))}^2 \ge\eps.
\eqno(4.4)
$$

For $n>\!>$ we consider the coverings of $V_1$ by 3 sets
$$
\textstyle\mathsurround=0pt
V^{(n)}_{1,1} \deff V_1 \bs \barr\Delta (y^*_1, {\varrho\over2}),
\qquad
V^{(n)}_{1,2}\deff \Delta (y^*_1, \varrho) \bs \barr \Delta (x_n, br_n),
\qquad
V^{(n)}_{1,3} \deff \Delta (x_n, 2br_n).
$$
Fix $n_0$ sufficiently big. Denote $V_{1,1} \deff V^{(n_0)}_{1,1}$, $V_{1,2}
\deff V^{(n_0)}_{1,2}$, and $V_{1,3} \deff V^{(n_0)}_{1,3}$. There exist
diffeomorphisms $\psi_n: V_1 \to V_1$ such that $\psi_n: V_{1,1} \to V^{(n)}
_{1,1}$ is the 
identity, $\psi_n: V_{1,2} \to V^{(n)}_{1,2}$ is a diffeomorphism,
and $\psi_n: V_{1,3} \to V^{(n)}_{1,3}$ is biholomorphic \wrt the complex
structures, induced from $C_{1,n}$.

Thus we have constructed the covering $\{ V_{1,1}, V_{1,2}, V_{1,3} \}$ of
$V_1$ and parametrizations $\sigma'_n \deff \sigma_{1,n} \scirc \psi_n: V_1
\to C_{1, n}$, such that the conditions of {\sl Theorem 4.2} are satisfied.
Moreover, $\area(u_n(\sigma'_n(V_{1,i}))) \le (N-1)\eps$ due to inequality
(4.4). Consequently, we can apply the inductive assumptions to the sequence
of curves $\sigma'_n(V_{1,i})$ and finish the proof by induction.

\medskip\noindent
{\sl Case 2):} $V_1$ is a cylinder, structures $\sigma_n^*j_n|_{V_1}$ vary
with $n$, but conformal radii of $(V_1, \sigma_n^*j_n)$ are bounded uniformly
in $n$. Applying {\sl Lemma 4.3} we can assume that structures $\sigma_n^*
j_n$ converge to a structure of an annulus with finite conformal radius. The
constructions of {\sl Case 1)} go through here with the following minor
modifications. Firstly, we find the set of the bubbling points $ y^*_i\in V_1$,
using the same patching construction and the characterization (4.1). Then
we find diffeomorphisms $\phi_n : V_1 \to V_1$, such that {\sl\.a)} $\phi_n$
converge to the identity map $\id: V_1 \to V_1$; {\sl\.b)} $\phi_n$ are
identical in fixed (\ie independent of $n$) annuli adjacent to the boundary
circles of $V_1$; {\sl\.c)} $\phi_n$ preserve every bubbling point, $\phi_n(
y^*_i) = y^*_i$; and finally {\sl\.d)} for the ``corrected'' parametrizations
$\ti\sigma_n \deff \sigma_n \scirc \phi_n$ the structures $\ti\sigma_n^*j_n|
_{V_1}$ are constant in a neighborhood of every bubbling point $y^*_i$. Then
we repeat remaining constructions of {\sl Case 1)} using the new
parametrizations $\ti\sigma_n$.

\medskip\noindent
{\sl Case 3): Every $C_{1,n} = \sigma_n(V_1)$ is isomorphic to the standard
node $\cala_0$}. Fix identifications $C_{1,n} \cong \cala_0$ such that the
induced  parametrization maps $\sigma_{1, n} : V_1 \to \cala_0$ are the same
for all $n$. Represent $\cala_0$, and hence every $C_{1,n}$, as the union of
two discs $\Delta'$ and $\Delta''$ with identification of the centers $0\in
\Delta'$ and $0\in \Delta''$ into the nodal point of $\cala_0$, still denoted
by $0$. Let $u'_n :\Delta' \to X$ and $u''_n :\Delta'' \to X$ be the
corresponding ``components'' of the maps $u_{1,n} : C_{1,n} \to X$. Find the
common collection of bubbling points $y^*_i$ for both maps
$u'_n :\Delta' \to X$ and $u''_n :\Delta'' \to X$.
If there are no bubbling points, then we obtain the convergence type
{\sl B)} and the proof can be finished by induction.
Otherwise, we  consider the first such point $y^*_1$, which lies, say, on
$\Delta'$. If $y^*_1$ is distinct from the nodal point $0 \in
\Delta'$, then we simply repeat all the constructions of  {\sl Case 1)}.

It remains to consider the case $y^*_1=0 \in \Delta'$. Now one should
modify the arguments of {\sl Case 1} in the following way. Construct the 
sequences of radii $r_n\lrar 0$, of points $x_n \lrar y^*_1=0$,
and of maps $v_n: \Delta(0, {\varrho \over 2r_n}) \to X$, $v_\infty: S^2 \to
X$ as in {\sl Case 1}. Set $R_n \deff |x_n|$, so that $R_n$ is the distance
from $x_n$ to point $0=y^*_1 \in \Delta'$. After rescaling $u_n$ to the maps
$v_n$, the point $0\in \Delta'$ will correspond to the point $z^*_n \deff
-{x_n \over r_n}$ in the definition domain $\Delta(0, {\varrho \over 2r_n})$
of the map $v_n$. We consider 2 subcases.

\smallskip\noindent
{\sl Subcase 3\/$'$): The sequence ${R_n\over r_n}$ is bounded}. This is
equivalent to boundednes of the sequence $ z^*_n$. Going to a
subsequence we may assume that the sequence $z^*_n$ converges to a point
$z^*\in \cc$. This point will be a nodal one for $(S^2, v_\infty)$.
As above, $v_\infty$ is nonconstant and $\norm{ dv _\infty} ^2 _{L^2( S^2)}=
\area( v_\infty( S^2) )\ge 3\eps$. Choose $b>0$ in such a way that
$$
\norm{ dv_\infty}_{L^2(\Delta (0,b))}^2 \ge 2\eps
\eqno(4.5)
$$
and $b \ge 2|z^*|+2$. Due to {\sl Corollary 3.3} for $n>\!>1$ we obtain
the estimate
$$
\norm{du'_n}_{L^2(\Delta' (x_n, br_n))}^2 =
\norm{ dv_n}_{L^2(\Delta (0,b))}^2 \ge \eps.
\eqno(4.6)
$$
Here $\Delta' (x, r)$ denotes the subdisc of $\Delta'$ with center $x$ and
radius $r$. Furthemore, for $n>\!>1$ we have the relation $z^*_n \in 
\Delta( 0, b-1)$, or equivalently, $0 \in  \Delta' (x_n, (b-1) r_n))$.

Define the coverings of $\cala_0$ by 4 sets
$$\mathsurround=0pt
\matrix\format\l\ \ &\l\\
W^{(n)}_1 \deff  \Delta'
\bs \barr\Delta' (0, {\textstyle{\varrho\over2}}),
&
W^{(n)}_2 \deff \Delta' (0, \varrho)
\bs \barr \Delta' (x_n, br_n),
\cr
\noalign{\vskip5pt}
W^{(n)}_3 \deff \Delta' (x_n, 2br_n)
\bs \barr\Delta'(0,{\textstyle {r_n \over 2}}),
&
W^{(n)}_4 \deff \Delta' (0, r_n) \cup
\Delta'',
\endmatrix
$$
and lift them to $V_1$ by putting $V^{(n)}_{1,i} \deff \sigma\inv_{1,n}(
W^{(n)}_i)$. Choose $n_0 >\!> 0$, such that $z^*_{n_0} \in \Delta( 0, b-1)$
and the relation (4.6) holds. Set $V_{1,i} \deff V^{(n_0)} _{1, i}$. Choose
diffeomorphisms $\psi_n: V_1 \to V_1$ such that $\psi_n: V_{1,1} \to V^{(n)}
_{1,1}$ is the identity map, $\psi_n: V_{1,2} \to V^{(n)}_{1,2}$ and $\psi_n:
V_{1,3} \to V^{(n)}_{1,3}$ are diffeomorphisms, and $\psi_n: V_{1,4} \to
V^{(n)}_{1,4}$ is corresponds to isomorphisms of nodes $W^{(n)}_4 \cong
\cala_0$. Set $\sigma'_n \deff \sigma_n \scirc \psi_n$. The choice above can
be done in such a way that the refined covering $\{ V_{1,i} \}$ of $V_1$ and
parametrization maps $\sigma'_n: V_1 \to C_{1,n}$ have the properties of {\sl
Theorem 4.2}. Moreover, relations (4.2) and (4.6) imply the estimate $\area(
u_n (\sigma'_n(V_{1,i})) \le (N-1)\,\eps$. This yields the inductive
conclusion in {\sl Subcase 3\/$'$)}.

\smallskip\noindent
{\sl Subcase 3\/$''$): The sequence $R_n \over r_n$ increase infinitely}.
This means that the sequence$ z^*_n$ is not bounded. Nevertheless $R_n \lrar
0$ since $x_n \lrar 0$. We proceed as follows. Construct of the radius 
$b$ as in  {\sl Case 1)}. For $n>\!>0$ define the coverings of $\cala_0$
by 6 sets
$$\mathsurround=0pt
\matrix\format\l\ \ &\l\\
W^{(n)}_1 \deff  \Delta'
\bs \barr\Delta' (0, {\varrho\over2}),
&
W^{(n)}_2 \deff \Delta' (0, \varrho)
\bs \barr \Delta' (x_n, 2R_n),
\cr\noalign{\vskip5pt}
W^{(n)}_3 \deff  \Delta' (x_n, 4R_n)
\bs \bigl(\barr \Delta' (x_n, {R_n \over 6})
\cup \barr \Delta' (0,{R_n\over 6}) \bigr)
&
W^{(n)}_4 \deff \Delta' (0, {R_n\over 3}) \cup \Delta'',
\cr\noalign{\vskip5pt}
W^{(n)}_5 \deff \Delta' (x_n, {R_n \over 3})
\bs \barr\Delta'(x_n, br_n),
&
W^{(n)}_6 \deff \Delta' (0, 2br_n ),
\endmatrix
$$
and lift them to $V_1$ by putting $V^{(n)}_{1,i} \deff \sigma\inv_{1,n}(
W^{(n)}_i)$. Choose $n_0 >\!> 0$, such that $R_{n_0} >\!> br_{n_0}$ , and set
$V_{1,i} \deff V^{(n_0)} _{1, i}$. Choose diffeomorphisms $\psi_n: V_1 \to
V_1$ such that $\psi_n: V_{1,1} \to V^{(n)} _{1,1}$ is the identity map,
$\psi_n: V_{1,2} \to V^{(n)}_{1,2}$, $\psi_n: V_{1,4} \to V^{(n)}_{1,4}$ and
$\psi_n: V_{1,5} \to V^{(n)}_{1,5}$ are diffeomorphisms, and finally,
$\psi_n: V_{1,6} \to V^{(n)}_{1,6}$ corresponds to isomorphisms of nodes
$W^{(n)}_6 \cong \cala_0$. Set $\sigma'_n \deff \sigma_n \scirc \psi_n$.
Again, this choice can be done in such a way that $\{ V_{1,i} \}$ and
parametrization maps $\sigma'_n: V_1 \to C_{1,n}$ have the properties of
{\sl Theorem 4.2}. As above, we get the estimate $\area(u_n (\sigma'_n(
V_{1,i} )) \le (N-1)\,\eps$ due to (4.2). Thus we get the inductive
conclusion for {\sl Subcase 3\/$''$)} and can proceed further.

\smallskip\noindent
{\sl Case 4):} $V_1$ is a cylinder, structures $\sigma_n^*j_n|_{V_1}$ vary
with $n$, and conformal radii of $(V_1, \sigma_n^*j_n)$ converge to
$+\infty$. Using {\sl Lemma 4.3} we can assume that structures $\sigma_n^*
j_n$ satisfy property \slii of this lemma.

Fix biholomorphisms $\sigma_n(V_1) \cong Z(0, l_n)$. If $\area(u_n(Z(a-1,
a))) \le \eps$ for any $n$ and any $a \in [1, l_n]$, then {\sl Lemma 3.7}
shows that $u_n: \sigma_n(V_1) \to X$ converge to a $J_\infty$-holomorphic 
map from node.

If not, then, after passing to a subsequence, we can find a sequence $\{a_n
\}$ with $a_n \in [1, l_n]$, such that $\area(u_n(Z(a_n \allowbreak -1, a_n)))
\ge \eps$. If $a_n$ is bounded, say $a_n \le a^+$, then we cover $Z(0,
l_n)$ by the sets $V_{1,1} \deff Z(0, a^+ +2)$ and  $V_{1,2} \deff Z(a^+ +1,
l_n)$. If $l_n - a_n$ is bounded, say $l_n - a_n \le a^+$, then we cover $Z(0,
l_n)$ by the sets $V_{1,1} \deff Z(0, l_n - a^+ +2)$ and  $V_{1,2} \deff
Z(l_n -a^+ +1, l_n)$. In the remaning case, when both $a_n$ and $l_n - a_n$
increase infinitely, we cover $Z(0, l_n)$ by 3 sets $V_{1,1}^{(n)} \deff
Z(0,a_n-1)$, $V_{1,2}^{(n)} \deff Z(a_n-2, a_n+1)$, and $V_{1,3}^{(n)} \deff
Z(a_n, l_n)$. 

Identify $V_1$ with the cylinder $Z(0, 5)$ and cover it by 2 or respectively 
3 successive cylinders $V_{1,i}$, \eg by 2 cylinders $V_{1,1} \deff Z(0, 2)$ 
and $V_{1,2} \deff Z(1, 5)$, or respectively by 3 cylinders $V_{1,1} \deff 
Z(0, 2)$, $V_{1,2} \deff Z(1, 4)$, and $V_{1,3} \deff Z(3, 5)$.
Find diffeomorphisms $\psi_n : V_1 \to V_1$, identical
in the neighborhood of the boundary of $V_1$ and such that $\psi_n(V_{1,i}) =
\sigma_n ^{-1} V^{(n)} _{1,i}$. Define the new parametrizations $\sigma'_n
\deff \sigma_n \scirc \psi_n$. Note that we may additionally assume that if
the conformal radius of $\sigma'_n(V_{1,i})$ is independent of $n$ then the
structure $\sigma'_n{}^*j_n \ogran_{V_{1,i}}$ is also is independent of $n$.

By the construction, we get the following property of the covering $\{ V_{1,i}
\}$ and new parametrizations $\sigma'_n$: For any $V_i$ either
we have the estimate
$$
\area(u_n(\sigma'_n(V_{1,i}))) \le (N-1)\,\eps.
$$
or the structures $\sigma'_n{}^*j_n|_{V_{1,i}}$ do not depend on $n$. Thus we
reduce our case to the situation which is covered either by the inductive
assumption or by {\sl Case 1)}.

\medskip
The proof of the theorem can be finished by induction. The fact that the
limit curve $(C_\infty, u_\infty)$ remains stable over $X$ is proved in {\sl
Lemma 4.5} below. \qed

\medskip
\state Remark. We explain here the meaning of the constructions used in the
proof. We start with {\sl Case 1)} where $J_n$-holomorphic maps from a fixed
domain $D\subset \cc$ are treated. The bubbling points $y^*_i$ appear in this 
case as those ones where the strong convergence of maps $u_n: D \to (X, J_n)$
fails. The patching construction of Sacks and Uhlenbeck insures us that the
``convergence failure'' set is finite and gives an effective estimate of
the number of bubbling points by the upper bound of the area, $l \le 3N$ in
our situation. The characterization property (4.2) of bubbling points is
essentially due to Sacks and Uhlenbeck; the only difference is that we use
the area of the map $u$, \ie $L^2$-norm of $du$, see {\sl Definition 1.4} , 
whereas in [S-U] the $\norm{du}_{L^\infty}$ is used. The next step, namely
the construction of maps $v_n$ as rescaling of the $u_n$ and the
existence of the limit $v_\infty$, is also due to Sacks and Uhlenbeck.

The explicit construction of the map $v_\infty$ suggests the interpretation 
of the curve $(S^2, v_\infty)$ as a ``bubbled sphere'' and $y^*_1$ as
the point where the ``bubbling'' occurs. Moreover, one obtains natural
partitions (one for each $n>\!>0$) of $D$ into three pieces: $D$ minus fixed
small neighborhood of $y^*_1$; disks $(\Delta(x_n, br_n), u_n)$ representing
pieces $(\Delta(0, b), v_\infty)$ and approximating a sufficiently big part
$(\Delta(0, b), v_\infty)$ of the bubbled sphere; and the ``part inbetween''.

These latter ``parts inbetween'' appear to be the annuli of infinitely 
growing conformal radii, considered in {\sl Case 4)}. Since neither
outer nor inner boundary circle should be preferred in some way, we consider
them as long cylinders $C_n = Z(0,l_n)$ with $l_n \lrar \infty$ according to 
{\sl Definition 3.2}. {\sl Lemma 3.7} provides a ``good'' convergence model
for long cylinders, stated above as convergence type {\sl A)}. If for a 
sequence $(C_n, u_n)$ such convergence fails, then there must exist
subannuli $A_n \subset C_n$ of a constant conformal radius, for which
$\area(u_n(A_n)) \ge \eps$.

In both cases --- a constant domain $D$ or long cylinders --- we proceed by
cutting the curves into smaller pieces. The situation we come to is simpler
in the following sense. The obtained curves either converge or have the
upper bound for the area smaller by the fixed constant $\eps$. Thus the 
induction leads finally to decomposition of the curves into pieces for
which one of the convergence types {\sl A)--C)} holds. The possibility to
glue these final pieces together is insured by the fact that the partitions
above are represented by appropriate coverings satisfying the conditions of
{\sl Theorem 4.2}.

Considering curves with nodes, an additional attention should be payed to the
case when bubbling appears at a nodal point. This situation is considered in
{\sl Case 3)}. The constructed points $x_n$ and radii $r_n$ describe the
``center'' and the ``size'' of energy localization of the bubbling,
represented by the sequence of the maps $v_n$ tending to $v_\infty$. So the
convergence picture depends on whether the energy localization occurs
near the nodal point ({\sl Subcase 3\/$'$)}) or away of it ({\sl Subcase
3\/$''$)}). As a result, the nodal point can either remain on the ``bubbled''
sphere $(S^2, v_\infty)$ or move into the ``part inbetween'', which is
represented by long cylinders.

In {\sl Subcase 3\/$'$)} we remove neighbourhoods of the nodal point from the
disks $(\Delta(0, b), v_n)$ and thus get 4 pieces of covering instead of 3
in {\sl Case 1)}. In {\sl Subcase 3\/$''$)} situation is more complicated,
since we must take into account the position of the nodal point in the long
cylinders --- the ``parts inbetween''. Thus we must consider now the sequence
of cylinders with one marked point, \ie the sequence of pants. The fact that
$R_n \lrar 0$ and ${r_n \over R_n} \lrar 0$ means that the conformal structure 
of those pants is not constant and converges to  one of the sphere with 3 
punctures.

In order to have the covering pieces with the convergence types {\sl A)--C)},
we choose an appropriate refinement of the covering. After it, we obtain 2
sequences of long cylinders, describing appearance of 2 new nodal points.
The first one corresponds to the part between the ``original'' nodal point
and bubbled sphere and is represented by $V_5$, whereas the other one,
represented by $V_2$, lies on the other side from the ``original'' nodal
point. Besides, we fix a neighborhood of the ``original'' nodal point which
has a constant complex structure and is topologically an annulus with the disc
$\Delta''$ attached to the nodal point. To satisfy the requirements of {\sl
Theorem 4.2}, we cover the neighborhood by 2 pieces, the pants $V_3$ and the
piece $V_4$ parametrizing the nodes $W^{(n)}_4$. This explains appearance of
6 pieces of covering in {\sl Subcase 3\/$''$)}.

\smallskip
\state Lemma 4.5. \it The limit curve $(C_\infty, u_\infty)$ constructed
in the proof of {\sl Theorem 1.1} is stable over $X$. \rm

\state Proof. If $(C_\infty, u_\infty)$ is unstable over $X$, then either
$C_\infty$ is a torus with $u_\infty$ constant, or $C_\infty$ should
have a component $C' \subset C_\infty$, such that $\ti C'$ is a sphere with
at most 2 marked points and $u_\infty(C')$ is a point.

The case of a constant map from a torus is easy to handle. In fact, in this
case all $C_n$ must be also tori with $\area(u_n(C_n))$ sufficiently small
for $n>\!>1$. Cover every $C_n$ by infinite cylinder $Z(-\infty, +\infty)$
and consider compositions $\ti u_n: Z(-\infty, +\infty) \to X$ of $u_n$ with
the covering maps. Since $\area(u_n(C_n))\approx 0$, {\sl Corollary 3.6} can
be applied to show that every $\ti u_n$ extends to a $J_n$-holomorphic
map from $S^2$ to $X$. Consequently, $\area(\ti u_n(Z(-\infty, +\infty)))$
must be finite. On the other hand, $\area(u_n(C_n))>0$ due to the stability
condition, and hence $\area(\ti u_n(Z(-\infty, +\infty)))$ must be infinite.
This contradiction rules out the case of a torus.

The same argumentations go through in the case, when $C_\infty$ is the sphere
with no marked points. Then the curves $C_n$ are also parametrized by
the sphere $S^2$. The condition of instability means that
$\area(u_\infty(C_\infty))=0$. Due to {\sl Corollary 3.3}, $\area(u_n(C_n))$
must be sufficiently small for $n>\!>1$. Now {\sl Lemma 3.8} and  the
stability of $(C_n,u_n)$ show that this is im possible.

Now consider the cases when the limit curve $C_\infty$ has a ``bubbled''
component $C'$, which  is the sphere with 1 or 2 marked points. If $C'$ has 1
marked point, then $C'$ must appear as a ``bubbled'' sphere $(S^2, v_\infty)$
in the constructions of {\sl Cases 1)--3)} in the proof of {\sl Theorem 1.1}.
But these constructions yield only non-trivial ``bubbled'' spheres, for which
$v_\infty \not= const$. Thus such component $C'$ must be stable.

\smallskip
In the remaining case, \ie if there exists a component $C'$ with 2 marked 
points, we consider a
domain $U\subset C_\infty$, which is the union of the component $C'$ and
neighborhoods of the marked point on $C'$. If $C'$ is an unstable component,
then $\area(u_\infty(C'))=0$ and we can achieve the estimate $\area( u_\infty(
U))< \eps$ taking $U$ sufficiently small. Set $\Omega \deff \sigma_\infty
\inv(U)$, where $\sigma_\infty: \Sigma \to C_\infty$ is the parametrization of
$C_\infty$. Let $\gamma_1$ and $\gamma_2$ be the pre-images of marked points
on $C'$. Then $\Omega$ must be a topological annulus, and $\gamma_i$, $i=1,2$,
disjoint circles generating the group $\pi_1 (\Omega) = \zz$. Further, $C'$
must be a ``bubbling'' component of $C_\infty$, \ie at least for one of the
circles $\gamma_i$, $i=1,2$, the images $\sigma_n(\gamma_i)$ are not nodal
points of $C_n$ but smooth circles.

If the both circles $\gamma_1$ and $\gamma_2$ are of this type, then $U_n 
\deff \sigma_n (\Omega)$ satisfies the conditions of {\sl Lemma 3.7}. In 
this case we should have the convergence type {\sl A)}, and hence the limit 
piece $\sigma_\infty(\Omega)$ should be isomorphic to the node $\cala_0$.

In the case when only one circle, say $\gamma_1$, corresponds to nodal points
on $C_n$, and for the other one the images $\sigma_n(\gamma_2)$ are smooth
circles, then the domains $\sigma_n(\Omega)$ must be isomorphic to the node
$\cala_0$. Furthermore, due to the condition $\area(u_\infty(\sigma_\infty
(\Omega))) < \eps$ we have $\area(u_n(\sigma_n (\Omega))) < \eps$.
Consequently, we must have the convergence type {\sl B)} and the
unstable component $C'$ could not appear. \qed


\bigskip\bigskip
\newsection{}{Curves with boundary on totally real submanifolds}

\smallskip\noindent
In this section we consider the behavior of pseudoholomorphic curves over an
almost complex manifold $(X,J)$ with boundary on totally real submanifold(s).
As in the "interior" case, we need to allow some type of boundary
singularity.

\state Definition 5.1. {\sl The set $\cala^+ \deff \{ (z_1,z_2) \in \Delta^2
\;:\; z_1 \cdot z_2=0, \im z_1\ge 0, \im z_2\ge 0\}$ is called the {\it
standard boundary node}. A curve $\barr C$ with boundary $\d C$ is called a
{\it nodal curve with boundary} if:

\sli $C$ is a nodal curve, possibly disconnected;

\slii $\barr C=C\cup \d C$ is connected and compact;

\sliii every boundary point $a\in \d C$ has a neighborhood homeomorphic
either to the half-disk $\Delta^+ \deff \{ z \in \Delta
\;:\; \im z\ge0 \}$, or to the standard boundary node $\cala^+$.

In the last case $a\in \d C$ is called a {\it boundary nodal point}, whereas
nodal points of $C$ are called {\it interior nodal points}. Both boundary and
interior nodal points are simply called nodal points.
}

\state Definition 5.2. {\sl Let $(X,J)$ be an almost complex manifold. A pair
$(\barr C, u)$ is called a {\it curve with boundary over $(X,J)$} if $\barr
C=C\cup \d C$ is a nodal curve with boundary, and $u: \barr C \to (X,J)$ is a
continuous $L^{1,2}$-smooth map, which is pseudoholomorphic on $C$.
}

\smallskip
A curve $(C,u)$ with boundary is stable if the same condition as in {\sl
Definition 1.5} on the automorphism groups of compact irreducible components
is satisfied.

\smallskip
\state Remark. One can see, that $\barr C$ has a uniquely defined real
analytic structure, such that the normalization $\barr C{}^{\sf nr}$ is a
real analytic manifold with boundary. More
precisely, the pre-image of every boundary nodal point $a_i$ consists of two
points $a'_i$ and $a''_i$. The normalization
map $s: \barr C{}^{\sf nr} \to \barr C$ glues each pair $a'_i$, $a''_i)$ of 
points on $\barr C{}^{\sf nr}$ together into nodal points $a_i= s(a'_i)= 
s(a''_i)$ on $\barr C$. This implies that the notion of an $L^{1,p}$-smooth 
map, $p>2$, and that of also a continuous $L^{1,p}$-smooth map $u: \barr C 
\to X$ re well defined.

\state Definition 5.3. {\sl We say that a real oriented surface with boundary
$(\Sigma, \d\Sigma)$ {\it parameterizes} a nodal curve with boundary $C$ if
there is a continuous map $\sigma :\barr\Sigma \to \barr C$ such that:

\sli if $a\in C$ is an interior nodal point, then $\gamma_a \deff \sigma\inv
(a)$ is a smooth imbedded circle in $\Sigma$;

\slii if $a\in\d C$ is a boundary nodal point, then $\gamma_a\deff \sigma
\inv (a)$ is a smooth imbedded arc in $\Sigma$ with end points on $\d\Sigma$,
transversal to $\d\Sigma$ at these points;

\sliii if $a,b \in \barr C$ are distinct (interior or boundary) nodal points,
then $\gamma_a\cap \gamma_b= \emptyset$;

\sliv $\sigma :\barr\Sigma \bs \bigcup_{i=1}^N\gamma_{a_i}\to \barr C \bs \{
a_1,\ldots ,a_N\} $ is a diffeomorphism, where $a_1,\ldots ,a_N$ are all
(interior and boundary) nodal points of $\barr C$.
}

\smallskip\rm
Recall that a real subspace $W$ of a complex vector space is called {\it
totally real} if $W \cap \isl W =0$. Similarly, a $C^1$-immersion $f:W \to X$
is called {\it totally real} if for any $w\in W$ the image $d f(T_wW)$ is
a totally real subspace of $T_{f(w)}X$.

Let $(\barr C, u)$ be a stable curve with boundary over an almost complex
manifold $(X,J)$ of a complex dimension $n$.

\state Definition 5.4. {\sl We say that $(\barr C, u)$ satisfies {\it totally
real boundary condition $\mib W$ of type $\bfbeta$} if

\sli $\bfbeta= \{\beta_i\}$ is a collection of arcs with disjoint interiors,
which defines a decomposition of the boundary $\d C= \cup_i \beta_i$;
moreover, we assume that every boundary nodal point is an endpoint for 4 arcs
$\beta_i$;

\slii ${\mib W} = \{(W_i, f_i)\}$ is a collection of totally real immersions
$f_i: W_i \to X$, one for every $\beta_i$;

\sliii there are given continuous maps $u^{(b)}_i :\beta_i \to W_i$ {\it
realizing conditions $\mib W$}, \ie $f_i \scirc u^{(b)}_i = u|_{\beta_i}$.
}

\state Remarks.~1. We shall consider (immersed) totally real submanifolds
only of {\sl maximal real dimension} $n= \dimc X$.

\state 2. If $\bfbeta$ is a collection of arcs as above, a parametrization
$\sigma: \barr \Sigma \to \barr C$ induces a collection of arcs $\sigma\inv(
\bfbeta) \deff \{ \sigma\inv(\beta_i)\;:\; \beta_i \in \bfbeta$ with the
properties similar to \sli of {\sl Definition 5.4}. Thus, $\sigma\inv(
\beta_i)$ have disjoint interiors, $\cup_i\sigma\inv(\beta_i) = \d\Sigma$,
and for any boundary node $a\in \barr C$ every endpoint of the arc $\beta_a=
\sigma\inv(a)$ is an endpoint of two arcs $\sigma\inv(\beta_i)$. Since
$\bfbeta$ is completely determined by $\sigma\inv(\bfbeta)$, we shall denote 
the both collections simply by $\bfbeta$ and shall not distinguish them
considering boundary conditions.

\medskip
A totally real boundary condition is a suitable elliptic boundary condition for
an elliptic differential operator $\dbar$ of Cauchy-Riemann type. In
particular, all statements about ``inner'' regularity and convergence for
pseudoholomorphic curves remain valid near ``totally real'' boundary. As in
``inner'' case, to get some ``uniform'' estimate at boundary one needs $W$ to
be ``uniformly totally real''.

\state Definition 5.5. {\sl Let $X$ be a manifold with a Riemannian metric
$h$, $J$ a continuous almost complex structure, $W$ a manifold, and $A_W
\subset W$ a subset. We say that an immersion $f: W \to X$ is {\it
$h$-uniformly totally real along $A_W$ with a lower angle $\alpha = \alpha(W,
A_W, f) > 0$}, \iff

\sli $df: TW \to TX$ is $h$-uniformly continuous along $A_W$;

\slii for any  $w \in A_W$ and any $\xi \not =0 \in T_wW$ the angle
$\angle_h\bigl(Jdf(\xi), df(T_wW) \bigr) \ge \alpha$.
}

\medskip
We start with a generalization of the First Apriori Estimate. Define the
half-disks $\Delta^+(r) \deff \{ z\in \Delta(r) \;:\; \im z \ge 0 \}$ with
$\Delta^+ = \Delta^+(1)$ and $\Delta^- \deff \{ z \in \Delta \; :\; \im z \le
0 \}$. Set $\beta_0 \deff (-1, 1) \subset \d\Delta^+$. Let $X$ be a manifold
with a Riemannian metric $h$, $A\subset X$ a subset, $J^*$ a continuous
almost complex structure, $f: W \to X$ a totally real immersion, and $A_W
\subset W$ a subset.

\state Lemma 5.1. {\it Suppose that $J^*$ is $h$-uniformly continuous on $A$
with the uniform continuity modulus $\mu_{J^*}$, and that $f : W \to X$ is
$h$-uniformly totally real along $A_W$ with a lower angle $\alpha_f>0$ and
the uniform continuity modulus $\mu_f$. Then for every $2< p<\infty $ there
exists an $\eps^b_1 =\eps^b_1(\mu_{J^*}, \alpha_f, \mu_f )$ (independent of
$p$) and $C_p= C(p, \mu_{J^*}, \alpha_f, \mu_f )$, such that the following
holds:

If $J$ is a continuous almost complex structure on $X$ with $\norm{ J- J^*}
_{L^\infty(A)} <\eps^b_1$, if $u\in C^0 \cap L^{1,2}(\barr\Delta^+ ,X)$ is
$J$-holomorphic map with $u(\Delta) \subset A$ and with the boundary
condition $u|_{\beta_0} =f \scirc u^b$ for some continuous $u^b: \beta_0
\to A_W \subset W$, then the condition $\norm{du}_{ L^2( \Delta^+ )}
<\eps^b_1$ implies the estimate
$$
\norm{du}_{L^p(\Delta^+({1\over 2}) )}\le C_p\cdot
\norm{du}_{L^2(\Delta^+ )}.\eqno(5.1)
$$
}

\state Proof. Suppose additionally that $\diam(u(\Delta^+))$ is sufficiently
small. Then we may assume that $u(\Delta^+)$ is contained in some chart $U
\subset \cc^n$, such that $\norm{ J- J\st} _{L^\infty(U)}$ is also small
enough. Let $z=(z_1, \ldots, z_n)$ be $J\st$-holomorphic coordinates in $U$,
such that $u(0)= \{z_i= 0\}$. Making an appropriate diffeomorphism of $U$, we
may additionally assume that $W_0 \deff f(W) \cap U$ lies in $\rr^n$ and that
$J = J\st$ along $W_0$.

Consider the trivial bundle $E \deff \Delta \times \cc^n$ over $\Delta$ with
complex structures $J\st$ and $J_u \deff J \scirc u$. We can consider $u$ as
a section of $E$ over $\Delta^+$ satisfying equation $\dbar_{J_u} u\deff \d_x
u+ J_u \d_y u=0$. Over $\beta_0$ we get a $J_u$-totally real subbundle $F
\deff \beta_0 \times \rr^n$, such that $u(\beta_0) \subset F$. Let $\tau$
denote a complex conjugation in $\Delta$ and also a complex conjugation in $E$
with respect to $J\st$. Extend $J_u$ on $E|_{\Delta^-}$ as the composition
$- \tau\scirc J_u \scirc \tau$. This means that for $z\in \Delta^-$ and we get
$$
J_u(z): v \mapsto  \tau v \in E_{\tau z}
\mapsto J_u(\tau z)(\tau v) \in E_{\tau z}
\mapsto  -\tau J_u(\tau z)(\tau v) \in E_z .
\eqno(5.2)
$$

Since $J_u=J \scirc u$ coincides with $J\st$ along $\beta_0$, this extension
is also continuous. Further, for $v\in L^1( \Delta^+, E)$ define the
extension $\ex(v)$ by setting $v(z)\deff \tau v(\tau z)$. This gives
continuous linear operators $\ex: L^p( \Delta^+, E) \to L^p( \Delta^+, E)$
for any $p\in [1, \infty]$. An important property of operator $\ex$ is that
if $v \in L^{1,p} (\Delta^+, E)$ with $1\le p\le \infty$ (resp.\ $v\in C^0(
\barr \Delta^+ )$) satisfies the boundary condition $v|_{\beta_0} \subset
F$, then $\ex v \in L^{1,p}(\Delta, \cc^n)$ (resp.\ $\ex v \in C^0 (\barr
\Delta)$). Let us denote by $L^{1,p} (\Delta^+, E,F)$ (resp.\ by $C^0(
\Delta^+, E,F) $) the spaces of all such $v$ with the boundary condition 
$v|_{\beta_0} \subset F$.

Since $\d_x(\tau v) =\tau (\d_x v)$ and $\d_y(\tau v) = - \tau (\d_y v)$ for 
$v\in L^{1,1}(\Delta^+, E)$, we get $\dbar_{J_u} (\ex v) =\ex (\dbar_{
J_u} v)$ for any $v \in L^{1,p} (\Delta^+, E,F)$. In particular, for $\ti u
\deff \ex u \in C^0 \cap L^{1,2} (\Delta, E)$ we have $\dbar_{J_u} \ti u=0$.

Starting from this point we can repeat the steps of the proof of 
{\sl Lemma 3.1}. \qed

\state Remark. We shall refer to the construction of a complex structure
$J_u$ in $E$ over $\Delta^-$ and (resp.\ of a section $\ti u$ of $E$ over
$\Delta^-$) as {\sl extension of $J \scirc u$ (resp.\ of $u$) by the
reflection principle}.

\smallskip\rm
Let $X$ be a manifold with a Riemannian metric $h$, $J^*$ a continuous almost
complex structures on $X$, $A\subset X$ a closed $h$-complete subset, such
that $J^*$ is $h$-uniformly continuous on $A$, and $f_0 :W \to X$ an
immersion, which is $h$-uniformly totally real on some close $f^*h$-complete
subset $A_W\subset W$.

\state Corollary 5.2. {\it Let $\{ J_n \}$ be a sequence of continuous almost
complex structures on $X$ such that $J_n$ converge $h$-uniformly on $A$ to
$J$, and $f_n: W \to X$ a sequence of totally real immersion such that $df_n$
converge $h$-uniformly on $A_W$ to $df_0$. 

Furthermore, let $u_n\in C^0\cap L^{1,2}( \Delta^+ ,X)$ be a sequence of 
$J_n$-holomorphic maps, such that $u_n(\Delta^+) \subset A$, $\norm{du_n}
_{L^2 (\Delta^+ )} \le \eps^b_1$, $u_n(0)$ is bounded in $X$, and $u_n| 
_{\beta_0} = f_n \scirc u^b_n$ for some continuous $u^b_n :\beta_0 \to A_W$.

Then there exists a subsequence $u_{n_k}$ which $L^{1,p}_\loc (\Delta^+
)$-converges to a $J$-holo\-morphic map $u_\infty$ for all $p< \infty$.
}

Here $\beta_0 = (-1,1) \subset \d \Delta^+$ and $L^{1,p}_\loc (\Delta^+
)$-convergence means $L^{1,p} (\Delta^+(r) )$-con\-ver\-gence for all $r<1$,
\ie convergence up to boundary component $\beta_0$.

\state Proof. The statement is a generalization of {\sl Corollary 3.3}. 
The proof of that statement goes through with an appropriate modification 
using the reflection principle.

\medskip
Consider now a generalization of the Second Apriori Estimate. Instead of 
``long cylinders'' we now have ``long strips" satisfying appropriate boundary 
conditions.

\state Definition 5.5. {\sl Define a {\it strip $\Theta(a,b) \deff (a,b)
\times [0,1]$} with the complex coordinate $\zeta \deff t - i\theta$, $t 
\in (a,b)$, $ \theta \in [0,1]$. Define also $\Theta_n \deff \Theta(n-1, n)$,
$\d_0 \Theta(a,b) \deff (a,b) \times \{0\}$, and $\d_1 \Theta(a,b) \deff (a,
b) \times \{1\}$.
}

We are interested in maps $u: \Theta(a,b) \to X$, which are holomorphic with
respect to the complex coordinate $\zeta$ on $\Theta(a,b)$ and a continuous
almost complex structure $J$ on $X$, and which satisfy boundary conditions
$$
u\ogran_{\d_0\Theta(a,b)} = f_0 \scirc u^b_0,
\hskip6em
u\ogran_{\d_1\Theta(a,b)} = f_1 \scirc u^b_1,
$$
with some $J$-totally real immersions $f_{0,1} :W_{0,1} \to X$ and continuous
maps $u_{0,1} :\d_{0,1}\Theta(a,b) \to W_{0,1}$. First we consider the 
linear case.

\state Lemma 5.3. {\it Let $W_0$ and $W_1$ be $n$-dimensional totally real
subspaces in $\cc^n=(\rr^n,J\st )$. Then there exist a constant $\gamma_W=
\gamma(n, W_0, W_1)$ with $0< \gamma_W <1$ such that for any holomorphic map
$u: \Theta(0,3) \to \cc^n$ with the boundary conditions
$$
u(\d_0 \Theta(0,3)) \subset W_0 \qquad
u(\d_1 \Theta(0,3)) \subset W_1
\eqno(5.2)
$$
we have the following estimate:
$$
\int_{\Theta_2} |du|^2 dt\,d\theta \le \msmall{\gamma_W\over2}
\left( \int_{\Theta_1} |du|^2 dt\,d\theta +
\int_{\Theta_3} |du|^2 dt\,d\theta \right).
\eqno(5.3)
$$
}

\state Proof. Let $L^{1,2}_W([0,1], \cc^n)$ be a Banach manifold $v(\theta) 
\in L^{1,2}([0,1], \cc^n)$, such that $v(0) \in W_0$ and $v(1) \in W_1$.
Consider a nonnegative quadratic form $Q(v) \deff \int_0^1 |\d_\theta v(
\theta)|^2 d\theta$. Since $Q(v) + \norm{v} ^2_{L^2} =\norm{v}^2 _{L^{1,2}}$
and the imbedding $L^{1,2}_W([0,1], \cc^n) \hook L^2([0,1],\cc^n)$ is compact,
we can decompose $L^{1,2}_W([0,1], \cc^n)$ into a direct Hilbert sum of
eigenspaces $\ee_\lambda$ of $Q$ \wrt $\norm{v}^2_{L^2}$. This means that
$v_\lambda$ belongs to $\ee_\lambda$ iff 
$$
\int_0^1 \<
\d_\theta v_\lambda(\theta), \d_\theta w(\theta) \> d\theta = \int_0^1
\lambda\< v_\lambda(\theta), w(\theta) \> d\theta ,\eqno(5.4)
$$
for any $w\in L^{1,2}_W([0,1], \cc^n)$. Here $\<\cdot,\cdot \>$ denotes 
a standard {\sl $\rr$-valued} scalar product
in $\cc^n$. Integrating by parts yields
$$
\int_0^1 \< \d^2_{\theta\theta} v_\lambda(\theta)
 + \lambda v_\lambda(\theta), w(\theta) \> d\theta + \<\d_\theta v_\lambda(
\theta), w(\theta) \>|_{\theta=1} - \<\d_\theta v_\lambda(
\theta), w(\theta) \>|_{\theta=0} =0.
$$
This implies that $v_\lambda$ belongs to $\ee_\lambda$ \iff
$\d^2_{\theta\theta} v_\lambda(\theta) + \lambda v_\lambda( \theta)=0$,
$\d_\theta v_\lambda(1) \perp W_1$, and $\d_\theta v_\lambda(0) \perp W_0$.
Since $J\st$ is $\<\cdot, \cdot\>$-orthogonal, we can conclude that
$J\d_\theta v_\lambda(\theta) \in \ee_\lambda$.

Positivity and compactness of $Q$ \wrt $\norm\cdot_{L^2}$ imply that all 
$\ee_\lambda$ are finite dimensional and  $\ee_\lambda = \{0\}$ for 
$\lambda<0$. Further, since $\d_\theta v=0$ for any $v\in \ee_0$, the space
$\ee_0$ consists of constant functions with values in $W_0 \cap W_1$.

Now let $u: \Theta(0,3) \to \cc^n$ be a holomorphic map with the boundary
condition (5.2). We can represent $u$ in the form $u(t, \theta)= \sum_\lambda
u_\lambda(t,\theta)$ with $u_\lambda(t, \cdot) \deff \pr_\lambda (u(t,\cdot))
\in L^{1,2}([0,3],\ee_\lambda)$. Since $J\d_\theta$ is an endomorphism of
every $\ee_\lambda$, every $u_\lambda(t,\theta)$ is also holomorphic. In
particular, $u_0$ is also holomorphic and constant in $\theta$. Thus $u_0$ is
constant.

Since $u$ is harmonic, $(\d^2_{tt} + \d^2_{\theta \theta})u =0$, we get
$\d^2_{tt} u_\lambda(t,\theta)\allowbreak = \lambda u_\lambda(t,\theta)$. For
$\lambda > 0$ this yields $u_\lambda(t,\theta) = e^{+\sqrt \lambda t}
v^+_\lambda(\theta) + e^{-\sqrt\lambda t} v^-_\lambda(\theta)$ with $v^\pm
_\lambda(\theta) \in \ee_\lambda$. Fixing an orthogonal $\rr$-basis of
$\ee_\lambda $ $v^i_\lambda$ we write every $u_\lambda$ in the form
$$
u_\lambda(t,\theta)= \sum_i (a^i_\lambda e^{+\sqrt \lambda t} +
b^i_\lambda e^{-\sqrt \lambda t}) v^i_\lambda(\theta)
$$
with real constants $a^i_\lambda$, $b^i_\lambda$. Since $u_0(t,\theta)$ is
constant, $\norm{du}^2_{L^2(\Theta_k)} = 2\norm{\d_\theta u}^2_{L^2(\Theta_k
)} = \sum_{\lambda,i} 2\lambda \int_{k-1}^k (a^i_\lambda e^{+\sqrt\lambda t}
+b^i_\lambda e^{-\sqrt\lambda t})^2 d\theta $. Here we use $(5.4)$ and
$L^2$-orthonormality of $v^i_\lambda $. This leads us to the problem of 
finding the smallest possible constant $\gamma$ in the inequality
$$
\int_1^2 (ae^{\alpha t} + be^{-\alpha t})^2 dt
\le \msmall{\gamma\over2} \left(
\int_0^1 (ae^{\alpha t} + be^{-\alpha t})^2 dt +
\int_2^3 (ae^{\alpha t} + be^{-\alpha t})^2 dt \right)
\eqno(5.5)
$$
with $a,\,b\in\rr$ for given $\alpha>0$. The integration gives
$$
a^2 e^{3\alpha}\msmall{ e^\alpha - e^{-\alpha} \over 2\alpha}
+ b^2 e^{-3\alpha}\msmall{ e^\alpha - e^{-\alpha} \over 2\alpha}
+2ab \le
$$
$$
\le \msmall{\gamma\over2} \left(
a^2 e^{3\alpha}\msmall{ (e^\alpha - e^{-\alpha}) (e^{2\alpha} + e^{-2\alpha})
 \over 2\alpha}+
b^2 e^{-3\alpha}\msmall{ (e^\alpha
- e^{-\alpha}) (e^{2\alpha} + e^{-2\alpha})
\over 2\alpha}
+ 4ab \right)
$$
or, equivalently,
$$
a^2 e^{3\alpha}\msmall{ (e^\alpha - e^{-\alpha})
(e^{2\alpha} + e^{-2\alpha} -2/\gamma ) \over 2\alpha}
+
b^2 e^{-3\alpha}\msmall{ (e^\alpha - e^{-\alpha})
(e^{2\alpha} + e^{-2\alpha} -2/\gamma) \over 2\alpha}
$$
$$
+ 4ab(1-1/\gamma)
\ge 0
$$
The determinant of the last quadratic form in $a,b$ is
$$
\left(\!\!\msmall{ (e^\alpha - e^{-\alpha})
(e^{2\alpha} + e^{-2\alpha} -2/\gamma )
 \over 2\alpha}\!\!\right)^{\!\!2} - 4(1-1/\gamma)^2
=
4\left(\!\!\msmall{ \sh \alpha
(\ch\,2\alpha -1/\gamma )
 \over \alpha} \!\!\right)^{\!\!2} - 4(1-1/\gamma)^2
$$
$$
\ge 4(\ch\,2\alpha -1/\gamma )^2 -4(1-1/\gamma)^2 =
4(\ch\,2\alpha -1)(\ch\,2\alpha +1 -2/\gamma).
$$
Thus inequality (5.5) holds for every $a,b \in\rr$ provided $\gamma \ge{2 
\over 1 + \ch\,2\alpha}<1$.

Note that the minimal positive eigenvalue $\lambda_1>0$ of $Q$ exists.
Thus we can set $\gamma_W \deff {2 \over 1 + \ch(2 \sqrt{\lambda_1} )}<1$ 
in estimate (5.3). \qed

\smallskip
\state Remark. A behavior of $\gamma_W$ as a function of $\lambda_1=
\lambda_1(W_0, W_1)$ shows that $\gamma_W<1$ can be chosen the same for all
pairs $(\wt W_0,\wt W_1)$ sufficiently close to $(W_0, W_1)$ provided $\dim(
\wt W_0 \cap \wt W_1) = \dim (W_0 \cap W_1)$. Vice versa, if we perform
sufficiently small deformation of $(W_0, W_1)$ into $(\wt W_0, \wt W_1)$ with
$\dim (\wt W_0 \cap \wt W_1) < \dim (W_0 \cap W_1)$, then some $v \in\ee_0(
W_0, W_1)$ will become an eigenvector $\ti v \not\in \ee_0 (\wt W_0, \wt W_1) 
= \wt W_0 \cap \wt W_1$, but with sufficiently small eigenvalue
$\lambda_1(\wt W_0,\wt W_1) >0$, so that the best possible $\gamma_{\wt W}$
will be arbitrary close to $1$. Thus the uniform separation of $\gamma_W$
from 1 under small perturbation of $(W_0, W_1)$ is equivalent to uniform
separation of $\lambda_1(W_0, W_1)$ from 0, which is equivalent to constancy
of $\dim(W_0 \cap W_1)$.

Another phenomenon, also connected with spectral behavior, is that 
for {\sl harmonic} $u(t, \theta)$, $\ee_0= \{\,const\, \}$ does not
imply $u_0 = const$, but merely $u_0(t, \theta) = v_0 + v_1 t$ with $v_0,
v_1 \in \ee_0$, so that inequality (5.3) in not true. This leads to much more
complicated bubbling with energy loss for harmonic and
harmonic-type maps, compare [S-U], [P-W], [Pa].

Note also that if $W_0$ and $W_1$ are {\sl affine} totally real subspaces of
$\cc^n$, then inequality (5.3) for holomorphic $u: \Theta(0,3)\to \cc^n$ with
boundary condition (5.2) is in general not true. An easy example is a natural
imbedding $u:\Theta(0,3) \hook \cc$, with nonconstant component $u_0 \equiv
u$ and with $\int_{\Theta(0,1)} |du|^2 = \int_{\Theta (1,2)} |du|^2 = \int_{
\Theta(2,3)} |du|^2 \not=0$. In general, those are $n$-dimensional totally
real affine planes $W_0$ and $W_1$ in $\cc^n$ with empty intersection. This
can happen if $W_0$ and $W_1$ are parallel or {\sl skew}. The later means
that the corresponding vector spaces $V_0$ and $V_1$ ($W_i = V_i + w_i$ for
some $w_i \in \cc^n$) are different. In both cases the intersection $V_0 \cap
V_1$ is not zero, since otherwise $W_0 \cap W_1 \not= \emptyset$ by
dimension argumentation.

\smallskip
The considerations above show which properties should be controlled to obtain
a reasonable statement in the nonlinear case.

\state Definition 5.6. \sl Let $X$ be a manifold with a Riemannian metric
$h$, $f_0 : W_0 \to X$ and $f_1 : W_1 \to X$ immersions, and $A_0 \subset
W_0$, $A_1 \subset W_1$ subsets. We say that {\it $f_0 : W_0 \to X$ and $f_1
:W_1 \to X$ are $h$-uniformly transversal along $A_0$ and $A_1$ with
parameters $\delta >0$ and $M$} if for any $x_0 \in A_0$ and $x_1 \in A_1$
one of the following conditions hold:

\smallskip
\sli $\dist_h(f_0(x_0), f_1(x_1)) > \delta$;

\slii there exists $x'_i \in A_i$ with $f_0( x'_0) =f_1(x'_1)$ and such that
$$
\dist_h(x_0,x'_0) + \dist_h(x_1,x'_1) \le M\,\dist_h(f_0(x_0), f_1(x_1)).
$$\rm

\noindent
{\bf Remark.} Roughly speaking, the condition excludes appearance of points
where $W_0$ and $W_1$ are "asymptotically parallel or skew" and ensures us 
existence of a uniform lower bound for the angle between $W_0$ and $W_1$.

\medskip
Let now $X$ be a manifold with a Riemannian metric $h$, $J^*$ a continuous
almost complex structure on $X$, $f_0 : W_0 \to X$ and $f_1 : W_1 \to X$
immersions, and $A\subset X$, $A_0 \subset W_0$, $A_1 \subset W_1$ subsets.
Suppose that $J^*$ is $h$-uniformly continuous on $A$, $df_i: TW_i \to TX$
are $h$-uniformly continuous on $A_i$, and that $f_i$ are $h$-uniformly
transversal along $A_i$ with parameters $\delta= \delta(f_0,f_1)>0$ and $M=
M(f_0,f_1)$.

\medskip
\state Lemma 5.4. \it There exist constants $\eps^b_2 =\eps^b_2 (\mu_{J^*},
\, f_0,\, f_1,\, \delta,\, M) >0$ and $\gamma^b = \gamma^b(\mu_{J^*},\, f_0,\,
f_1,\, \delta,\, M) <1$ such that for any continuous almost complex $\ti J$
with $\norm{\ti J-J^*}_{L^\infty(A)}<\eps_2^b$, any immersions $\ti f_i :W_i
\to X$ with $\dist(\ti f_i, f_i)_{C^1(A_i)} <\eps_2^b$, and any $\ti J
$-holomorphic map $u \in C^0 \cap L^{1,2}( \Theta(0,5), X)$ with $u(
\Theta(0,5)) \subset A$, $u|_{\d_i\Theta(0,5)} = f_i \scirc u^b_i$ for some
continuous $u^b_i : \d_i\Theta(0,5) \to A_i \subset W_i$ the conditions

\itemitem\sli $\norm{ du}_{ L^2(\Theta_i)} <\eps_2^b$;

\itemitem\slii $\ti f_i: W_i \to X$ are $h$-uniformly transversal along
$A_i$ with the same parameters $\delta$ and $M$

\noindent
imply the estimate
$$
\norm{du}^2_{L^2(\Theta_3)}\le \msmall{\gamma^b \over 2} \cdot
\left( \norm{du}^2_{L^2(\Theta_2)} +  \norm{du}^2_{L^2(\Theta_4)} \right).
$$\rm

\state Proof. Suppose the statement of the lemma is false. Then there should
exist a sequence of continuous almost complex structures $J_k$ with $\norm{
J_k- J^*} _{L^\infty (A)} \lrar 0$, a sequence of immersions $f_{k,i}: W_i
\to X$ with $f_{k,i} \lrar f_i$ in $C^1(A_i)$, such that $f_{k,i} :W_i \to X$
are $h$-uniformly transversal with the same parameters $\delta$ and $M$, and
a sequence of $J_k$-holomorphic maps $u_k \in C^0\cap L^{1,2}(\Theta(0,5),
X)$ with $u_k(\Theta(0,5)) \subset A$ and $u_k |_{\d_i\Theta(0,5)} = f_{k,i}
\scirc u^b_{n,i}$ for some continuous $u^b_{n,i} : \d_i\Theta(0,5) \to A_i
\subset W_i$, such that $\norm{du_k}^2_{L^2(\Theta(0,5))} \lrar 0$ and
$$
\norm{du_k}^2_{L^2(\Theta_3)}\ge \msmall{\gamma_k \over 2} \cdot
\left( \norm{du_k}^2_{L^2(\Theta_2)} +  \norm{du_k}^2_{L^2(\Theta_4)} \right)
$$
with $\gamma_k = 1- 1/k$. {\sl Lemmas 4.1} and {\sl 5.1} provide that in this
case $\diam_h(u_k( \Theta(1,4))) \allowbreak \lrar 0$.

Since $f_{k,i} :W_i \to X$ are $h$-uniformly transversal with the same
parameters $\delta$ and $M$, there should exist sequences $x_k \in A$, $x_{k,
0} \in A_0$, and $x_{k,1} \in A_1$ such that $x_k = u_0(x_{k,0}) =u_1(x_{k,
1})$ and $u_k(\Theta(1,4)) \subset B(x_k, r_k)$ with $r_k \lrar 0$. The $h
$-uniform continuity of $J^*$ implies that there exist $C^1$-diffeomorphisms
$\phi_k: B(x_k, r_k) \to B(0, r_k) \subset \cc^n$ with $\norm{J_k - \phi^*_k
J\st}_{L^\infty(B(x_k, r_k))} + \norm{h -\phi^*_k h\st} _{L^\infty(B(x_k,
r_k))} \lrar 0$.

Using $\phi_k$ we transfer our situation into $B(0, r_k) \subset \cc^n$ and
rescale it. Namely we set $\alpha_k \deff \norm{du_k}_{L^2(\Theta_3)}$ and
define diffeomorphisms $\psi_k \deff {1\over \alpha_k} \scirc \phi_k: B(x_k,
r_k) \to B(0, R_k) \subset \cc^n$ with $R_k \deff \alpha_k\inv \cdot r_k$.
Note that by {\sl Lemmas 4.1} and {\sl 5.1} we have $\alpha_k = \norm{du_k}
_{L^2(\Theta(2,3))} \le C \diam_h(u_k( \Theta(1,4))) \le C' r_k$, so that
$R_k$ are uniformly bounded from below.

In $B(0, R_k)$ we consider Riemannian metrics $h_k \deff \alpha^{-2}_k \cdot
\psi_{k\,*}h_k$ (\ie pushed forward and $\alpha^{-2}_k$-rescaled $h_k$),
almost complex structures $J^*_k \deff \psi_{k\,*} J_k$, and $J^*_k
$-holomorphic maps $u^*_k \deff \psi_k \scirc u_k : \Theta(1,4) \to
B(0,R_k)$. Note that here we consider $h$ as a metric tensor, so multiplying
$h$ by $\alpha^{-2}$ we increase $h$-norms and $h$-distances in $\alpha^{-1}$
and not in $\alpha^{-2}$ times.

Then $\norm{du^*_k}_{L^2(\Theta_3, h_k)}=1$, $\norm{du^*_k}^2_{L^2(\Theta_2,
h_k)} + \norm{du^*_k}^2_{L^2(\Theta_4, h_k)} \le {2k\over k-1}$, and
$$
\norm{J^*_k - J\st}_{L^\infty(B(0, R_k), h_k)} =
\norm{J_k - \phi^*_kJ\st}_{L^\infty(B(x_k, r_k), h)} \lrar 0.
$$
The last equality uses an obvious relation
$$
\msmall{|F(\xi)|_{\alpha^{-2} \cdot h} \over |\xi|_{\alpha^{-2} \cdot h} }
= \msmall{\alpha\inv \cdot |F(\xi)|_h \over \alpha\inv \cdot |\xi|_\cdot h }
= \msmall{|F(\xi)|_h \over |\xi|_h }
$$
for any linear $F: T_xX \to T_xX$ and $\xi\not = 0\in T_xX$. In a similar way,
we also obtain $\norm{h_k - h\st}_{L^\infty(B(0, R_k), h_k)}\lrar 0$.

Going to a subsequence, we may additionally assume that the tangent spaces
$d\psi_k \scirc df_i (T_{x_{k,i}} W_i)$, $i=1,2$, converge to some spaces
$W_i^* \subset \cc^n$. Since $W_i$ are uniformly totally real, $W_i^*$ are
also totally real linear subspaces in $\cc^n$. Since the maps $df_i : TW_0
\to TX$ are uniformly continuous on $A_i \subset W_i$, $f_{k,i} \lrar f_i$ in
$C^1(A_i)$, and since $r_k \lrar 0$, the images $W^*_{k,i} \deff \psi_k
\scirc f_{k,i}(B_{W_i}(x_{k,i}, r_k))$ of the balls $B_{W_i}(x_{k,i}, r_k)
\subset W_i$ are imbedded submanifolds of $\cc^n$ with $0 \in W^*_{k,i}$,
which converge to $W^*_i$ in Hausdorff topology. Moreover, we can consider
$W^*_{k,i}$ as graphs of maps $g_{k,i}$ from subdomains $U_{k,i} \subset
W^*_i \cap B(0,R_k)$ to $W^*_i {}^\perp$ and for any fixed $R \le \inf
\{R_k\}$ the restrictions $g_{k,i} |_{W^*_i \cap B(0,R)}$ converge to zero
map from $W^*_i \cap B(0,R)$ to $W^*_i {}^\perp$.

The apriori estimates for the maps $u^*_k : \Theta(1,4) \to \cc^n$ imply
that for any $p<\infty$ the maps $u^*_k$ converge in weak- $L^{1,p}$-topology
to some $J\st$-holomorphic map $u^*: \Theta(1,4) \to \cc^n$. Further, since
$u^*_k$ satisfy totally real boundary conditions $u^*_k|_{\d_i\Theta(1,4)}
\subset W^*_{k,i}$, the same is true for $u^*$, \ie $u^*|_{\d_i\Theta(1,4)}
\subset W^*_i$. Nice behavior of $W^*_{k,i}$ shows that on $\Theta_3$ we
have also a strong convergence, and hence $\norm{du^*}_{L^2(\Theta_3)} = \lim
\norm{du^*_k}_{L^2(\Theta_3)} =1$. In particular, $u^*$ is not constant. On
the other hand, $\norm{du^*}^2_{L^2(\Theta_2)} + \norm{du^*}^2_{L^2(\Theta_4)}
\le \lim \norm{du^*_k}^2_{L^2(\Theta_2)} + \norm{du^*_k}^2_{L^2(\Theta_4)} \le
2$. The obtained contradiction with {\sl Lemma 5.3} shows that {\sl Lemma 5.4}
is true. \qed

\medskip
Let $X$, $h$, $J$, $A$, $f_i: W_i \to X$, $A_i$, and the constant $\eps^b_2$
and $\gamma^b$ be as in {\sl Lemma 5.4}. Suppose that $\ti J$ is a continuous
almost complex structure on $X$ with $\norm{\ti J-J}_{L^\infty(A)} <\eps^b_2$,
$\ti f_i: W_i \to X$ are totally real immersions with $\dist(\ti f_i, f_i)
_{C^1(A_i)} \le \eps^b_2$, such that $f_i$ are $h$-uniformly transversal
along $A_i$ with the same parameters $\delta$ and $M$ as $f_i: W_i \to X$.

\medskip
\state Corollary 5.5. {\it Let $u \in C^0\cap L^{1,2}(\Theta(0,l),X)$ be a
$\ti J$-holomorphic map, such that $u(\Theta(0,l))\subset A$, $u|_{\d_i
\Theta(0,l)} = \ti f_i \scirc u^b_i$ for some continuous $u^b_i: \d_i\Theta(
0,l) \to A_i \subset W_i$, and such that $\norm{du}_{ L^2(Z_k)}<\eps_2$ for
any $k=1,\ldots,l$.

Let $\lambda_b>1$ be the (uniquely defined) real number with $\lambda_b 
={\gamma^b \over 2}(\lambda_b^2+ 1)$. Then for $2\le k\le l-1$ one has
$$
\norm{du}^2_{L^2(\Theta_k)}
\le \lambda_b^{-(k-2)} \cdot \norm{du}^2_{L^2(\Theta_2)} +
\lambda_b^{-(l-1-k)} \cdot \norm{du}^2_{L^2(\Theta_{n-1})}.
\eqno(5.?9)
$$
}

\state Proof. The same as in {\sl Lemma 4.5}. \qed

\medskip
Immediate corollary of this estimate is a lower bound of energy on
nonconstant ``infinite strip''.

\smallskip
\state Lemma 5.6. \it In the setting of {\sl Corollary 5.5}, let
$u \in C^0\cap L^{1,2}(\Theta(-\infty, +\infty),X)$ be a nonconstant $\ti
J$-holomorphic map, such that $u(\Theta(-\infty, +\infty))\subset A$ and
$u|_{\d_i\Theta(0,l)} = \ti f_i \scirc u^b_i$ for some continuous $u^b_i :
\d_i\Theta(-\infty, +\infty) \to A_i \subset W_i$. Then $\norm{du}_{ L^2(
\Theta_k)} > \eps^b_2$ for some $k$. In particular, $\norm{du}_{ L^2(
\Theta(-\infty, +\infty))} > \eps^b_2$.\rm

\state Proof. {\sl Corollary 5.5\/} shows that if $\norm{du}_{L^2(
\Theta_k)} \le \eps^b_2$ for all $k$, then $\norm{du}_{ L^2(\Theta_k)}=0$,
\ie $u$ is constant. \qed

\medskip
Another consequence of {\sl Corollary 5.5} is a generalization of the Gromov's
result about removability of boundary point singularity, see [G]. An
important improvement is the fact that the statement remains valid also when
one has {\sl different} boundary conditions on the left and on the
right from a singular point. One can see such a point $x$ as a {\sl corner
point} for the corresponding pseudoholomorphic curve. Typical examples appear
in symplectic geometry where one takes Lagrangian submanifolds as boundary
conditions.

\smallskip
Define the punctured half-disk by setting $\check\Delta^+ \deff \Delta^+\bs
\{0\}$. Define $I_- \deff (-1,0) \subset \d\check\Delta^+$ and
$I_+ \deff (0,+1) \subset \d\check\Delta^+$.

\smallskip
\state Corollary 5.7. {\sl (Removal of boundary point singularities).
\it Let $X$ be a manifold with a Riemannian metric $h$, $J$ a continuous
almost complex structure, $f_i : W_i \to X$, $i=1,2$, totally real immersions,
and $A_i \subset W_i$ subsets. Let $u:(\check\Delta^+, J\st) \to(X,J)$ be a 
pseudoholomorphic map. Suppose that

\item\sli $J$ is uniformly continuous on $A \deff u(\check\Delta)$ \wrt
$h$, and closure of $A$ is $h$-complete;

\item\slii $u$ satisfies boundary conditions of the form $u|_{I_+} = f_0
\scirc u^b_+$ and $u|_{I_-} = f_1 \scirc u^b_-$ with come continuous $u^b_+ :
I_+ \to A_0 \subset W_0$ and $u^b_- : I_- \to A_1 \subset W_1$;

\item\sliii $f_i$ are $h$-uniformly totally real on $A_i$ and $h$-uniformly
transversal along $A_i$;

\item\sliv there exists $k_0$, such that for all half-annuli $R^+_k\deff \{
z\in \Delta^+ :{1\over e^{\pi(k+1)}}\le | z| \le {1\over e^{\pi k}}\}$ with
$k\ge k_0$ one has $\norm{du}^2_{L^2(R^+_k)}\le \eps^b_2$, $\eps^b_2$ being
from {\sl Lemma 5.4}.

\smallskip\noindent
Then $u$ extends to the origin $0\in \Delta^+$ as an $L^{1,p}$-map for some
$p>2$.
}

\state Proof. Using the holomorphic map $\exp : \Theta(0,\infty) \to
\check\Delta^+$, $\exp(\theta + \isl t) \deff e^{\pi(-t + \isl\theta)}$, we
can reduce our situation to the case of pseudoholomorphic map $u^* \deff u
\scirc \exp$ form ``infinite strip'' $\Theta(0,\infty)$. By {\sl Corollary
5.5}, for $k\ge k_0$ we obtain estimate $\norm{ du^*} _{L^2(\Theta_k)} \le
\lambda_b^{-(k-k_0)/2} \norm{ du^*}_{L^2(\Theta_{k_0})}$ with some $\lambda_b
>1$. This is equivalent to the estimate $\norm{ du} _{L^2(R^+_k)} \le
\lambda_b^{-(k-k_0)/2} \norm{ du}_{L^2(R^+_{k_0})}$. {\sl Lemmas 4.1} and
{\sl5.1} and scaling property of $L^p$-norms provide the estimate
$$
\norm{ du}_{L^p(R^+_k)} \le C e^{-k(\log\lambda_b/2 + \pi(2/p-1))}
$$
Thus $du \in L^p(\Delta^+)$ for any $p$ with $\log\lambda_b/2 > \pi(1-2/p)$,
which means $p < {4\pi \over 2\pi - \log \lambda_b}\cdot$ \qed

\smallskip
\state Remark. Unlike the ``interior" and smooth boundary cases, it is 
possible
that the map $u$ as in {\sl Corollary 5.7\/} is not $L^{1,p}$-regular in the
neighborhood of ``corner point'' $0\in \Delta^+$ for some $p>2$. For example,
the map $u(z) = z^\alpha$ with $0<\alpha <1$ satisfies totally real boundary
conditions $u(I_+) \subset \rr$, $u(I_-) \subset e^{\alpha\pi \isl}\rr$ and
is $L^{1,p}$-regular only for $p< p^* \deff {2\over1-\alpha}\cdot$

\medskip
As in the ``interior" case, for the proof of the boundary compactness theorem 
we shall need a description of a convergence of a sequence of ``long strip''. 
Let $X$ be a manifold with a Riemannian metric $h$, $J$ a continuous almost
complex structure, $A\subset X$ a closed $h$-complete subset, such that $J$
is $h$-uniformly continuous on $A$, and let $\{J_n\}$ be a sequence of almost
complex structures converging $h$-uniformly on $A$ to $J$. Let also $f_0: W_0
\to X$ and $f_1: W_1\to X$ be immersions, $A_i \subset W_i$ subsets, such
that $df_i$ are uniformly $h$-uniformly totally real on $A_i$ and $f_i$ are
$h$-uniformly transversal along $A_i$. Let $f_{n,i} : W_i \to X$ be totally
real immersions, which $C^1$-converge to $f_i$ on $A_i$, such that $f_{n,0}$
and $f_{n,1}$ are $h$-uniformly transversal along $A_i$ with uniform in $n$
parameters $\delta$ and $C^*$. Finally, let $\{l_n\}$ be a sequence of
integers with $l_n\to \infty$, and $u_n: \Theta(0,l_n) \to X$ a sequence of
$J_n$-holomorphic maps, satisfying boundary conditions $u_n|_{\d_i\Theta(0,
l_n)} = f_{n,i}\scirc u^b_{n,i}$ with some continuous $u^b_{n,i}:
\d_i\Theta(0,l_n) \to A_i \subset W_i$.

\state Lemma 5.8. {\it In the described situation, suppose additionally that
$u_n(\Theta(0,l_n))\subset A$ and $\norm{du_n} _{L^2 (\Theta_k)}\le \eps^b_2$
for all $n$ and $k\le l_n$. Take a sequence $k_n\to \infty$ such that
$k_n<l_n-k_n\to \infty$. Then:

\item{\sl1)} $\norm{du_n}_{L^2(\Theta(k_n,l_n-k_n))}\to 0$ and
$\diam\bigl(u_n (\Theta(k_n,l_n-k_n))\bigr)\to 0$.

\item{\sl2)}
There is a subsequence $\{ u_n \}$, still denoted $\{ u_n \}$, such that both
$u_n |_{\Theta(0,k_n)}$ and $u_n |_{\Theta(k_n,l_n)}$ converge in $L^{1,p}
$-topology on compact subsets in $\check\Delta^+ \cong \Theta(0, +\infty)$ to a
$J^*$-holomorphic maps $u^-_\infty$ and $u^+_\infty$. Moreover, both
$u^+_\infty$ and $u^-_\infty$ extend to origin and $u^+_\infty(0)=
u^-_\infty(0)$.
}

\medskip
Let us turn to the Gromov compactness theorem for curves with boundary on
totally real submanifolds. To give a precise statement we need to modify the
definition of the Gromov convergence ({\sl Definition 1.6}). The reason to do 
it is the following. Considering {\sl open} curves $C_n$ with changing complex
structures, we want to fix some kind of a common ``neighborhood of infinity"
$i_n: C^*\hookrightarrow C_n$ of every $C_n$. Thus we can imagine that all
changes of complex structure take place "outside of infinity", \ie in
relatively compact part $C_n \bs i_n(C^*) \Subset C_n$. This is done to
insure that $C_n$ do not approach to infinity in an appropriate moduli space.

On the other hand, it is more natural to consider curves $(\barr C_n, u_n)$
with totally real boundary conditions as compact objects without ``infinity".
In fact, in this case the behavior of $u_n$ near the boundary $\d C_n$ can be
controlled. The obtained apriori estimates near ``totally real boundary''
can be viewed as a part of such a "control". So for curves with totally real
boundary conditions we can hope to extend the Gromov convergence up to
boundary.

Further, as in the ``inner case'', an appropriate modification of the Gromov
convergence in this case should allow boundary bubbling and appearance of
boundary nodes. This means, however, that the structure of the boundary can
change during approach to the limit curve and cannot be considered as fixed.
Instead of it one should fix a type of boundary conditions. We shall
consider the following general situation.

Let $u_n: \barr C_n \to X$ be a sequence of stable $J_n$-holomorphic over $X$
with parametri\-zations $\delta_n: \barr\Sigma \to \barr C_n$. Let also
$\bfbeta =\{ \beta_i \}_{i=1}^m$ be a collection of arcs $\beta_i$ in $\d
\Sigma$, $\{ W_i\} _{i=1}^m$ a collection of real $n$-dimensional manifolds,
$f_{n,i} : W_i \to X$ a sequence of totally real immersions and $u^b_{n,i}:
\beta_{n,i} \to W_i$ a sequence of continuous maps from $\beta_{n,i} \deff
\delta_n( \beta_i)$. Assume that $\cup_{i=1}^m \beta_i= \d\Sigma$ and that
the interiors of $\beta_i$ are mutually disjoint and do not intersect the
pre-images of boundary nodal points of $C_n$. Then ${\mib W}_n \deff \{ (W_i,
f_{n,i} )\} _{i=1}^m$ are totally real boundary conditions on $(\barr C_n,
u_n)$ of the same type $\bfbeta$.

\state Definition 5.6. {\sl In the situation above we say that the sequence
of boundary conditions ${\mib W}_n$ of the same type $\bfbeta$ {\it
converges $h$-uniformly transversally} to $J^*$-totally real boundary
conditions ${\mib W}$ on subsets $A_i \subset W_i$ if

\sli ${\mib W}= \{ (W_i, f_i) \}_{i=1}^m$ where $f_i: W_i \to X$ are
$J^*$-totally real immersions;

\slii $f_{n,i}$ converge to $f_i$ in $C^1$-topology and this
convergence is $h$-uniform on $A_i$;

\sliii for any $n$ immersions $\{f_{n,i}\}_{i=1}^m$ are mutually
$h$-uniformly transversal along $A_i$ with parameters $\delta>0$ and $M$, and
this parameters are independent of $n$.
}

\smallskip
Note that the condition \sliii implies that the limit immersions $f_i$ are
also mutually $h$-uniformly transversal along $A_i$ with the same parameters
$\delta>0$ and $M$.

\state Definition 5.7. {\sl We say that the sequence $(\barr C_n, u_n)$ {\it
converges up to boundary} to a stable $J^*$-holomorphic curve $(\barr
C_\infty, u_\infty)$ over $X$ if the parametrizations $\sigma_n: \barr\Sigma
\to \barr C_n$ and $\sigma_\infty: \barr\Sigma \to \barr C_\infty$ can be
chosen in such a way that the following holds:

\sli $u_n\scirc \sigma_n$ converges to $u_\infty\scirc \sigma_\infty$ in
$C^0( \barr\Sigma, X)$-topology;

\slii if $\{ a_k \}$ is the set of the nodes of $C_\infty$ and $\{ \gamma_k
\}$, $\gamma_k \deff \sigma_\infty\inv(a_k)$ are the corresponding circles
and arcs in $\barr\Sigma$, then on any compact subset $K\comp \barr\Sigma \bs
\cup_k\gamma_k$ the convergence $u_n\scirc \sigma_n\to u_\infty\scirc
\sigma_\infty$ is $L^{1,p}(K, X)$ for all $p< \infty$;

\sliii for any compact subset $K\comp \barr\Sigma \bs \cup_k\gamma_k$ there
exists $n_0=n_0(K)$ such that $ \sigma_n^{-1}(\{ a_k \}) \cap K= \emptyset$
for all $n\ge n_0$ and complex structures $\sigma_n^*j_{C_n}$ smoothly
converge to $\sigma_\infty^*j_{C_\infty}$ on $K$.
}

\state Theorem 5.9. {\it Fix a metric $h$ on $X$, an $h$-complete subset
$A \subset X$, and subsets $A_i \subset W_i$. Suppose that:

\item{\sl a)} $J_n$ are continuous almost complex structures on $X$,
converging $h$-uniformly on $A$ to a continuous almost complex structure
$J^*$;

\item{\sl b)} $u_n(C_n) \subset A$ and $\area [u_n (C_n)]\le M$ with a
constant $M$ independent of $n$;

\item{\sl c)} ${\mib W}_n \deff \{(W_i, f_{n,i}) \}_{i=1}^m$ are totally real
boundary conditions of the same type $\bfbeta = \{ \beta_i \}_{i=1}^m$, such
that ${\mib W}_n$ converge $h$-uniformly transversally to a boundary condition
${\mib W}= \{(W_i, f_i) \}_{i=1}^m$ on subsets $A_i \subset W_i$;

\item{\sl d)} immersions $f_i: W_i \to (X,J^*)$ are $h$-uniformly totally
real along $A_i$;

\item{\sl e)} there exist maps $u^b_{i,n}: \beta_i \to  A_i\subset W_i$,
realizing boundary conditions ${\mib W}_n$.

\smallskip
Then there exits a subsequence of $\{( \barr C_n, u_n )\}$, still denoted
$\{( \barr C_n, u_n )\}$, and para\-metri\-sations $\sigma_n: \barr \Sigma
\to \barr C_n$, such that $(C_n, u_n, \sigma_n)$ converges up to boundary to
a stable $J^*$-holomorphic curve $(\barr C_\infty, u_\infty, \sigma_\infty)$
over $X$.

If, in addition, $A_i \subset W_i$ are $f_i^*h$-complete, then the limit
curve $(\barr C_\infty, u_\infty)$ satisfies real boundary conditions $\mib
W$ with maps $u^b_i: \beta_i \to A_i \subset W_i$.
}

\medskip
Our main idea of the proof is to apply arguments used in the proof of
{\sl Theorem 1.1}. To do so we replace every
pair $(C_n, u_n)$ by a triple $(C^d_n, \tau_n, u^d_n)$ where $C^d_n$ is the
{\sl Schottky double} of $C_n$ with an antiholomorphic involution $\tau_n$
and $u^d_n: C_N^d \to X$ a $\tau_n$-invariant map. Then we shall change
all the constructions of the proof to make them $\tau_n$-invariant in an
appropriate sense. In particular, the convergence $(C^d_n, \tau_n, u^d_n)
\lrar (C^d_\infty, \tau_\infty, u^d_\infty)$ will be equivalent to the
convergence $(C_n, u_n) \lrar (C_\infty, u_\infty)$.

\smallskip
We start with construction of the Schottky double of a nodal curve $\barr C$
with boundary. Take two copies $\barr C^+ \equiv \barr C$ and $\barr C^-$ of
$\barr C$. Equip $\barr C^-$ with the opposite complex structure, so that the
identity map $\tau: \barr C^+ \to \barr C^-$ becomes now antiholomorphic.
Glue $\barr C^+$ and $\barr C^-$ together along their boundaries identifying
$\d C^+$ and $\d C^-$ by means of the identity map $\tau: \d C^+ \buildrel
\cong \over \lrar \d C^-$. The union $C^d \deff \barr C^+ \cup_{\d C} \barr
C^-$ possesses the unique structure of a closed nodal curve, which is 
compatible with imbeddings $\barr C^\pm \hook C^d$. The boundary $\d C$ 
becomes the fixed point set of $\tau$.

The map $\tau$ induces an antiholomorphic involution of $C^d$ which we also
denote by $\tau$. We call the obtained curve $C^d$ the {\sl Schottky double} 
of $\barr C$. Note that every boundary nodal point $a_i\in \d C$ defines a 
$\tau$-invariant nodal point $a_i$ on $C^d$, whereas an inner nodal point $b_i
\in C$ defines a pair of nodal points $b_i^\pm$ on $C^d$ interchanged by 
$\tau$. If $\sigma: \barr\Sigma \to \barr C$ is a parametrization of $\barr C$,
then we obtain in an obvious way the double $\Sigma^d$ with the involution 
$\tau: \Sigma^d \to \Sigma^d$ and the parametrization $\sigma^d: \Sigma^d 
\to C^d$ compatible with the involutions.

\state Remark.
The introduced notation $C^d$ for the {\sl Schottky double} of a nodal curve
$\barr C$ with boundary coincides with the one for the {\sl holomorphic
double}, used in {\sl Section 2}. Since in this section only the Schottky
double is considered, this should not lead to confusion.

\smallskip
Suppose additionally that an almost complex structure $J$ on $X$ and a
$J$-holo\-mor\-phic map $u: \barr C \to X$ are given. Suppose also that the
curve $(\barr C,u)$ satisfies the totally real boundary conditions $\mib W$ of
type $\bfbeta$. In particular, $\bfbeta$ defines the certain system of arcs
$\{\beta_i\}$ on $\d C$. In order to take into account the type of boundary
conditions, we fix the ends of $\beta_i$ which are not boundary nodal points
of $\barr C$, and declare these points to be marked points of $C^d$. Note that
these ones and the nodal points are the only ``corner'' points of $(\barr C,
u)$. The latter means, that in a neighborhood of these points the map $u$ can
be $L^{1,p}$-smooth not for all $p< \infty$. The example in {\sl Remark} after
{\sl Corollary 5.7} explains the notion ``corner point''. Considering the
Schottky double, we shall always equip $C^d$ with this set of marking points.
Note also that every boundary circle of $\barr C$ contains at least one nodal
or marked point as above.

For $(\barr C,u)$ as above, we extend the $J$-holomorphic map $u: \barr C \to
X$ to a map $u^d: C^d \to X$ by setting $u^d(x) \deff u(\tau(x))$ for $x\in
C^-$. By the construction, $u^d$ is $\tau$-invariant, $u^d \scirc \tau= u^d$,
but $u^d$ {\sl is not $J$-holomorphic} (with the only trivial exception $u
\equiv const$). However, the analysis already done in this section provides
necessary $L^{1,p} $-estimates for $u^d$, at least for some $p^*>2$.

\smallskip
In the situation of {\sl Theorem 5.9}, such an exponent $p^*>2$ can be chosen
the same for all curves $(\barr C_n, u_n)$, it depends only on the topology
of $\barr C_n$ and the geometry of immersions $f_n: W_n \to X$. In
particular, every $u_n^d$ is continuous.

\smallskip
Next step of the proof is to find a $\tau_n$-invariant decomposition of
$C^d_n$ into pants. This implies that the corresponding graph $\Gamma_n$
becomes $\tau_n$-invariant. In the construction which follows we shall use
the fact that $\tau_n$ is an isometry on the union of the non-exceptional
components of $C^d_n$. This is provided by uniqueness of the intrinsic
metric.

\state Lemma 5.10. \it Let $C$ be a nodal curve with boundary, $\sigma: \barr
\Sigma \to \barr C$ a parametrization, and $\{ x_i \}_{i=1}^m$ a set of
marked points on boundary $\d C$. Let $C^d$ be the Schottky double of $C$
with the anti-holomorphic involution $\tau$.

Then there exists a $\tau$-invariant decomposition of $C^d\bs \mapo$ into
pants, such that the intrinsic length of corresponding boundary circles is
bounded by a constant $l^+$ depending only on genus $g$ of $\Sigma^d$ and
the number of marked points $m$.

Moreover, every {\sl short} geodesic appears as a boundary circle of some
pants of the decomposition.\rm

\state Remark. Recall (see {\sl Remark} on page 25) that a closed geodesic
$\gamma$ is called {\sl short} if $\ell(\gamma) < l^*$, where $l^*$ is the
universal constant $l^*$ with the following property: For any simple closed
geodesics $\gamma'$ and $\gamma''$ on the conditions $\ell(\gamma') <l^*$ and
$\ell( \gamma'') <l^*$ imply $\gamma' \cap \gamma'' = \emptyset$.

\state Proof. Since genus of the parameterizing real surface $\Sigma^d$ and
the number of marked points is fixed, we obtain a uniform upper bound on
possible genus and the number of marked points of non-exceptional components
of $C^d$, as also on the number of exceptional components. This implies that
there exists a decomposition of every non-exceptional component $C_i$ of
$C^d$ into pants $S_\alpha$, such that the intrinsic length of boundary
circles of $S_\alpha$ is bounded by the constant $l^+$ depending only on $g$
and $m$. The idea of the proof of our lemma is to show that the construction
of such a decomposition, given in [Ab], Ch.II, \S\.3.3, can be modified to
produce a $\tau$-invariant decomposition.

Let us first describe the construction itself, say, for a given smooth curve
$C^*$ with marked points $\{x_i\}$ of non-exceptional type. The procedure is
done inductively by choosing at every step a non-trivial simple closed
geodesic $\gamma_{J^*} \subset C^* \bs\mapo$, disjoint from already chosen
geodesic $\gamma_j$, $j<J^*$. Moreover, at every step there exists a geodesic
$\gamma_{J^*}$ as above whose intrinsic length is bounded by a constant
$l^+_{J^*}$ depending only on genus of $C^*$, the number of marked
points, and the maximum of the lengths of the already chosen geodesics
$\gamma_j$, $j<J^*$.

\smallskip
Take any non-exceptional component $C^d_i$ of $C^d$. Two cases can happen:
either $C^d_i$ is $\tau$-invariant, or $\tau \bigl( C^d_i \bigr)$ is another
component $C^d_{i'}$. These cases are distinguished by the property whether
$C^d_i$ intersects the boundary $\d C$ (first case) or not (second one).

The existence of $\tau$-invariant decomposition into pants for every pair of
non-exceptional components $C^d_i$ and $\tau(C^d_i) \not= C^d_i$ is obvious.
We choose an appropriate decomposition of $C^d_i$ and transfer it on $\tau(
C^d_i)$ by means of $\tau$.

\smallskip
It remains to consider the case of a $\tau$-invariant non-exceptional
component $C^d_i$.

Suppose that on some step we have already chosen a $\tau$-invariant set
$\{\gamma_1, \ldots, \allowbreak \gamma_{J^*-1}\}$ of simple disjoint
geodesics on $C^d_i \bs \mapo$. Take a simple geodesic $\gamma$ of the length
$\ell(\gamma) \le l^+_{J^*}$, where $l^+_{J^*}$ is the upper bound introduced
above. By the construction of the double $C^d$, the fixed point set of $\tau$
on $C^d_i$ is $C^d_i \cap \d C$ and is non-empty. Denote $C_i \deff C \cap
C^d_i$, so that $C^d_i \cap \d C= \d C_i$. Note that any boundary circle of
$C_i$ contains at least one marked point of $C^d_i$. Consequently, it has an
infinite length \wrt the intrinsic metric on $C^d_i \bs\mapo$. Thus the
chosen geodesic $\gamma$ can not lie on $\d C_i$. Only 3 cases can happen.

\smallskip\noindent
{\sl Case 1.} $\gamma$ is disjoint from $\d C_i$. Then $\gamma$ lies either in
$C_i$ or in $\tau (C_i)$. In any case, $\gamma\cap \tau(\gamma) = \emptyset$.
Thus we can set $\gamma_{J^*} = \gamma$ and $\gamma_{J^*+1}= \tau(\gamma)$,
getting the $\tau$-invariant set $\{\gamma_1, \ldots, \gamma _{J^*+1}\}$ of
simple disjoint geodesic. This will be the next 2 steps of our construction.

\smallskip\noindent
{\sl Case 2.} $\gamma \cap \d C_i \not=0$ and $\gamma$ is $\tau$-invariant.
We set $\gamma_{J^*}=\gamma$ and proceed inductively further. Note that in
this case $\gamma \cap \d C_i$ consists of 2 points, in which $\gamma$ is
orthogonal to $\d C_i$.

\smallskip\noindent
\line{\vtop{\hsize=.5\hsize\noindent
{\sl Case 3.} This time $\gamma \cap \d C_i \not=0$, but $\gamma \not=
\tau(\gamma)$. Define arcs $\gamma^+ \deff \gamma \cap \barr C_i$ and
$\gamma^- \deff \gamma \cap \tau(\barr C_i)$, the parts of $\gamma$ inside
and outside of $C_i$ (see Fig.~7). Consider the following free homotopy
classes of closed circles on $C^d_i$:

1) $[\ti\gamma_1] \deff [\gamma^+ \cup \tau(\gamma^-)]$;

2) $[\ti\gamma_2] \deff [\gamma^+ \cup \tau(\gamma^-)]$;

3) $[\ti\gamma_3] \deff [\gamma^- \cup \tau(\gamma^-)]$;

4) $[\ti\gamma_4] \deff [\gamma^+ \cup \tau(\gamma^-) \cup \gamma^- \cup
\tau(\gamma^+)]$.
}
\hss
\vtop{\xsize=.46\hsize\hsize=\xsize\nolineskip\rm
\putm[.47][-.01]{\ti\gamma_1}%
\putm[.02][.32]{\ti\gamma^+_2}%
\putm[.95][.42]{\ti\gamma^+_3}%
\putm[.15][.25]{\gamma^+}%
\putm[.61][.23]{\tau(\gamma^-)}%
\putm[.40][.17]{\ti\gamma^+_4}%
\noindent
\epsfxsize=\xsize\epsfbox{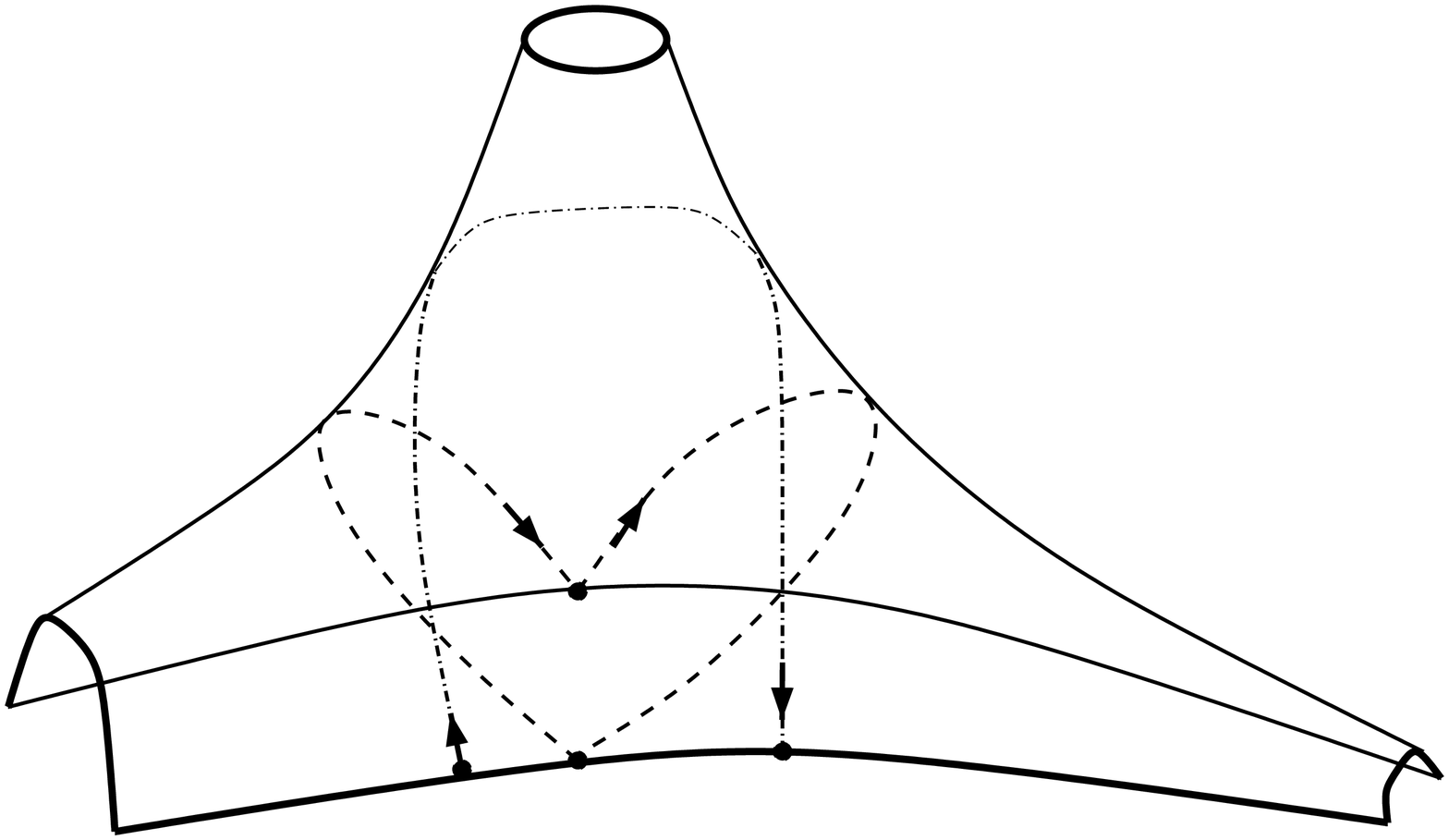}
\smallskip
\centerline{Fig.~7. Geodesics on $C_i$.}
}}

\vskip2pt\noindent
The last expression means that we move along corresponding arcs in the
prescribed order, as it is shown on Fig.~7. Note that only one part of
$C^d_i$ is drawn, namely $C_i$. The rest of the picture is symmetric \wrt
the involution $\tau$. Thus we can see only the half of geodesics in classes
$[\ti\gamma_i]$, $i=2,3,4$.

Each of classes $[\ti\gamma_k]$ either is represented by a closed geodesic or
corresponds to a wind around some marked point of $C^d_i$. To shorten
notations, we say in the last case that the class $[\ti\gamma_i]$ corresponds
to a marked point of $C^d_i$.

If one of the classes $[\ti\gamma_k]$, $k=1,2,3$, is represented by the
geodesic $\ti\gamma_k$, which is different and disjoint from the already
chosen geodesics $\gamma_j$, $j<J^*$, then we can set $\gamma_{J^*} =
\ti\gamma_k$. If $k=1$ we set also $\gamma_{J^*} = \ti\gamma_1$ and
$\gamma_{J^*+1} = \tau(\ti\gamma_1)$. Then we proceed inductively further.

\smallskip
To finish the proof it remains to consider the following situation: Under
conditions of {\sl Case 3}, each of the classes $[\ti\gamma_k]$, $k=1,2,3$,
either corresponds to a marked point, or is represented by a closed geodesic
$\ti\gamma_k$, which intersects or coincides with one from the already chosen
geodesics $\gamma_j$, $j<J^*$.

We claim that a proper intersection can not happen, \ie each class
$[\ti\gamma_k]$, $k=1,2,3$, either corresponds to a marked point, or is
represented by an already chosen geodesic $\gamma_j$, $j<J^*$. To show this
we note that $\gamma_j \cap \tau(\gamma) =\emptyset$ for all $j<J^*$.
Otherwise we could have a contradiction with the conditions $\gamma_j \cap
\gamma =\emptyset$ and $\tau$-invariance of the set of the geodesics
$\gamma_j$, $j<J^*$. Consequently, each class $[\ti\gamma_k]$ is represented
by a circle $\alpha_k \subset C^d_i \bs\mapo$, $k=1,2,3,4,$ with $\alpha_k
\cap \gamma_j = \emptyset$.

Now assume that the proper intersection of $\ti\gamma_k$ and some $\gamma_j$,
$j<J^*$, does have place. Let $\ell_k\deff \ell(\ti\gamma_k)$ be the intrinsic
metric of $\ti\gamma_k$. As in the proof of {\sl Lemma 2.2} construct the
annulus $A= \{ (\rho, \theta) : |\rho| < {\pi^2 \over \ell} \} \times \{0\le
\theta \le 2\pi\}$ with the metric $({\ell_k \over 2\pi} / \cos {\ell_k \rho
\over 2\pi})^2 (d\rho^2 + d\theta^2)$ and an isometric covering of $C^d_i
\bs\mapo$ by $A$, which sends the geodesic $\beta_k \deff \{ \rho=0\} \subset
A$ onto $\ti\gamma_k \subset C^d_i$. Find a lift of $\gamma_j$ to a geodesic
line $L_j \subset A$ with $L_j\cap \beta_k \not= \emptyset$, and a lift of the
circle $\alpha_k$ to a circle $\ti\alpha_k \subset A$ homotopic to $\beta_k$.
Then the intersection $L_j \cap \beta_k$ must consist of exactly one point,
and consequently, the homology intersection index $[L_j] \cdot [\beta_k]$ is
equal to $\pm1$. This would imply that $[L_j] \cdot [\ti\alpha_k]= [L_j] 
\cdot [\beta_k] \not =0$ and consequently $L_j \cap \ti\alpha_k \not= 
\emptyset$. But this would contradict to $\gamma_j \cap \alpha_k =\emptyset$.

Summing up, we see that in our situation we must have a picture of {\sl 
Fig.~7}. Namely, the both geodesics $\gamma$ and $\tau(\gamma)$ lie in a
$\tau $-invariant domain $\Omega$ on $C^d_i$ with 4 components of the
boundary. These components of $\d\Omega$ are either marked points or geodesics
corresponding to the classes $[\ti\gamma_1]$, $[\tau(\ti\gamma_1)]$, $[\ti
\gamma_2]$, $[\ti\gamma_3]$. Finally, every boundary circle of $\Omega$ is one
of the geodesics $\gamma_j$. We conclude, that the class $[\ti\gamma_4]$ is
represented by a $\tau$-invariant geodesic $\ti\gamma_4$. This geodesic can be
chosen at this step of construction of $\tau$-invariant decomposition of 
$C^d_i$ into pants.

Note that by construction for the intrinsic length of $\gamma_{J^*}$ we get
$\ell(\gamma_{J^*}) \le 2\ell(\gamma) \le 2l^+_{J^*}$. This means that in our
construction we do not lose control of the intrinsic length of chosen
geodesics. This provides the existence of a constant $l^+$ stated in the
lemma.

Finally, the definition of a {\sl short} geodesic provides that the geodesic
$\gamma$ in {\sl Case 3} above cannot be short. This implies that the set of
short geodesic on $C^d$ is disjoint. Since the involution $\tau$ is an
isometry, the set of short geodesic on $C^d$ is also $\tau$-invariant. Thus
in our construction of decomposition into pants we can start with this set of
geodesics. This shows the last statement of the lemma. \qed

\smallskip
\state Remark. To explain the meaning of {\sl Lemma 5.10}, let us consider
pants $S$ with a complex structure $J_S$ and an anti-holomorphic involution
$\tau$ acting on $S$. It is easy to see that only two types of such an
action, illustrated by Figs.\ 8\.a) and 8\.b), are possible.

\medskip\smallskip
\line{%
\vbox{\xsize=.47\hsize
\hsize=\xsize\nolineskip
\putm[.14][.11]{\gamma_1}%
\putm[.14][.60]{\gamma_2}%
\putm[.92][.37]{\gamma_3}%
\putm[.52][.36]{\beta}%
\putm[.52][.10]{S^+}%
\putm[.52][.57]{S^-}%
\noindent
\epsfxsize=\hsize\epsfbox{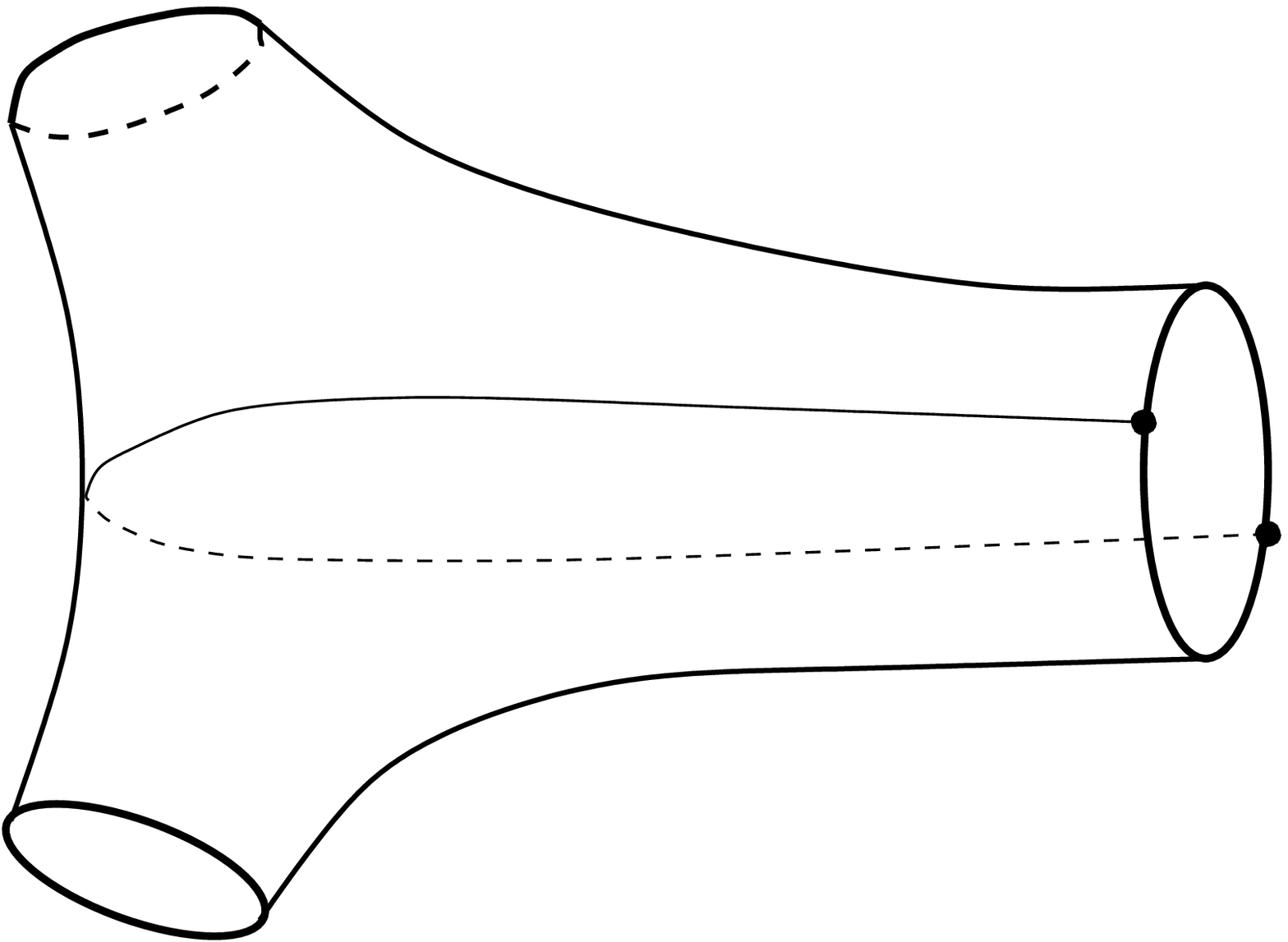}}%
\hfil
\vbox{\xsize=.47\hsize
\hsize=\xsize\nolineskip
\putm[.0][-.035]{\gamma_1}%
\putm[.28][.45]{\gamma_2}%
\putm[.90][-.045]{\gamma_3}%
\putm[.60][.22]{\beta_1}%
\putm[.52][.14]{\beta_2}%
\putm[.28][.12]{\beta_3}%
\noindent
\epsfxsize=\hsize\epsfbox{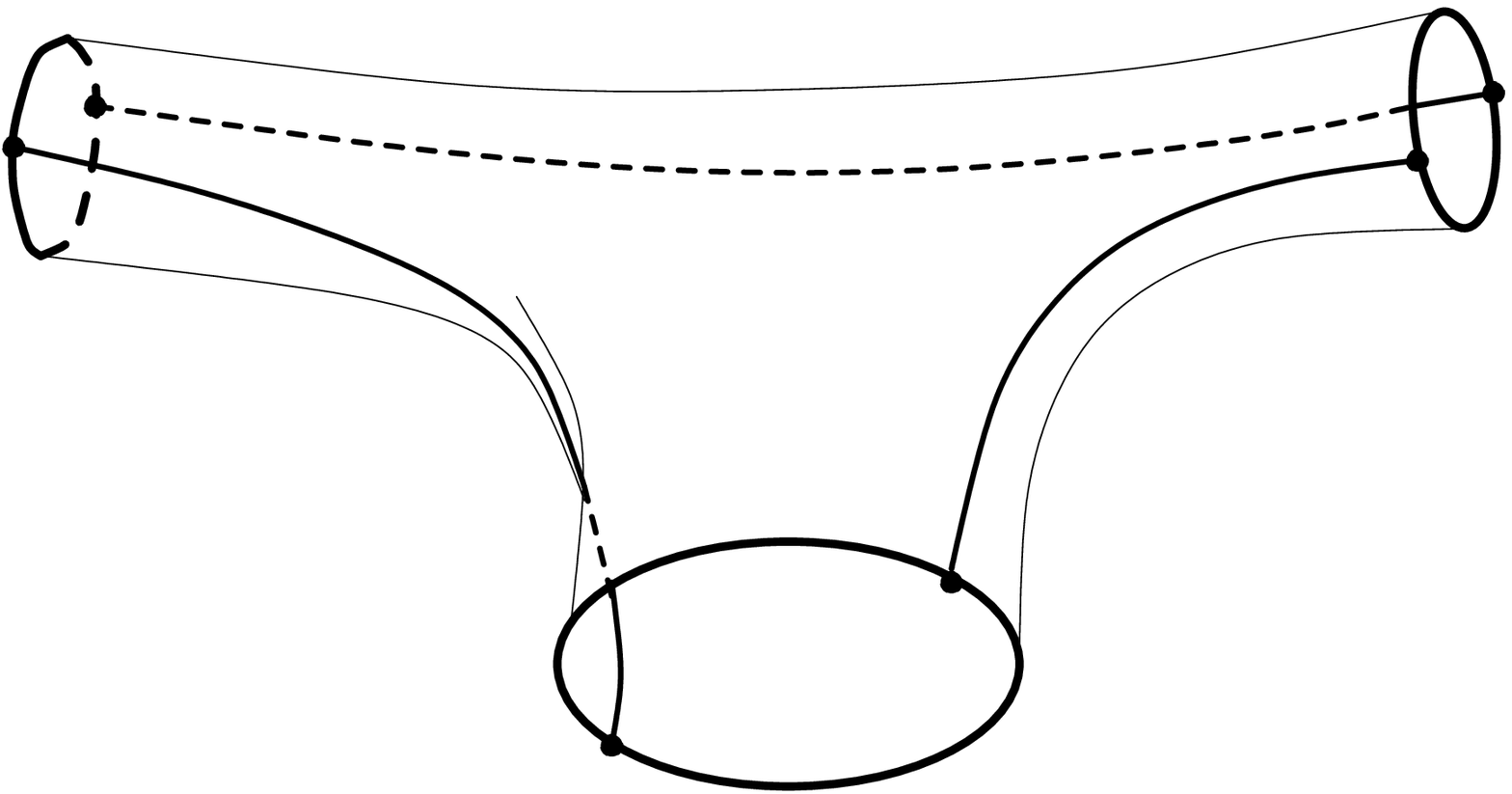}}}
\smallskip
\line{%
\vtop{\hsize=.47\hsize
\centerline{Fig.~8\.a)}
}%
\hfil
\vtop{\hsize=.47\hsize
\centerline{Fig.~8\.b)}
}%
}

\smallskip
\medskip
In the first case, Fig.~8\.a), the involution $\tau$ interchanges two
boundary components $\gamma_1$ and $\gamma_2$ of $S$ and leaves the third one
$\gamma_3$ invariant. The fixed point set $\beta$ of $\tau$ is a geodesic arc
with both ends on the $\tau$-invariant boundary component $\gamma_3$. This
case includes subcases when some boundary components of $S$ are not geodesics
but marked points (\ie punctures). In particular, if $\gamma_3$ is a marked
point, then the set $\beta$ is an (infinite) geodesic line with both ends
approaching to $\gamma_3$. The set $\beta$ divides $S$ into two parts, $S^+$
and $S^-$ (see Fig.~8\.a)), which are interchanged by $\tau$. Topologically,
each part $S^\pm$ is an annulus.

\smallskip
In the second case, Fig.~8\.b), all three boundary components $\gamma_1$,
$\gamma_2$, and $\gamma_3$ are invariant. The fixed point set of $\tau$
consists of geodesic arcs $\beta_1$, $\beta_2$, and $\beta_3$. These are the
shortest simple geodesics between $\gamma_2$ and $\gamma_3$, resp.\ $\gamma_3$
and $\gamma_1$, and resp.\ $\gamma_2$ and $\gamma_3$. If some boundary
component of $S$ is not a geodesic but a marked point, then corresponding
arcs have ends of infinite length approaching to this boundary component. The
arcs $\beta_k$, $k=1,2,3$, divide $S$ into two parts, $S^+$ and $S^-$ (see
Fig.~8\.b)), which are interchanged by $\tau$. In this case, each part
$S^\pm$ is topologically a disc.

\smallskip
We call pieces $S^\pm$ {\sl half-pants of first or second type} respectively.
Note that in both cases $\tau$-invariant arcs $\beta$ or $\beta_i$ are
orthogonal to corresponding boundary circles $\gamma_j$.

\medskip
Return to the situation of a nodal curve $\barr C$ with boundary and marked
points. Let $C^d$ be the Schottky double and $\tau$ the anti-holomorphic
involution. Suppose that $C^d \bs\mapo$ is non-exceptional. Use {\sl Lemma
5.10} and find a $\tau$-invariant decomposition into pants $C^d = \cup_j
S_j$. Set $S_j^+ \deff S_j \cap \barr C$. Then we obtain a decomposition
$\barr C= \cup_j S^+_j$, such that the pieces $S^+_j$ are either pants (which
means $S^+_j = S_j$), or half-pants of the first or the second type. This
decomposition is a suitable one for the situation of the Gromov convergence up
to boundary of curves with totally real boundary conditions. In particular, we
obtain arcs $\beta_{j,k}$ as $\tau$-fixed point sets of $S^+_j$, which define
a decomposition $\d C= \cup_{j,k} \beta_{j,k}$ of the boundary of $\barr C$.
The collection $\bfbeta' \deff \{ \beta_{j,k} \}$ of these arcs satisfies the
condition \sli of {\sl Definition 5.4}, but it can be different from the
collection $\bfbeta =\{ \beta_i \}$ which was given. The reason is that in
construction of the pants-decomposition $C^d = \cup_j S_j$ we can subdivide
original arcs $\beta_i \in \bfbeta$ into smaller pieces, so that every arc
$\beta_i \in \bfbeta$ is a union of arcs $\beta_{j,k}$ from $\bfbeta'$. This
means compatibility of $\bfbeta$ and $\bfbeta'$.

\medskip\smallskip
The next step is to establish a generalization of {\sl Theorem 4.2}. Assume 
that the hypothesis of {\sl Theorem 5.9} are fulfilled. For each curve $C_n$ 
denote by $C^d_n$ its Schottky double and by $\tau_n$ the corresponding 
involution.

\state Lemma 5.11. \it In the situation above, after passing to a
subsequence, there exist parametrizations $\sigma^d_n: \Sigma^d \to C^d_n$,
a finite covering $\calv$ of\/ $\Sigma^d$ by open sets $\{ V_\alpha \}$, and
a set $\{x^*_1, \ldots, x^*_m\}$ of marked points on $\Sigma$, such that the
conditions {\sl(a), (c)--(f)} of {\sl Theorem 4.2} and the following
additional conditions {\sl(b')} and {\sl(h\/)} are satisfied:

{\sl(b')} $\sigma_n\{x^*_1, \ldots, x^*_m\}$ is the set of marked points
on $C^d_n$ corresponding decomposition of the boundary $\d C_n$ into arcs
$\beta_{n,i}$; moreover, each such point $x^*_j$ lies in a single piece of
covering $V_\alpha$ which is a disc;

{\sl(h\/)} there exist an involution $\tau : \Sigma^d \to \Sigma^d$ which is
compatible with the covering $\calv$ and with parametrizations $\sigma^d_n$,
\ie $\calv$ is $\tau$-invariant and $\tau_n \scirc \sigma^d_n = \sigma^d_n
\scirc \tau$. In particular, each marked point $x^*_i$ of $\Sigma^d$ is fixed
by $\tau$. \rm

\state Remark. The condition {\sl(g)} of {\sl Theorem 4.2} is trivial in this
case since $C^d_n$ and $\Sigma^d$ are closed.

\state Proof. One can use the proof of {\sl Theorem 4.2} with minor
modifications. Note that starting points of that proof were the intrinsic
metric on non-exceptional components of nodal curves $C_n$ there and the
decomposition of $C_n$ into pants. Now the existence of a $\tau$-invariant
decomposition of the curves $C^d_n$ into pants is provided by {\sl Lemma
5.10}, whereas the $\tau$-invariance of the intrinsic metrics follows from
the fact, that any (anti)holomorphic isomorphism of curves with marked points
is an isometry \wrt the intrinsic metric. Thus the constructions of the proof
of {\sl Theorem 4.2} yield $\tau$-invariant objects. The condition {\sl(b')}
does not bring much difficulty. \qed

\medskip
Now we are ready to finish

\state Proof of Theorem 5.9. As it was mentioned, our main idea is to modify
the construction used in the proof of {\sl Theorem 1.1} to make them
$\tau$-invariant. Main work is already done: we have necessary apriori
estimates, the construction of a $\tau$-invariant pants-decomposition of the
double $C^d_n$ of the curve $\barr C_n$, and the appropriate covering $\calv$
of the real surface $\Sigma^d$ parametrizing the doubles $C^d_n$.

As in the proof of {\sl Theorem 1.1} we consider the curves $C_{ \alpha ,n}
\deff \sigma^d_n(V_\alpha)$. Due to presence of the involutions $\tau_n$ the
geometrical situation in now different. This involves new phenomena and needs
additional considerations and constructions. In particular, the pieces $C_{
\alpha ,n}$ are divided into 2 groups depending on that whether they are
disjoint from the boundary $\d C_n$ or intersect it. In the last case $C_{
\alpha ,n}$ is $\tau_n$-invariant. In this case we shall use the notation
$C^+_{ \alpha ,n} \deff C_{ \alpha,n} \cap C_n$ for the part of $C_{ \alpha
,n}$ lying in $C_n$. Besides, we denote $V^+_\alpha \deff V_\alpha \cap
\Sigma$. Then $V_\alpha$ appear as the union of domains $V^+_\alpha$ and
$\tau(V^+_\alpha)$, interchanged by $\tau$. Similar is true for $C_{ \alpha
,n}$.

\smallskip
To prove the theorem, we want to construct a refined covering $\wt\calv$ of
$\Sigma$ and refined parametrizations $\ti\sigma_n: \Sigma \to C_n$ such
that for every $V_\alpha \in \wt\calv$ the sequence $(C_{ \alpha ,n},
u_{ \alpha ,n})$ with $C_{ \alpha ,n} \deff \ti\sigma_n(V_\alpha)$
one the following convergence types holds:

\smallskip{\parindent=1.5\parindent\sl
\item{A\/$'$)} $C_{\alpha, n}$ are annuli of infinitely growing conformal
radii $l_n$ disjoint from $\d C_n$, and the conclusions of {\sl Lemma 3.7}
hold;

\item{A\/$''$)} $C_{\alpha, n}$ are $\tau_n$-invariant annuli of infinitely
growing conformal radii $l_n$ and the conclusions of {\sl Lemma 5.8} are
valid for $\Theta(0,l_n) \cong C^+_{\alpha, n} \deff C_{\alpha, n} \cap C_n$;

\item{B\/$'$)} every $C_{\alpha, n}$ is disjoint from $\d C_n$ and isomorphic
to the standard node $\cala_0 =\Delta \cup_{ \{0\} } \Delta$, such that the
compositions $V_\alpha \buildrel \sigma _{\alpha, n} \over \lrar C_{\alpha,
n} \buildrel \cong \over \lrar \cala_0$ define the same parametrisations of
$\cala_0$ for all $n$; furthermore, the induced maps $\ti u_{\alpha, n}:
\cala_0 \to X$ strongly converge;

\item{B\/$''$)} every $C_{\alpha, n}$ is $\tau_n$-invariant and $C^+_{\alpha,
n} \deff C_{\alpha, n} \cap C_n$ is isomorphic to the standard boundary node
$\cala^+ _0 =\Delta^+ \cup_{ \{0\} } \Delta^+$, such that for $V^+_\alpha
\deff V_\alpha \cap \Sigma$ the compositions $V^+_\alpha \buildrel \sigma^+
_{\alpha, n} \over \lrar C^+_{\alpha, n} \buildrel \cong \over \lrar
\cala^+_0$ define the same parametrisations of $\cala^+_0$ for all $n$;
furthermore, the induced maps $\ti u^+_{\alpha, n}: \cala^+_0 \to X$ strongly
converge;

\item{C)} the structures $\sigma_n^*j\vph_n \ogran_{V_\alpha}$ and the maps
$u_{\alpha, n}\scirc \sigma _{\alpha, n}: V_\alpha \to X$ strongly
converge.
}

\noindent
In the case {\sl B\/$''$)} the strong convergence of maps $\ti u^+_n: \cala^+
_0 \to X$ is the one in the $L^{1,p^*}$-topology for some $p^*>2$ {\sl up to
the boundary intervals containing the nodal point}. An equivalent requirement
is the usual $L^{1,p^*}$-convergence of the doubles $\ti u^d_n: C_{\alpha, n}
\to X$ on compact subsets of $C_{\alpha, n} \cong \cala_0$.

\smallskip
To obtain a desired refinement we use the same inductive procedure as in the
proof of {\sl Theorem 1.1}. To insure convergence near boundary $\d C_n$, we
take a new value for the constant determining the inductive step. We choose a
positive $\eps^b$, such that $\eps^b \le \eps$ and such that all apriori
estimates of this section are valid for maps with area\.\.$\le 3\eps^b$. This
will yield the convergence of type {\sl A)--C)} for sequences of curves with
totally real boundary conditions and with the upper bound $\eps^b$ on area.

In fact, essential modifications of constructions of {\sl Theorem 1.1} are
needed only if the covering piece $V_\alpha$ is $\tau$-invariant. Indeed, if
$V_\alpha$ is not $\tau$-invariant, then we can apply all the argumentations
and constructions used in {\sl Cases 1)--4)} in the proof of {\sl Theorem
1.1}, and then ``transfer'' them onto $\tau(V_\alpha)$ by means of $\tau$.
This gives the inductive step preserving $\tau$-invariance.

Hence, it remains consider the situation when the covering piece $V_\alpha$
is $\tau$-invariant. As in {\sl Theorem 1.1}, we must consider 4 cases:

{\sl\smallskip
Case 1$_b$): $C_{\alpha, n}$ have constant complex structure, different from
the one of the standard node;

\smallskip
Case 2$_b$): $C_{\alpha, n}$ are annuli of changing conformal radii $R_n$,
such that $R_n \to R <\infty$

\smallskip
Case 3$_b$): $C_{\alpha, n}$ are isomorphic to the standard node, so that
$C^+_{\alpha, n}$ are isomorphic to the standard boundary node $\cala_0^+$;

\smallskip
Case 4$_b$): $C_{\alpha, n}$ are annuli of infinitely growing conformal radii
$R_n$.
}

The subindex $(\cdot)_b$ indicates that we consider the cases when $V_\alpha$
intersects the {\it b}oundary of $\Sigma$. As it was mentioned the last
property is equivalent to the $\tau$-invariantness of $V_\alpha$. References
to {\sl Cases 1)--4) \it without} the subindex will mean the corresponding
parts of the proof of {\sl Theorem 1.1}.

\smallskip\noindent
{\sl Case 1$_b$)}. Without loss of generality we may assume that $V_\alpha$
is a domain with a fixed complex structure and a fixed antiholomorphic
involution $\tau$, and that $u_{\alpha, n}: V_\alpha \to X$ is a sequence of
$\tau$-invariant maps which are  (anti)holomorphic outside the set of $\tau
$-invariant points of $V_\alpha$. If we have the convergence of type {\sl C)}
 there is nothing to do. Otherwise we fix a $\tau$-invariant metric on
$V_\alpha$ compatible with the complex structure. Repeating the constructions
from {\sl Case 1)} we distinguish the ``bubbling'' points
$y^*_1, \ldots, y^*_l$ where the strong convergence fails.

Take the first such point $y^*_1$. Suppose $y^*_1$ is disjoint from
$\d\Sigma$, Then we may assume that $y^*_1\in V^+_\alpha$. Thus we can repeat
the rest of the constructions from {\sl Case 1)}. The only correction needed
at this place is that the neighborhood $\Delta( y^*_1, \varrho)$ of $y^*_1$
must be small enough and lie in $V^+_\alpha$. Transfering all these
constructions into $\tau(V_\alpha)$, we realize the inductive step preserving
$\tau$-invariance.

It remains to consider the case when $y^*_1 \in \d\Sigma$. This means that
$y^*_1$ is $\tau$-invariant. Let $z$ be a holomorphic coordinate in a
neigborhood of $y^*_1$ on $V_\alpha$, such that $z=0$ in $y^*_1$, the
involution $\tau$ corresponds to the conjugation $z \mapsto \bar z$, and $\im
z >0$ in $\Sigma$. Find the sequences $r_n\lrar0$ of radii and $x_n \to
y^*_1$ using the constructions from {\sl Case 1)}. Note that the sequence
$\tau(x_n)$ have the same property. Thus, replacing some points $x_n$ by
$\tau(x_n)$, we may additionally assume that all $x_n$ lie in $\barr
V^+_\alpha$. Let $v_n:\Delta (0,{\varrho \over 2r_n}) \to (X,J_n)$ be the
rescalings of maps $u_n$ defined by $v_n(z) \deff u_n(x_n + {z\over r_n})$.
Argumentations of {\sl Case 1} shows that there exists the limit $v_\infty:
\cc \to X$ of (a subsequence of) $\{ v_n \}$ which extends to a map $v_\infty
: S^2 \to X$.

Denote by $\rho_n$ the distance from $x_n$ to $\d\Sigma$ and by $\ti x_n$ the
point on $\d\Sigma$ closest to $x_n$. Then $x_n = \ti x_n + \isl \rho_n$ in the
coordinate $z$ introduced above. Besides, $\lim \ti x_n = y^*_1$. We consider
2 subcases according to possible behavior of $\rho_n$ and $r_n$.

\smallskip\noindent
{\sl Subcase 1$'_b$): $\bigl\{{ \rho_n \over r_n }\bigr\}$ is bounded}.
Passing to a subsequence, we may assume that ${\rho_n \over r_n}$ converges.
Fix an upper bound $b$ for the sequence ${\rho_n \over r_n}$. In particular,
$b \ge \lim {\rho_n \over r_n}$.

For $n>\!>1$ define maps $v_n:\Delta (0,{\varrho \over 2r_n}-b) \to (X,J_n)$
and $\ti v_n:\Delta (0,{\varrho \over 2r_n}-b) \to (X,J_n)$ setting $v_n(z)
\deff u_n (x_n +r_n z)$ and $\ti v_n(z)\deff u_n(\ti x_n +r_n z)$
respectively. Then every $\ti v_n$ is the shift of the map $v_n$ by
$\isl{\rho_n\over r_n}$, \ie $\ti v_n(z)= v_n\bigl(z + \isl {\rho_n\over
r_n}\bigr)$. The arguments of {\sl Case 1)} show that $v_n$ converge on
compact subsets of $\cc$ to a non-constant map. Consequently, $\ti v_n$ also
converge on compact subsets of $\cc$ to a non-constant map $\ti v_\infty: \cc
\to X$. Moreover, since $\area( \ti v_\infty( \cc))$ is finite, $\ti v_\infty$
extends to a map $\ti v_\infty: S^2 \to X$. By the choice of $\eps^b$,
$\area(\ti v_\infty (S^2)) \ge 3\eps^b$. Changing the choice of the constant
$b$, we can additionally assume that $\area(\ti v_\infty( \Delta(0,b)) \ge
2\eps^b$. Then for all sufficiently big $n$ we get
$$
\area(\ti v_n(\Delta(0,b))  \ge \eps^b
\eqno(5.?10)
$$

For $n>\!>1$ we define the coverings of $V_\alpha$ by 3 sets
$$
V^{(n)}_{\alpha,1} \deff  V_\alpha
\bs \barr\Delta (0, {\textstyle{\varrho\over2}}),
\qquad
V^{(n)}_{\alpha,2} \deff \Delta (0, \varrho)
\bs \barr \Delta (\ti x_n, br_n),
\qquad
V^{(n)}_{\alpha,3} \deff \Delta (\ti x_n, 2br_n).
$$
Fix $n_0$ sufficiently big. Denote $V_{\alpha,1} \deff V^{(n_0)}_{\alpha,1}$,
$V_{\alpha,2} \deff V^{(n_0)}_{\alpha,2}$, and $V_{\alpha,3} \deff V^{(n_0)}
_{\alpha,3}$. There exist diffeomorphisms $\psi_n: V_1 \to V_1$ such that
$\psi_n: V_{\alpha,1} \to V^{(n)}_{\alpha,1}$ is identity, $\psi_n:
V_{\alpha, 2} \to V^{(n)}_{\alpha,2}$ is a diffeomorphism, and $\psi_n:
V_{\alpha,3} \to V^{(n)}_{\alpha,3}$ is biholomorphic \wrt the complex
structures, induced from $C_n$ by means of $\sigma^d_n$. Note that the sets
$V^{(n)}_{\alpha,i}$ are $\tau$-invariant. Moreover, we can choose the maps
$\psi_n$ in such a way that $\psi_n$ are also $\tau$-invariant.

The covering $\{ V_{\alpha,1}, V_{\alpha,2}, V_{\alpha,3} \}$ of $V_1$ and
parametrizations $\ti\sigma_n \deff \sigma_{\alpha,n} \scirc \psi_n: V_1 \to
C_{\alpha, n}$ satisfy the conditions of {\sl Lemma 5.11}. Moreover,
inequality (5.?1) implies $\area(u_n(\ti\sigma_n(V_{\alpha,i}))) \le
(N-1)\eps^b$. Consequently, we can apply the inductive assumptions for the
sequence of curves $\ti\sigma_n(V_{\alpha,i})$ and finish the proof by
induction.

\smallskip\noindent
{\sl Subcase 1$''_b$): $\bigl\{{ \rho_n \over r_n}\bigr\}$ is unbounded}.
Passing to a subsequence, we may assume that ${\rho_n \over r_n}$ increases
infinitely. However, $\rho_n \lrar 0$ since $x_n \lrar y^*_1 \in \d\Sigma$.

Define maps $v_n:\Delta (0,{\varrho \over 2r_n}) \to (X,J_n)$ setting $v_n(z)
\deff u_n(x_n +r_n z)$. As in {\sl Case~1)}, $v_n$ converge on compact
subsets of $\cc$ to a non-constant map $v_\infty: \cc \to X$, which extends
to a map from the whole sphere $S^2$. Choose $b>0$ satisfying (5.?1).

For $n>\!>1$ we define the coverings of $V_\alpha$ by 5 sets
$$\mathsurround=0pt
\matrix\format\l\ \ &\l\\
V^{(n)}_{\alpha,1} \deff  V_\alpha
\bs \barr\Delta (0, {\textstyle{\varrho\over2}}),
&
V^{(n)}_{\alpha,2} \deff \Delta(0, \varrho) \bs
\barr \Delta(0, 2\rho_n)
\cr
\noalign{\vskip5pt}
\rlap{$
V^{(n)}_{\alpha,3} \deff \Delta(0, 4\rho_n) \bs
\bigl(\barr \Delta (x_n, br_n) \cap  \barr \Delta (\tau(x_n), br_n) \bigr)
$}&
\cr
\noalign{\vskip5pt}
V^{(n)}_{\alpha,4} \deff \Delta (x_n, 2br_n),
&
V^{(n)}_{\alpha,5} \deff \Delta(\tau(x_n), 2br_n).
\endmatrix
$$
Fix $n_0$ sufficiently big. Denote $V_{\alpha,i} \deff V^{(n_0)}_{\alpha,i}$,
$i=1,\ldots,5$. Then for every $n>\nobreak\!>1$ there exists a diffeomorphism
$\psi_n: V_1 \to V_1$ with the following properties:

\sli $\psi_n$ maps $V_{\alpha,i}$ onto $V^{(n)}_{\alpha,i}$ diffeomorphically;

\slii $\psi_n: V^{(n)}_{\alpha,1} \to V^{(n)}_{\alpha,1}$ is the identity;

\sliii $\psi_n: V_{\alpha, 2} \to V^{(n)}_{\alpha,2}$ and $\psi_n:
V_{\alpha,3} \to V^{(n)}_{\alpha,3}$ are diffeomorphisms;

\sliv $\psi_n: V_{\alpha,3} \to V^{(n)}_{\alpha,3}$ is biholomorphic \wrt
the complex structures, induced from $C_n$ by means of $\sigma^d_n$; and,
finally

\slv $\psi_n$ are $\tau$-invariant: $\tau \scirc\psi_n =\psi_n \scirc\tau$.
\newline\noindent
Note that the last property is obtained due to the fact that the sets
$V^{(n)} _{\alpha, i}$ are $\tau$-invariant. The rest constructions are the
same as in {\sl Subcase 1$'_b$)}.

\medskip\noindent
{\sl Case 2$_b$)}. Consider the parametrizations $\sigma_n: V_\alpha \to
C_{\alpha,n}$. Without loss of generality we may assume that the complex
strictures $\sigma_n^* j_n\ogran_{V_\alpha}$ are constant near boundary $\d
V_\alpha$ and converge to some complex structure. If we have the convergence
of type {\sl C)}, \ie the strong convergence, there is nothing to do.
Otherwise there exists only a finite set of points $\{y^*_1,\ldots, y^*_l\}$
where the strong convergence fails. Changing the parametrizations $\sigma_n$,
we may additionally assume that the strictures $\sigma_n^* j_n \ogran
_{V_\alpha}$ are constant in the neighborhood of these points. Then we repeat
the argumentations of {\sl Case 1$_b$)}.

\medskip\noindent
{\sl Case 3$_b$)}. Fix identifications $C_{\alpha,n} \cong \cala_0$ such that
every $C_{\alpha,n}^+$ is mapped onto $\cala_0^+$ and such that the induced
parametrization maps $\sigma_{\alpha, n} : V_\alpha \to \cala_0$ are the same
for all $n$ and $\tau$-invariant. Fix the standard representation of $\cala_0$
as the union of two discs $\Delta'$ and $\Delta''$ with identification of the
centers $0\in \Delta'$ and $0\in \Delta''$ into the nodal point of $\cala_0$,
still denoted by $0$. Let $\Delta' (x, r)$ denote the subdisc of $\Delta'$
with the center $x$ and the radius $r$.

Denote by $u'_n :\Delta' \to X$ and $u''_n :\Delta'' \to X$ the corresponding
``components'' of the maps $u_{\alpha,n} : C_{\alpha,n} \to X$. Find the
common collection of bubbling points $y^*_i$ for both sequences of maps $u'_n
:\Delta' \to X$ and $u''_n :\Delta'' \to X$. If there are no bubbling points,
then we obtain the convergence type {\sl B)} and the proof can be finished by
induction. Otherwise consider the first such point $y^*_1$, which lies, say,
on $\Delta'$. If $y^*_1$ is distinct from the nodal point $0 \in \Delta'$,
then we simply repeat all the construction {\sl Case 1$_b$)}.

It remains to consider the case $y^*_1=0 \in \Delta'$. The following
modifications of the argumentations are needed. Repeat the construction of
the radii $r_n\lrar 0$ and the points $x_n \lrar y^*_1=0$ from {\sl Case
1$_b$}. Then $\{x_n \}$ is a sequence in the half-disk $\delta^{\prime+}
\deff \{ z\in \Delta' : \im z \ge0 \}$. Set $\ti x_n \deff \re(x_n)$, $\rho_n
\deff \im(x_n)$ and $R_n \deff |x_n|$. Thus $x_n= \ti x_n +\isl\rho_n$, $R_n$
is the distance from $x_n$ to the point $0=y^*_1 \in \Delta'$, whereas
$\rho_n$ is the distance from $x_n$ to the interval $]-1,1[ \subset \Delta'$,
the set $\tau$-invariant points of $\Delta'$. Thus $\rho_n \le R_n$. Fix
$\varrho>0$ such that the disc $\Delta'(0,\varrho)$ contains no bubbling points
$y^*_i\not=0 \in \Delta'$.

Depending on the behavior of the sequences $r_n$, $\rho_n$ and $R_n$, we
consider the 4 subcases.

\medskip\noindent
{\sl Subcase 3\/$'_b$): The sequence $\bigl\{{ R_n \over r_n }\bigr\}$ is
bounded.} Then the sequences $\bigl\{{ {\rho_n\over r_n} }\bigr\}$ and $\bigl
\{ {\ti x_n\over r_n} \bigr\}$ are also bounded. Passing to a subsequence we
may assume that the corresponding limits exist. Let $b$ be some upper bound
for the sequence $\bigl\{{ R_n \over r_n }\bigr\}$. Consider the maps $\ti
v_n: \Delta(0, {\varrho \over 2r_n}-b) \to X$ defined by $\ti v_n(z) \deff
u_n(\ti x_n+ {z\over r_n})$. Then $\ti v_n$ are $\tau$-invariant, $\ti v_n
\scirc \tau= \ti v_n$, $\ti v_n$ converge to a nonconstant map $\ti v_\infty:
\cc \to X$ on compact subsets of $\cc$, and $\ti v_\infty$ extends to a map
$\ti v_\infty: S^2 \to X$.

Since  $\ti v_\infty$ is nonconstant, $\norm{ d\ti v _\infty} ^2 _{L^2( S^2)}
= \area( \ti v_\infty( S^2) )\ge 3\eps_b$. Choose $b>0$ in such a way that
$$
\norm{d\ti v_\infty}_{L^2(\Delta (0,b))}^2 \ge 2\eps_b
\eqno(5.?15)
$$
and $b \ge 2\lim {R_n \over r_n}+2$. Due to {\sl Corollary 5.2} for $n>\!>1$
we obtain the estimate
$$
\norm{du'_n}_{L^2(\Delta' (\ti x_n, br_n))}^2 =
\norm{ d\ti v_n}_{L^2(\Delta (0,b))}^2 \ge \eps_b.
\eqno(5.?16)
$$
Note that $0\in \Delta' (\ti x_n, (b-1)r_n))$ for $n>\!>1$ by the choice of
$b$.

Define the coverings of $\cala_0$ by 4 sets
$$\mathsurround=0pt
\matrix\format\l\ \ &\l\\
W^{(n)}_1 \deff  \Delta'
\bs \barr\Delta' (0, {\textstyle{\varrho\over2}}),
&
W^{(n)}_2 \deff \Delta' (0, \varrho)
\bs \barr \Delta' (\ti x_n, br_n),
\cr
\noalign{\vskip4.5pt}
W^{(n)}_3 \deff \Delta' (\ti x_n, 2br_n)
\bs \barr\Delta'(0,{\textstyle {r_n \over 2}}),
&
W^{(n)}_4 \deff \Delta' (0, r_n) \cup
\Delta'',
\endmatrix
$$
and lift them to $V_\alpha$ by putting $V^{(n)}_{\alpha,i} \deff \sigma\inv
_{\alpha,n}(W^{(n)}_i)$. Choose $n_0 >\!> 0$, such that $|x_n| < (b-1)r_n$
and the relation (5.?16) holds for all $n \ge n_0$. Set $V_{\alpha,i} \deff
V^{( n_0)} _{\alpha, i}$. Fix diffeomorphisms $\psi_n: V_\alpha \to V_\alpha$
such that $\psi_n: V_{\alpha,1} \to V^{(n)} _{\alpha,1}$ is the identity map,
$\psi_n: V_{\alpha,2} \to V^{(n)}_{\alpha,2}$ and $\psi_n: V_{\alpha,3} \to
V^{(n)} _{\alpha,3}$ are diffeomorphisms, and $\psi_n: V_{\alpha,4} \to
V^{(n)} _{\alpha,4}$ correspond to isomorphisms of nodes $W^{(n)}_4 \cong
\cala_0$. Set $\sigma'_n \deff \sigma_n \scirc \psi_n$. The choice above can
be done in such a way that the refined covering $\{ V_{\alpha,i} \}$ of
$V_\alpha$ and parametrization maps $\sigma'_n: V_\alpha \to C_{\alpha,n}$
have the properties of {\sl Lemma 5.11}. Relation (5.?16) implies the
estimate $\area(u_n (\sigma'_n(V_{\alpha,i})) \le (N-1)\,\eps$. This provides
the inductive conclusion for {\sl Subcase 3\/$'_b$)}.

\medskip\noindent
{\sl Subcase 3\/$''_b$): The sequence $\bigl\{{R_n \over r_n }\bigr\}$
increases infinitely but $\bigl\{{\rho_n \over r_n }\bigr\}$ remains
bounded}.  Note that in this subcase we still have the relation $R_n \lrar
0$, or equivalently, $x_n \lrar 0$. On the other hand, $\lim {\rho_n \over
R_n} =0$.  This implies that for $\ti R_n \deff |\ti x_n|$ we have $\lim
{\ti R_n \over R_n} =1$ since $R_n^2 = \ti R_n^2 + \rho_n^2$.

We proceed as follows. Define the maps $\ti v_n: \Delta(0, {\varrho \over
2r_n}-b) \to X$ setting $\ti v_n(z) \deff u'_n(\ti x_n + {z\over r_n})$. Then
$\ti v_n$ have the same properties as in {\sl Subcase 3\/$'_b$)}. Choose $b>0$
obeying the relation (5.?15). Then for $n>\!>0$ we get the property (5.?16).

For $n>\!>0$ define the coverings of $\cala_0$ by 6 sets
$$\mathsurround=0pt
\matrix\format\l\ \ &\l\\
W^{(n)}_1 \deff  \Delta'
\bs \barr\Delta' (0, {\textstyle{\varrho\over2}}),
&
W^{(n)}_2 \deff \Delta' (0, \varrho)
\bs \barr \Delta' (\ti x_n, 2\ti R_n),
\cr\noalign{\vskip5pt}
W^{(n)}_3 \deff  \Delta' (\ti x_n, 4\ti R_n)
\bs \bigl(\barr \Delta' (\ti x_n, {\ti R_n \over 6})
\cup \barr \Delta' (0, {\ti R_n \over 6}) \bigr)
&
W^{(n)}_4 \deff \Delta' (0, {\ti R_n \over 3}) \cup \Delta'',
\cr\noalign{\vskip5pt}
W^{(n)}_5 \deff \Delta' (\ti x_n, {\textstyle{\ti R_n \over 3}})
\bs \barr\Delta'(\ti x_n, br_n),
&
W^{(n)}_6 \deff \Delta' (0, 2br_n ),
\endmatrix
$$
and lift them to $V_\alpha$ by putting $V^{(n)}_{\alpha,i} \deff \sigma\inv
_{\alpha,n}(W^{(n)}_i)$. Choose $n_0 >\!> 0$, such that $R_{n_0} >\!> b
r_{n_0} $, and set $V_{\alpha,i} \deff V^{(n_0)} _{\alpha, i}$. Choose
diffeomorphisms $\psi_n: V_\alpha \to V_\alpha$ such that $\psi_n: V
_{\alpha, 1} \to V^{(n)} _{\alpha,1}$ is the identity map, $\psi_n: V
_{\alpha, 2} \to V^{(n)}_{\alpha,2}$, $\psi_n: V_{\alpha,4} \to V^{(n)}
_{\alpha,4}$ and $\psi_n: V_{\alpha,5} \to V^{(n)}_{\alpha,5}$ are
diffeomorphisms, and finally, $\psi_n: V_{\alpha,6} \to V^{(n)}_{\alpha,6}$
corresponds to isomorphisms of nodes $W^{(n)}_6 \cong \cala_0$. Set
$\sigma'_n \deff \sigma_n \scirc \psi_n$. Note that the choices can be done
in such a way that $\{ V_{\alpha,i} \}$ and parametrization maps $\sigma'_n:
V_\alpha \to C_{\alpha,n}$ have the properties of {\sl Lemma 5.11}. As above,
we get the estimate $\area(u_n (\sigma'_n(V_{\alpha,i} )) \le (N-1)\, \eps$
due to (5.?16). Thus we get the inductive conclusion for {\sl Subcase
3\/$''_b$)} and can proceed further.

\smallskip\noindent
{\sl Subcase 3\/$'''_b$): The sequence $\bigl\{{\rho_n \over r_n} \bigr\}$
increases infinitely, but $\bigl\{{ R_n \over \rho_n} \bigr\}$ remains
bounded}. Then $\bigl\{{R_n \over r_n} \bigr\}$ also increases infinitely, but
both sequences $\{R_n \}$ and $\bigl\{{ \rho_n }\bigr\}$ converge to 0.
We may also assume
that $\bigl\{{ \rho_n \over R_n }\bigr\}$ and $\bigl\{{ \ti x_n \over R_n }
\bigr\}$ also converge. Set $a_1 \deff \lim {\ti x_n \over R_n}$, $a_2 \deff
\lim{ \rho_n \over R_n}$, $a\deff a_1 + \isl a_2$, and $\barr a \deff a_1
- \isl a_2$. Note that $0 < a_2 \le 1$ and that the involutions $\tau_n$ in
$C_{\alpha, n}$ correspond to the complex conjugation $z\to \barr z$ in
$\Delta'$. In particular, $\barr x_n = \tau_n(x_n)$.

Consider maps $v_n: \Delta(0, {\varrho \over 2r_n}) \to X$ defined by $v_n(z)
\deff u'_n(x_n + {z\over r_n})$. Then the sequence $\{ v_n \}$ converges on
compact subsets to a nonconstant map which extends to the map $v_\infty :S^2
\to X$. Moreover, we can fix sufficiently big $b>0$ such that for $n>\!>0$ we
get the property (5.?16).

For $n>\!>0$ define the coverings of $\cala_0$ by 8 sets
$$\mathsurround=0pt
\matrix\format\l\ \ \ \ \ &\l\\
W^{(n)}_1 \deff  \Delta'
\bs \barr\Delta' (0, {\textstyle{\varrho\over2}}),
&
W^{(n)}_2 \deff \Delta' (0, \varrho)
\bs \barr \Delta' (0, 2R_n),
\cr\noalign{\vskip6pt}
\rlap{$W^{(n)}_3 \deff  \Delta' (0, 4R_n)
\bs \bigl(\, \barr \Delta' ( a\, R_n , {a_2 R_n \over 4})
\cup \barr \Delta' ( \barr a\, R_n , { a_2 R_n \over 4})
\cup \barr \Delta' (0, {a_2 R_n \over 4})
\,\bigr)\qquad$}
&
\cr\noalign{\vskip6pt}
W^{(n)}_4 \deff \Delta' (0, {a_2 R_n \over 3}) \cup \Delta''.
\cr\noalign{\vskip5pt}
W^{(n)}_5 \deff \Delta' (a R_n, {a_2 R_n \over 3})
\bs \barr\Delta'(x_n, br_n),
&
W^{(n)}_6 \deff \Delta' (x_n, 2br_n ),
\cr\noalign{\vskip5pt}
W^{(n)}_7 \deff \Delta' (\barr a R_n, {a_2 R_n \over 3})
\bs \barr\Delta'(\barr x_n, br_n),
&
W^{(n)}_8 \deff \Delta' (\barr x_n, 2br_n ),
\endmatrix
$$
and lift them to $V_\alpha$ by putting $V^{(n)}_{\alpha,i} \deff \sigma\inv
_{\alpha,n}(W^{(n)}_i)$. Fix sufficiently big $n_0 >\!> 0$, and set
$V_{\alpha, i} \deff V^{(n_0)} _{\alpha, i}$. Choose diffeomorphisms $\psi_n:
V_\alpha \to V_\alpha$ mapping $V _{\alpha, i}$ diffeomorphically onto
$V^{(n)} _{\alpha, i}$ such that the assertions of {\sl Lemma 5.11} are
fulfilled. As above, we get the estimate $\area(u_n (\sigma'_n(V_{\alpha,i} ))
\le (N-1)\, \eps$. This gives the inductive conclusion for {\sl Subcase
3\/$'''_b$)}.

\smallskip\noindent
{\sl Subcase 3\/$''''_b$): The sequences $\bigl\{{\rho_n \over r_n} \bigr\}$
and  $\bigl\{{ R_n \over \rho_n} \bigr\}$ increases infinitely}.
Thus $\lim {\ti R_n \over R_n} =1$.
We consider the sequence of maps $\{ v_n \}$. It is defined in the same way
as in the previous subcase and has the same properties. In particular,
$\{ v_n \}$ converges to the map $v_\infty :S^2 \to X$ and there exists
a sufficiently big $b>0$ such that for $n>\!>0$ we get the property (5.?16).

For $n>\!>0$ define the coverings of $\cala_0$ by 10 sets
$$\mathsurround=0pt
\matrix\format\l\ \ \ \ \ &\l\\
W^{(n)}_1 \deff  \Delta'
\bs \barr\Delta' (0, {\textstyle{\varrho\over2}}),
&
W^{(n)}_2 \deff \Delta' (0, \varrho)
\bs \barr \Delta' (0, 2R_n),
\cr\noalign{\vskip6pt}
\rlap{$W^{(n)}_3 \deff  \Delta' (0, 4R_n)
\bs \bigl(\, \barr \Delta' ( 0, {\ti R_n \over 4})
\cup \barr \Delta' ( \ti x_n, {\ti R_n \over 4})
\bigr),$}
&
\cr\noalign{\vskip6pt}
W^{(n)}_4 \deff \Delta' (0, {\ti R_n \over 3}) \cup \Delta'',
&
W^{(n)}_5 \deff \Delta' (\ti x_n, {\ti R_n \over 3})
\bs \barr\Delta'(\ti x_n, 2\rho_n),
\cr\noalign{\vskip6pt}
\rlap{$W^{(n)}_6 \deff \Delta' (\ti x_n, 4\rho_n )
\bs \bigl(\, \barr \Delta' ( x_n, {\rho_n \over 4})
\cup \barr \Delta' ( \barr x_n, {\rho_n \over 4})
\bigr), $}
\cr\noalign{\vskip5pt}
W^{(n)}_7 \deff \Delta' (x_n, {\rho_n \over 3})
\bs \barr\Delta'(x_n, br_n),
&
W^{(n)}_8 \deff \Delta' (x_n, 2br_n ),
\cr\noalign{\vskip5pt}
W^{(n)}_9 \deff \Delta' (\barr x_n, {\rho_n \over 3})
\bs \barr\Delta'(\barr x_n, br_n),
&
W^{(n)}_{10} \deff \Delta' (\barr x_n, 2br_n ).
\endmatrix
$$
The rest ``manipulations'' with $W^{(n)}_i$ are the same as in the previous
subcases. As result, we obtain the covering of $V_\alpha$ by sets $V_{\alpha,
i} \deff \sigma\inv _{\alpha, n_0}(W^{(n_0)}_i)$ with an appropriate $n_0 >\!>
0$ and refined parametrizations $\sigma'_n: V_\alpha \to C_{\alpha, n}$, for
which the assertions of {\sl Lemma 5.11} are fulfilled.
As above, we get the estimate $\area(u_n (\sigma'_n(V_{\alpha,i} ))
\le (N-1)\, \eps$. This gives the inductive conclusion for {\sl Subcase
3\/$''''_b$)}.

\smallskip\noindent
{\sl Case 4$_b$): $V_\alpha$ is a cylinder, such that conformal radii
of $(V_\alpha, \sigma_n^*j_n)$ increase infinitely}. We can simply repeat the
contructions made in {\sl Case 4)} from the proof of {\sl Theorem 1.1}.
An additional attention is  needed to preserve $\tau$-invariantness.

\smallskip
The proof of theorem  can be now finished by induction.
\qed

\baselineskip=13pt 
\medskip
\state Remark. Here we give some explanation of the geometrical meaning of
the constructions of the proof of {\sl Theorem 5.9} and describe the picture
of the bubbling. We restrict ourselves to {\sl Case 3$_b$)} as the most
complicated one, the constructions of other cases can be treated similarly.
The reflection principle allows us to reduce the {\sl Case 3$_b$)} to
consideration of $\tau$-invariant maps $u^d_n: \cala_0 \to X$ from the standard
node which are $J_n$-holomorphic on $\cala^+_0$. The situation is most
different from the situations of {\sl Theorem 1.1} when the bubbling appears
in the nodal point. In this case we must take into consideration not only
parameters $r_n$ describing the size of energy localization of the bubbled
sphere, but also additional parameters $R_n$ and $\rho_n$. These ones describe
the position of the localization centers $x_n$ \wrt the nodal point and the
set of $\tau$-invariant points of $\cala_0$. Depending on the behavior of $r_n$
$\rho_n$, and $R_n$ we can have 4 different types of the bubbling and
corrsponding {\sl Subcases 3\/$'_b$)--3\/$''''_b$)}.

\smallskip
\line{%
\vtop{\hsize=.52\hsize
In {\sl Subcase 3\/$'_b$)} the bubbling take place in the nodal point, so that
the nodal point remains on the bubbled sphere (see region $W^{(n)}_3$ on
Fig.~9). Furthermore, the bubbled sphere contains another one nodal point.
This one  appears in the limit of long cylinders $W^{(n)}_2$. Note $W^{(n)}_2$
can either strongly converge to a boundary node, or have additional bubblings.
}
\hfil
\vtop{\hsize=.43\hsize\xsize=\hsize\nolineskip\rm
\putm[-.03][.17]{W^{(n)}_1}%
\putm[.23][.09]{W^{(n)}_2}%
\putm[.53][.17]{W^{(n)}_3}%
\putm[.85][.14]{W^{(n)}_4}%
\noindent
\epsfxsize=\xsize\epsfbox{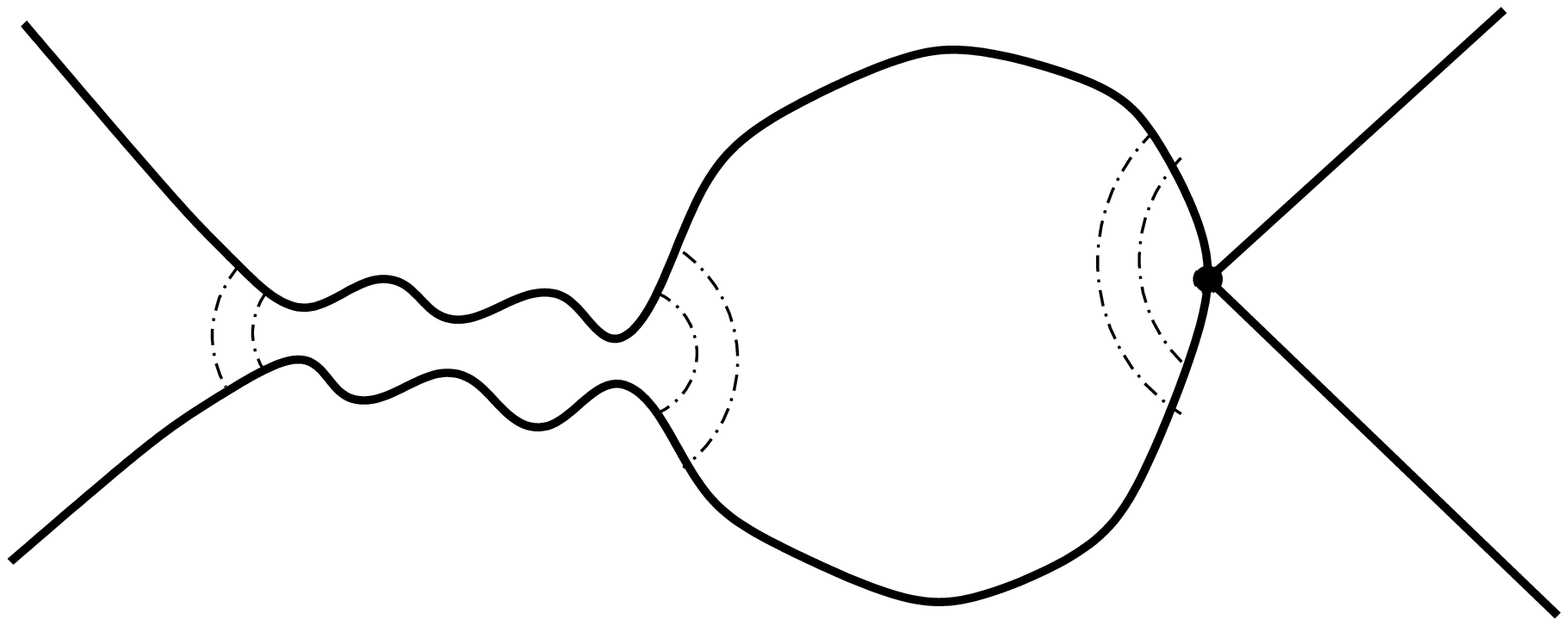}
\smallskip\smallskip
\putt[.0][.0]{
\centerline{Fig.~9. Bubbling in {\sl Subcase 3\/$'_b$)}.}
}
}}

\smallskip
Turning back from ``doubled'' description by $\tau$-invariant objects to the
original  maps $u_n: \cala^+_0 \to X$ with totally real boundary condition we
obtain the following picture. Since every covering piece $W^{(n)}_i$ is
$\tau$-invariant, for $\cala_0^+$ we get the covering piece $W^{(n)\,+}_i
\deff W^{(n)}_i \cap \cala_0^+$. Thus we obtain a bubbled dics represented by
$W^{(n)\,+}_3$ instead of the bubbled sphere represented by $W^{(n)}_3$,
the sequence of long strips $W^{(n)\,+}_2$ instead of the  sequence of long
cylinder $W^{(n)\,+}_2$ and so on.

\medskip
\line{%
\vtop{\hsize=.533\hsize \baselineskip=13.2pt
In {\sl Subcase 3\/$''_b$)} the bubbling happens at the boundary but away
from the nodal point. In the limit we obtain 2 bubbled spheres. The first one
is the limit of the re\-scaled maps $v_n$ (region $W^{(n)}_6$ on the Fig.~10).
The appearance of the second sphere can be explained as follows. The part of
the node $\cala_0$ between the first bubbled sphere $W^{(n)}_6$ and the
``constant part'' $W^{(n)}_1$ of the node is a long cylinder, represented by
pieces $W^{(n)}_2$, $W^{(n)}_3$, and $W^{(n)}_5$.
}
\hfil \vtop{\hsize=.44\hsize\xsize=\hsize\nolineskip\rm
\putm[-.02][.40]{W^{(n)}_1}%
\putm[.21][.31]{W^{(n)}_2}%
\putm[.53][.43]{W^{(n)}_3}%
\putm[.50][.06]{W^{(n)}_6}%
\putm[.63][.22]{W^{(n)}_5}%
\putm[.85][.415]{W^{(n)}_4}%
\noindent
\epsfxsize=\xsize\epsfbox{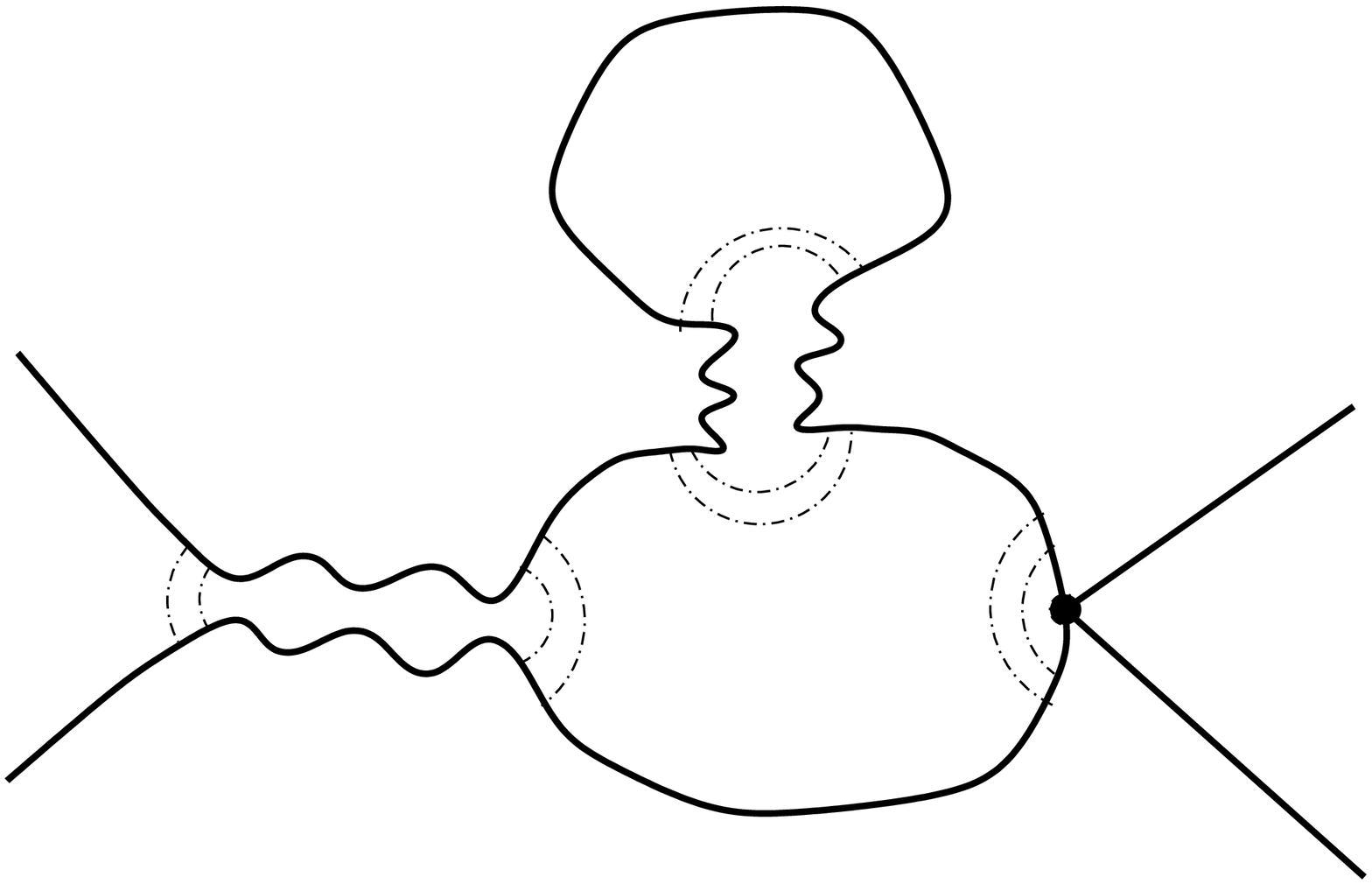} 
\smallskip\smallskip 
\putt[.0][.0]{
\centerline{Fig.~10. Bubbling in {\sl Subcase 3\/$''_b$)}.}
}
}}

Howeverer, because of the presence of the nodal point (piece $W^{(n)}_4$ on
the figure), this ``part inbetween'' is topologically not a cylinder (\ie an
annulus) but pants. Furthermore, the complex structures on the pants are not
constant. To get pants with constant structure (piece $W^{(n)}_3$) we cut off
the annuli $W^{(n)}_2$ and $W^{(n)}_5$. Since $\lim R_n =0 = \lim {r_n\over
R_n}$, the conformal radii of these annuli increase infinitely. This shows
that $W^{(n)}_2$ and $W^{(n)}_5$ are sequences of long cylinders and that the
sequence $W^{(n)}_3$ defines in the limit a sphere with 3 nodal points.

As in {\sl Subcase 3\/$'_b$)} every covering piece $W^{(n)}_i$ is $\tau
$-invariant, whereas $W^{(n)\,+}_i \deff W^{(n)}_i \cap \cala_0^+$ is the 
``half'' of $W^{(n)}_n$. Thus for sequence of undoubled maps $u_n: \cala^+_0 
\to X$ we get the following picture of the bubbling . The limit contains 2
bubbled discs represented by $W^{(n)\,+}_6$ and $W^{(n)\,+}_3$, a boundary node
$W^{(n)\,+}_4$, and possibly futher bubbled pieces which can appear in the
limit of long strips $W^{(n)\,+}_2$ and $W^{(n)\,+}_5$. Note also that the
action of the involution $\tau$ on the pants $W^{(n)}_3$ is described by
Fig.~8\.b).

\medskip
In {\sl Subcase 3\/$'''_b$)} the bubbling takes place near, but not at the
boundary. Indeed, since ${\rho_n \over r_n}\lrar \infty$, the bubbled sphere
which appears as the limit of the sequence $\{v_n \}$ is not $\tau$-invariant.
To see this phenomenon we note that for any fixed $b>0$ the covering pieces
$W^{(n)}_5= \Delta' (x_n,2br_n)$ representing sufficient big part of this
sphere lie in $\cala^+_0$ for $n>\!>0$. This implies that the sequence $v_n
\scirc \tau$ converges to another bubbled sphere, which is $\tau$-symmetric
to the first one and represented by $W^{(n)}_7$.

Another one bubbled sphere, represented by $W^{(n)}_3$, appears from
pants between the first two spheres and the disc $\Delta'$. Since $\{ {R_n
\over \rho_n}\}$ remains bounded, the original nodal point remains on
this latter sphere.

\line{%
\vtop{\hsize=.53\hsize\baselineskip=13pt plus 1.5pt
The corresponding bubbling picture for undoubled maps $u_n: \cala^+_0 \to X$
is shown on on Fig.~11. The boundary of $\cala^+_0$ is drawn by thick line.
We obtain the bubbled sphere represented by $W^{(n)}_6$, the sequence of long
cylinders $W^{(n)}_8\!$, the bubbled disc $W^{(n)+}_3\!\!$, and the sequence
of long strips $W^{(n)+}_2$. Note that both sequences of long cylinders
and long strips can yield  further bubblings in the limit.
}
\hfil
\vtop{\hsize=.44\hsize\xsize=\hsize\nolineskip\rm
\putm[-.02][.34]{W^{(n)+}_1}%
\putm[.21][.26]{W^{(n)+}_2}%
\putm[.51][.42]{W^{(n)+}_3}%
\putm[.84][.34]{W^{(n)+}_4}%
\putm[.77][.07]{W^{(n)}_5}%
\putm[.50][.065]{W^{(n)}_6}%
\noindent
\epsfxsize=\xsize\epsfbox{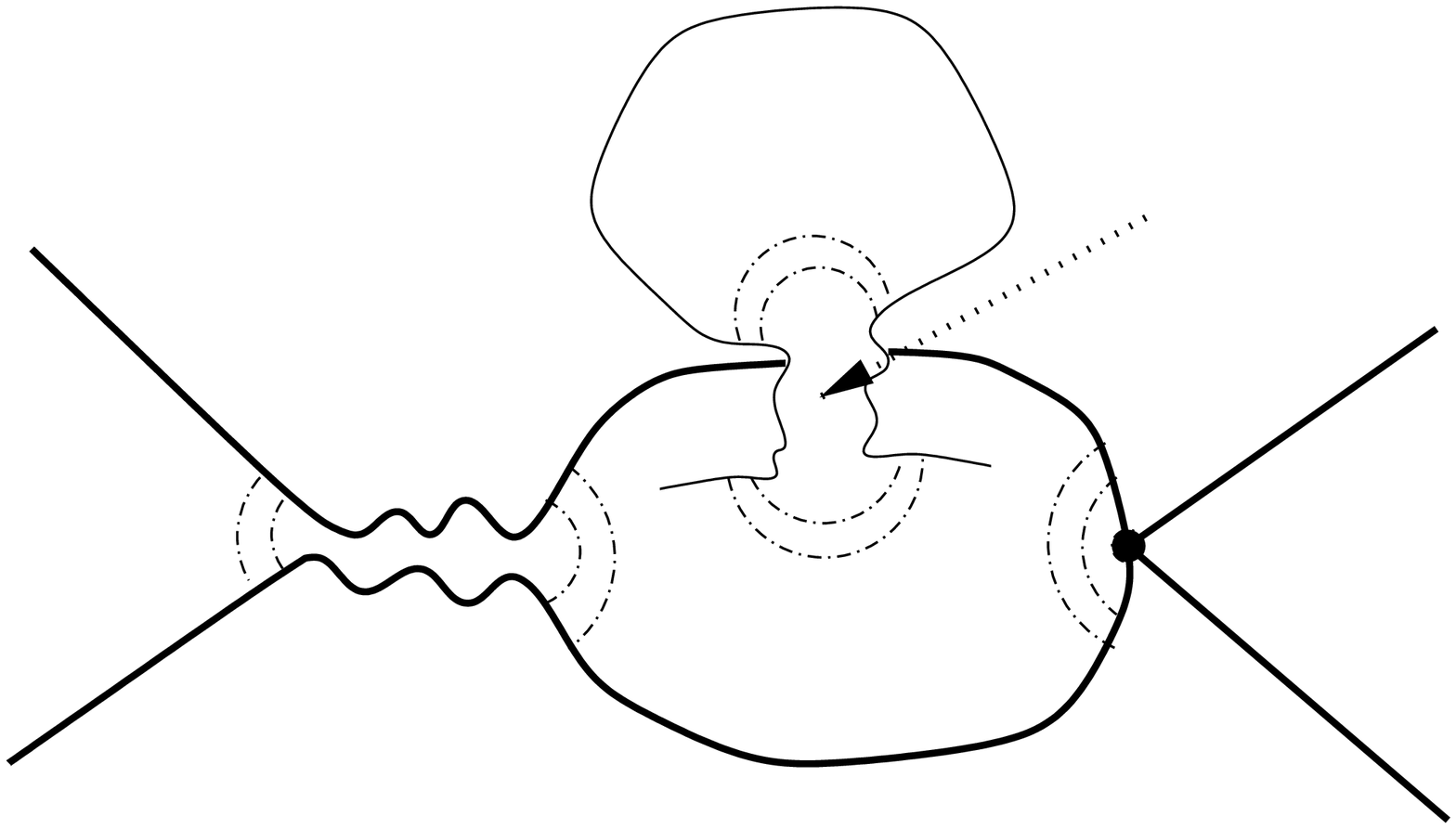}
\smallskip\smallskip
\putt[.0][.0]{
\centerline{Fig.~11. Bubbling in {\sl Subcase 3\/$'''_b$)}.}
}
\vskip\baselineskip\vskip11pt
}}

\medskip\smallskip
\line{%
\vtop{\hsize=.53\hsize\baselineskip=13pt 
The bubbling picture in {\sl Subcase 3\/$''''_b$)} is similar to the one of
previous subcase, so we explain only the difference. It comes from
the fact that the sequence $\{ {R_n \over \rho_n}\}$ is now unbounded, \ie
$\lim {\rho_n \over R_n}=0$. Informally speaking, this means that the long
cylinder from {\sl Subcase 3\/$'''_b$)} (piece $W^{(n)}_5$ on Fig.~11)
moves to the boundary of the bubbled disc (piece $W^{(n)+}_3$ on Fig.~11).
The precedure of additional re\-scaling divides every such dics into two
new discs connected by a strip, pieces $W^{(n)+}_3\!\!$, $W^{(n)+}_6\!\!$, and
$W^{(n)+}_5$ on Fig.~12 respectively. The infinite growth ${R_n \over
\rho_n} \to \infty$ means that $W^{(n)+}_5$ form a sequence of long strips.
}
\hfil
\vtop{\hsize=.44\hsize\xsize=\hsize\nolineskip\rm
\putm[-.02][.70]{W^{(n)+}_1}%
\putm[.22][.62]{W^{(n)+}_2}%
\putm[.51][.75]{W^{(n)+}_3}%
\putm[.86][.71]{W^{(n)+}_4}%
\putm[.82][.36]{W^{(n)+}_5}%
\putm[.51][.39]{W^{(n)+}_6}%
\putm[.80][.07]{W^{(n)}_7}%
\putm[.51][.065]{W^{(n)}_8}%
\noindent
\epsfxsize=\xsize\epsfbox{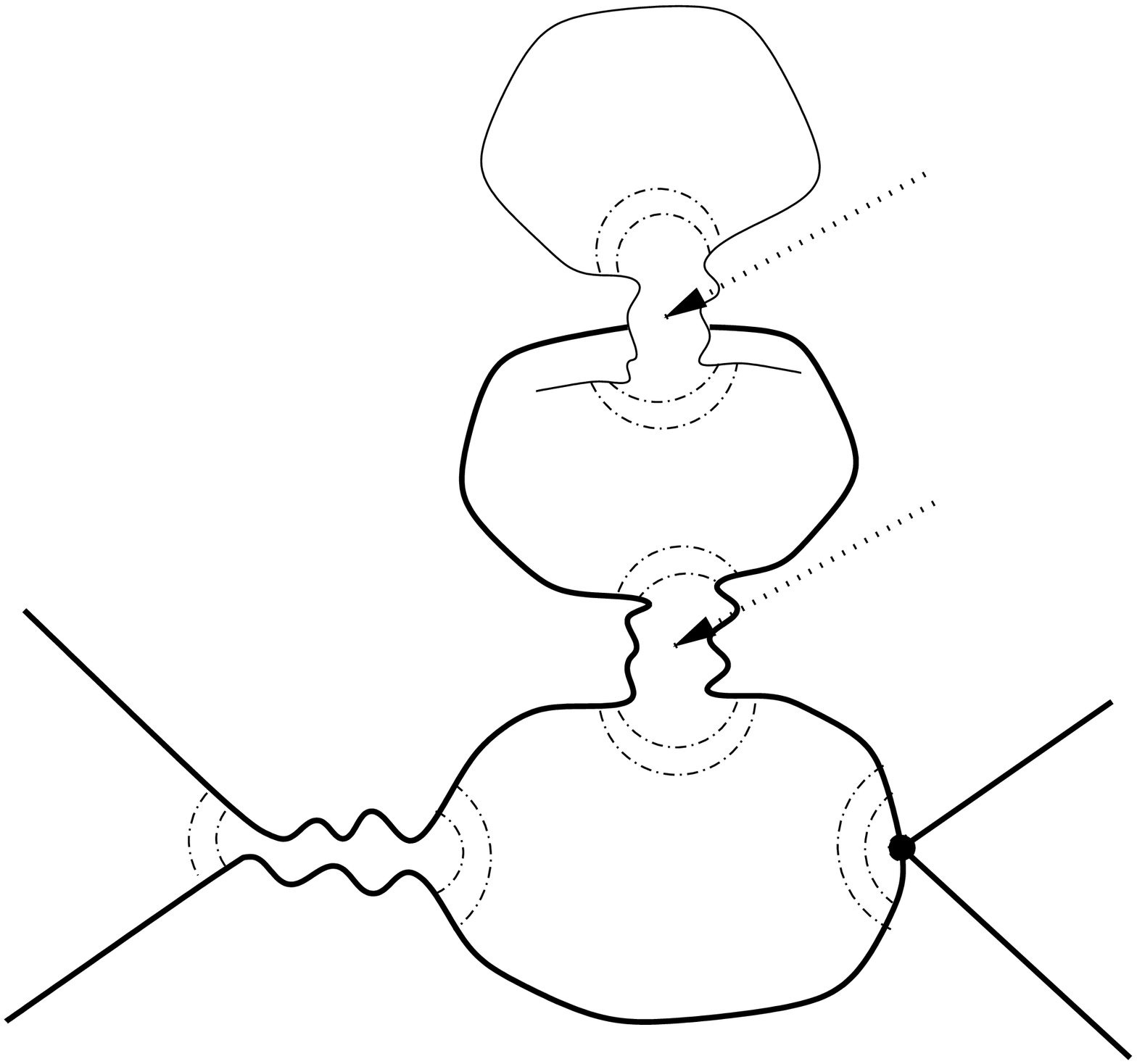}
\smallskip\smallskip
\putt[.0][.0]{
\centerline{Fig.~16. Bubbling in {\sl Subcase 3\/$''''_b$)}.}
}
\vskip\baselineskip\vskip11pt
}}

\spaceskip=4pt plus3.5pt minus 1.5pt
\font\csc=cmcsc10

\newdimen\length
\newdimen\lleftskip
\lleftskip=3.8\parindent
\length=\hsize \advance\length-\lleftskip
\def\entry#1#2#3#4\par{\parshape=2  0pt  \hsize%
\lleftskip \length%
\noindent\hbox to \lleftskip%
{\bf[#1]\hfill}
{\csc{#2}} 
{\sl{#3}} 
#4 \medskip
}
\ifx \twelvebf\undefined \font\twelvebf=cmbx12\fi

\bigskip\bigskip
\centerline{\bf R E F E R E N C E S}

\nobreak
\bigskip

\entry{Ab}{Abikoff W.:}{The Real Analytic Theory of Teichm\"uller Space.}
Springer-Verlag (1980).

\entry{D-G}{Dethloff G., Grauert H.:}{Deformation of compact Riemann surfaces
$Y$ of genus $p$ with distinguished points $P_1, \ldots P_m \in Y$.}
Int.\ Symp.\ ``Complex geometry and analysis'' in Pisa/Italy, 1988;
Lect.\ Notes in Math., {\bf1422}(1990), 37--44.

\entry{D-M}{Deligne P., Mumford D.:}{The irreducibility of the space of 
curves of a given genus.} IHES Math.\ Publ.,  {\bf36}(1969), 75--109. 

\entry{G}{Gromov M.:}{Pseudoholomorphic curves in symplectic manifolds}
Invent.\ Math., {\bf82}(1985), 307--347.

\entry{I-S}{Ivashkovich S., Shevchishin S.:}{Pseudoholomorphic curves and
envelopes of meromorphy of two-spheres in $\cc\pp^2$.} Preprint,
Sonderforschungsbereich 237 "Unordnung und grosse Fluktuationen",
(1995), available as e-print {\tt math.CV/9804014}.

\entry{K}{Kontsevich M.:}{Enumeration of rational curves via torus actions.}
Proc.\ Conf.\ ``The moduli spaces of curves'' on Texel Island, Netherland.
Birkh\"auser Prog.\ Math., {\bf129}(1995), 335--368.

\entry{K-M}{Kontsevich M., Manin Yu.:}{Gromov-Witten classes, quantum 
cohomology, and enumerative geometry.} Comm.\ Math.\ Phys., {\bf164}(1994),
525--562.

\entry{M}{Mumford D.:}{A remark on Mahler's compactness theorem.} Proc.\
Amer.\ Math.\ Soc., {\bf28}(1971), 289-294.

\entry{Pa}{Parker T.:}{Bubble tree convergence for harmonic maps} J.\ Diff.\
Geom., {\bf44}(1996), 595--633.

\entry{P-W}{Parker T., Wolfson J.:}{Pseudoholomorphic maps and bubble trees}
J.\ Geom.\ Anal., {\bf3}(1993), 63--98.

\entry{S-U}{Sacks J., Uhlenbeck K.:}{Existence of minimal immersions of
two-spheres} Annal.\ Math., {\bf113}(1981), 1--24.

\entry{S}{Sikorav J.-C.:}{Some properties of holomorphic curves in almost
complex manifolds.}In "Holomorphic curves in symplectic geometry".
Edited by M. Audin, J. Lafontaine. Birkh\"auser, Progress in
Mathematics v.\.117, Ch.V, 165-189.

\bigskip
\bigskip
\settabs 2 \columns
\+Institute of Applied Problems
&Institute of Applied problems\cr
\+of Mechanics and Mathematics
&of Mechanics and Mathematics\cr
\+Ukrainian Acad. Sci.,
&Ukrainian Acad. Sci.,\cr
\+vul. Naukova 3b, 290053 L'viv
&vul. Naukova 3b, 290053 L'viv \cr
\+Ukraine
&Ukraine \cr
\medskip
\+ U.F.R. de Maht\'ematiques
& Ruhr-Universit\"at Bochum \cr
\+ Universit\'e de Lille-I
& Mathematisches Institut \cr

\+ Villeneuve d'Ascq Cedex
& Universit\"atsstrasse 150 \cr
\+ 59655 France
& NA 4/67  44780 Bochum  Germany \cr
\+ ivachkov\@gat.univ-lille1.fr
& sewa\@cplx.ruhr-uni-bochum.de   \cr

\enddocument
\bye